\documentclass{report}
\usepackage[left=1.5in, right=1in, top=1in, bottom=1in]{geometry}
\usepackage{graphicx}
\usepackage{latexsym}
\usepackage{amsfonts}
\usepackage{amsthm}
\usepackage{amsmath}
\usepackage{amssymb}
\usepackage{copyright}

\newtheorem{thm}{Theorem}[section]
\newtheorem{prop}{Proposition}[section]
\newtheorem{lemma}{Lemma}[section]
\newtheorem{defn}{Definition}[section]

\newtheorem{algo}{Algorithm}[section]
\newtheorem{corollary}{Corollary}[thm]

\newcommand{\reals}{\mathbb{R}}
\newcommand{\complex}{\mathbb{C}}
\newcommand{\integers}{\mathbb{Z}}
\newcommand{\naturals}{\mathbb{N}}
\newcommand{\rationals}{\mathbb{Q}}

\newcommand{\skipline}{\vspace{12pt}}
\newcommand{\comment}[1]{}

\begin{document}
\pagenumbering{roman}
\pagestyle{plain}

% DOCUMENT INFO

\begin{center}
{\Large \bf ALGEBRAIC COMBINATORICS OF MAGIC SQUARES}\\
\skipline
By\\
\skipline
MAYA AHMED\\
B.Sc. (Bombay University) 1988\\
M.Sc. (IIT, Bombay) 1991\\
M.S. (University of Washington) 1997\\
\skipline

DISSERTATION\\
\skipline
Submitted in partial satisfaction of the requirements for the degree of\\
\skipline
DOCTOR OF PHILOSOPHY\\
\skipline
in\\
\skipline
MATHEMATICS\\
\skipline
in the\\
\skipline
OFFICE OF GRADUATE STUDIES\\
\skipline
of the\\
\skipline
UNIVERSITY OF CALIFORNIA,\\
\skipline
DAVIS\\
\vspace{1in}
Approved:\\
\skipline
{JES\'US DE LOERA}\\
\skipline
{ANNE SCHILLING}\\
\skipline
{CRAIG TRACY}\\
\skipline
Committee in Charge\\
\skipline
2004\\

\end{center}

% this starts the double-spacing
\renewcommand{\baselinestretch}{1.6}\small\normalsize

%%%%%%%%%%%%%%%%%%%%%%%%%%%%%%%%%%%%%%%%%%%%%%%%%%%%%%%%%%%%%%%%%%
%%%%%%%%%%%%%%%%%%%%%%%%%%%%%%%%%%%%%%%%%%%%%%%%%%%%%%%%%%%%%%%%%%
\newpage 
\begin{copyrightpage}
 \makecopyright
 \end{copyrightpage}
\newpage 

\begin{center}
{\LARGE Algebraic Combinatorics of Magic Squares}
\end{center}

\begin{center}
\underline{\bf \Large  Abstract} 
\end{center}
The problem of constructing magic squares is of classical interest and
the first known magic square was constructed around 2700 B.C. in
China.  Enumerating magic squares is a relatively new problem. In
1906, Macmahon enumerated magic squares of order 3.

In this thesis, we describe how to construct and enumerate magic
squares as lattice points inside polyhedral cones using techniques
from Algebraic Combinatorics. The main tools of our methods are the
Hilbert Poincar\'e series to enumerate lattice points and the Hilbert
bases to generate lattice points.  With these techniques, we derive
formulas for the number of magic squares of order 4. We extend
Halleck's work on Pandiagonal squares and Bona's work on magic cubes
and provide formulas for $5 \times 5$ pandiagonal squares and $3
\times 3 \times 3$ magic cubes and semi-magic cubes.

Benjamin Franklin constructed three famous squares which have several
interesting properties. Many people have tried to understand the
method Franklin used to construct his squares (called Franklin
squares) and many theories have been developed along these lines. Our
method is a new method to construct not only the three famous squares
but all other Franklin squares. We provide formulas for counting the
number of Franklin squares and also describe several symmetries of
Franklin squares.

Magic labelings of graphs are studied in great detail by Stanley and
Stewart.  Our methods enable us to construct and enumerate magic
labelings of graphs as well. We explore further the correspondence of
magic labelings of graphs and symmetric magic squares. We define
polytopes of magic labelings of graphs and digraphs, and give a
description of the faces of the Birkhoff polytope as polytopes of
magic labelings of digraphs.

%\medskip

%%%%%%%%%%%%%%%%%%%%%%%%%%%%%%%%%%%%%%%%%%%%%%%%%%%%%%%%%%%%%%%%%%
%%%%%%%%%%%%%%%%%%%%%%%%%%%%%%%%%%%%%%%%%%%%%%%%%%%%%%%%%%%%%%%%%%

\newpage
{\Large \bf ACKNOWLEDGMENTS} \\

I take this opportunity to thank my advisor Jes\'us De Loera for many
many things. It was his efforts that secured me admissions to UC Davis
and helped me pursue my doctorate degree. He began by teaching me to
ace exams and gradually taught me to solve seemingly impossible
research problems.  I thank him for being an excellent teacher who can
make difficult topics easy to understand. I thank him for introducing
me to a lot of wonderful mathematics especially algebraic
combinatorics. I also thank him for his exceptionally high standards
though it was daunting at times. The list is long, so I will just say
thanks for making me the mathematician I am today.

I also take this opportunity to thank all the teachers who have had a
substantial influence in my mathematical training. I thank Parvati
Subramaniam from S.I.E.S high school who was the first teacher I met
who could teach mathematics. She made mathematics easy for me and it
always stayed easy after-wards. I thank Edward Curtis from University
of Washington for making me an Algebra fanatic with the excellent
summer course he taught. I thank Ronald Irving at University of
Washington for launching my research career. He gave me the necessary
confidence with all his encouragement and enthusiasm. I thank Bernd
Sturmfels at UC Berkeley from whom I learned Computational Algebraic
Geometry and also the thrill of real modern mathematics. His influence
is present in all my work.

Its my great pleasure to thank the people at UC Davis for making my
stay here memorable. I thank Craig Tracy for all the excellent courses
he taught. Craig Tracy always found topics that felt like it was
tailor-made for each student in the class. I understand a lot of
combinatorics because of him. I thank Abigail Thompson and Motohico
Mulase for their timely encouragement. A special thanks to Joel Hass,
Matthew Franklin, and Anne Schilling for making my oral exam a
pleasant experience. I also thank Greg Kuperberg and Alexander
Shoshnikov for their help.  In fact, I thank everyone at UC Davis for
being so friendly.

I thank my coauthor Raymond Hemmecke from whom I have learned a lot of
good mathematics. A special thanks to my officemate Lipika Deka for
all the good times and my friend Ruchira Dutta for all her useful
comments on my thesis work. I thank my colleague Ruriko Yoshida for
help with the program LattE.

I thank my wonderful family: my father Neeliyara Devasia for being the
kindest soul on earth; my mother Mariakutty Devasia for all her hard
work; my elder brother Santosh Devasia for always being there for me;
my younger brother Vinod Devasia for all the good cheer; my
sisters-in-law Jessy and Blossom for being nice, my niece Lovita and
nephew Brian for just being fun; my in-laws, the Ahmed family: Shafi,
Sofia, Shabbir, Juzar, Cynthia, Jasmine, Alicia, and Mohammed for all
their help and good times; and finally, my husband Mohsin Ahmed for
being my best friend, my guardian angel, my greatest critique, and my
all time sponsor\footnote{Partially supported by NSF grant 0309694 and
0073815.}.

{\small
\begin{verse}
I can no other answer make but thanks, \\
And thanks; and ever thanks; and oft good turns \\
Are shuffled off with such uncurrent pay:\\
But, were my worth as is my conscience firm, \\
You should find better dealing. \\
-- William Shakespeare.
% twelfth night, Act 3, Scene 3.
\end{verse}
}

%%%%%%%%%%%%%%%%%%%%%%%%%%%%%%%%%%%%%%%%%%%%%%%%%%%%%%%%%%%%%%%%%%%%%%%%%%%%%%%%%%%%%%%%%%%%%%%%%%%%%%%%%%%%%%%%%%%%%%%%%%%%%%%%%%%%%%%%%%%%%%%%%%%%%%%%%%

\newpage
\large
\tableofcontents

% BEGIN TEXT OF THESIS
% here you start page-numbering in the upper right
                               
%%%%%%%%%%%%%%%%%%%%%%%%%%%%%%%%%%%%%%%%%%%%%%%%%%%%%%%%%%%%%%%%

\newpage
\pagestyle{myheadings} 
\pagenumbering{arabic}
\addcontentsline{toc}{chapter}{{\bf Foreword}}
\newpage
\large 
\chapter*{Foreword} \label{forewordchapter}
\thispagestyle{myheadings}
This thesis is a story of fusion of classical mathematics and modern
computational mathematics. We redefine the methods of solving
classical problems like constructing a magic square and exhibit the
power of computational Algebra. For example, only three Franklin
squares were known since such squares were difficult to construct. Our
results make it possible to construct any number of Franklin squares
easily. A detailed description of our methods is given in Section
\ref{method}. More significantly, our study leads to an elegant
description of the faces of the Birkhoff polytope as polytopes of
magic labelings of digraphs. The main theorems are introduced in
Chapter \ref{introductionchapter} and the subsequent chapters explore
the details of these results.

This thesis is based on my three papers
\begin{itemize}
\item Polyhedral cones of magic cubes and squares (joint work with
J. De Loera and R. Hemmecke) \cite{adh};

\item How many squares are there,
Mr. Franklin?: Constructing and Enumerating Franklin Squares \cite{ma};

\item Magic graphs and the faces of the
Birkhoff polytope \cite{maya}.
\end{itemize}

\newpage
\markright{  \rm \normalsize CHAPTER 1. \hspace{0.5cm}
 Introduction}
\large 
\chapter{Introduction} \label{introductionchapter}
\thispagestyle{myheadings}
\section{Magic Squares.} \label{magicsquares}
{\small 
\begin{verse}
\noindent Miracles seem very wonderful to the people who witness them, and very
simple to the people who perform them. That does not matter: if they confirm
or create faith they are true miracles. 

\noindent -- George Bernard Shaw.
% Sain Joan, Scene II 
\end{verse}
}

\begin{figure}[h]
 \begin{center}
     \includegraphics[scale=0.2]{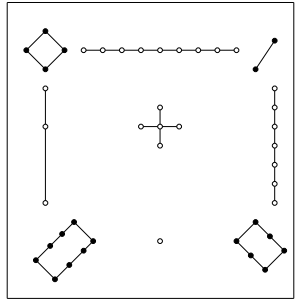}
\caption{The Loh-Shu  magic square \cite{andrews}.} \label{lohshu}
 \end{center}
 \end{figure}

\begin{figure}[h]
 \begin{center}
     \includegraphics[scale=0.9]{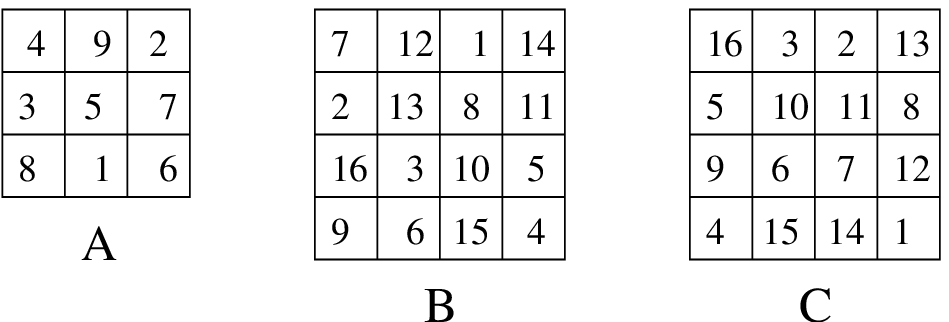}
\caption{Loh-Shu (A), Jaina (B), and the D\"urer (C) Magic squares.} 
\label{ancient}
 \end{center} \end{figure}

A {\em magic square} is a square matrix whose entries are nonnegative
integers, such that the sum of the numbers in every row, in every
column, and in each diagonal is the same number called the {\em magic
sum}. The study of magic squares probably dates back to prehistoric
times \cite{andrews}. The Loh-Shu magic square is the oldest known
magic square and its invention is attributed to Fuh-Hi (2858-2738
B.C.), the mythical founder of the Chinese civilization
\cite{andrews}. The odd numbers are expressed by white dots, i.e.,
yang symbols, the emblem of heaven, while the even numbers are in
black dots, i.e., yin symbols, the emblem of earth (see Figures
\ref{lohshu} and \ref{ancient} A).

%%%% Figure of Durer magic  square %%%%%%%%%%%%%%%%%%%%%%%

%%% Figure of Durer magic  square %%%%%%%%%%%%%%%%%%%%%%%

%%%% Figure of Melancholia %%%%%%%%%%%%%%%%%%%%%%%
\begin{figure}[h]
 \begin{center}
     \includegraphics[scale=0.5]{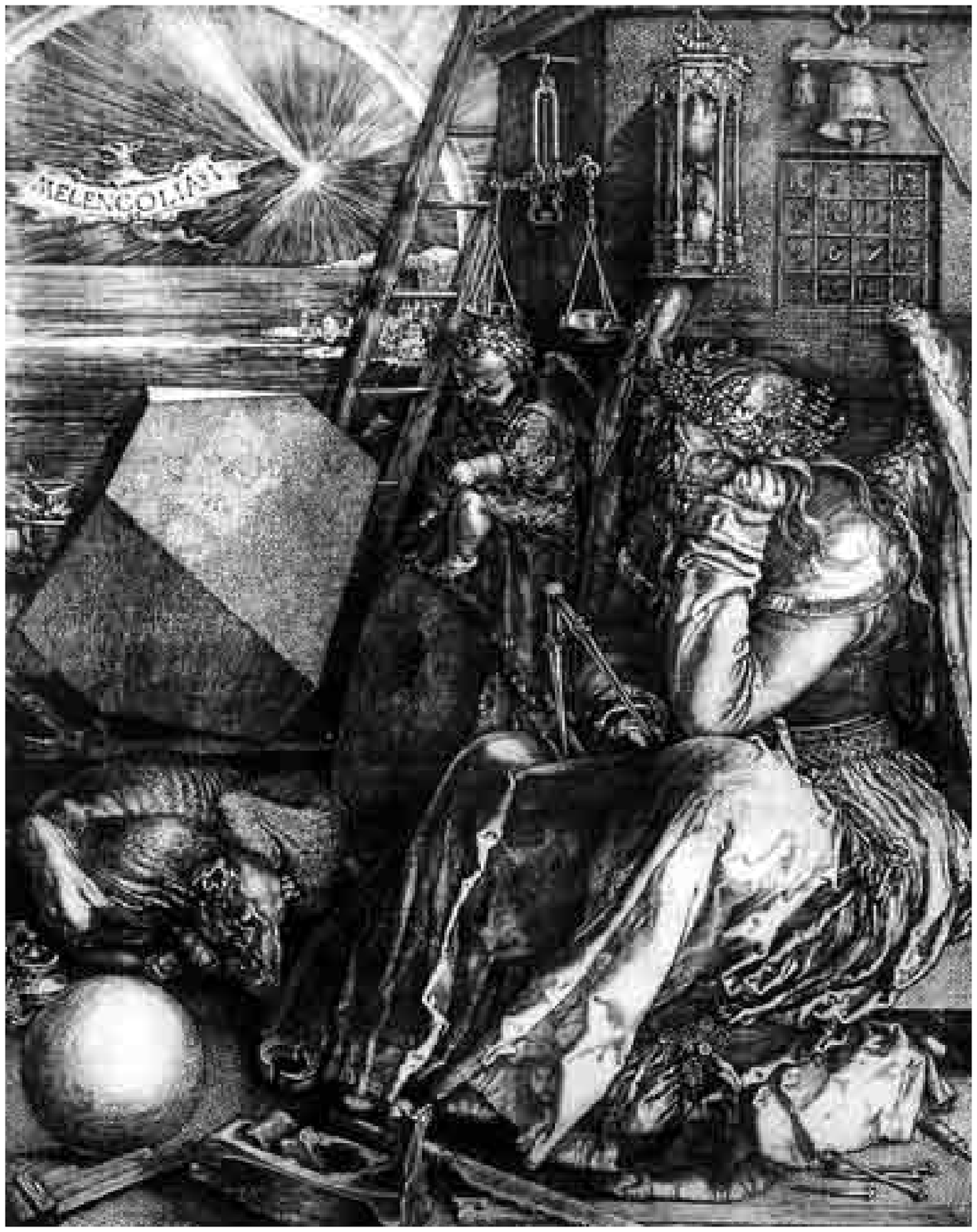}
\caption {Magic square in a 1514 engraving by Albrecht D\"urer
entitled {\em Melancholia} \cite{andrews}.} \label{durer}
 \end{center}
 \end{figure}

\begin{figure}[h]
 \begin{center}
     \includegraphics[scale=0.5]{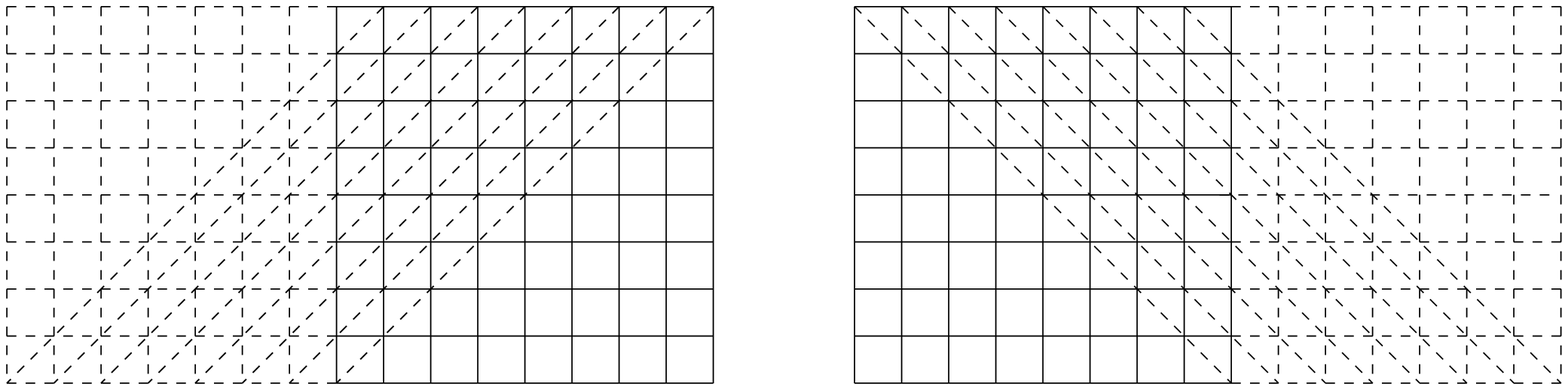}
\caption{The pandiagonals.} \label{pans}
 \end{center}
 \end{figure}

Like chess and many of the problems founded on the figure of the
chess-board, the problem of constructing a magic square probably
traces its origin to India \cite{schubert}. A $4 \times 4$ magic
square is found in a Jaina inscription of the twelfth or thirteenth
century in the city of Khajuraho, India (see Figure \ref{ancient} B).
This square shows an advanced knowledge of magic squares because all the 
pandiagonals also add to the common magic sum (see Figure \ref{pans}). 

%%%%%%% pan-diagonals %%%%%%%%%%%%%%%%%
%%%%%%%%%%%%%%%%%%%%%%%%%%%%%%%%%%%%%%%%%%%%%%%%

From India, the problem found its way among the Arabs, and by them it
was brought to the Roman Orient.  It is recorded that as early as the
ninth century magic squares were used by Arabian astrologers in their
calculations of horoscopes etc. which might be the reason such squares
are called ``magic'' \cite{andrews}. Their introduction into Europe
appears to have been due to Moschopulus, who lived in Constantinople
in the early part of the fifteenth century.  The famous Cornelius
Agrippa (1486-1535) constructed magic squares of the orders
3, 4, 5, 6, 7, 8, 9, which were associated by him with the seven
astrological ``planets''; namely, Saturn, Jupiter, Mars, the Sun,
Venus, Mercury, and the Moon. A magic square engraved on a silver
plate was sometimes prescribed as a charm against the plague and a
magic square appears in a well known 1514 engraving by Albrecht
D\"urer entitled {\em Melancholia} (see Figures \ref{ancient} C and
\ref{durer}). Magic squares, in general, were considered as mystical
objects with the power to ward off evil and bring good fortune. The
mathematical theory of the construction of these squares was taken up
in France only in the seventeenth century, and since then it has
remained a favorite topic of study throughout the mathematical
world. See \cite{andrews}, \cite{ball}, \cite{zen}, or \cite{schubert}
to read more about the history of magic squares.

Constructing and enumerating magic squares are the two fundamental
problems in the topic of magic squares. Let $M_n(s)$ denote the number
of $n \times n$ magic squares of magic sum $s$.  In 1906, MacMahon
\cite{macmahon} enumerated magic squares of order 3:

\[ 
M_3(s) = \left\{ \begin{array}{ll}
{\frac{2}{9}}s^2 + {\frac{2}{3}}s +1 & \mbox{if 3 divides $s$}, \\
0 & \mbox{otherwise.}
\end{array}
\right.
\]

Later, MacMahon \cite{macmahon}, Anand et al. \cite{anandgupta}, and Stanley
\cite{stanley2}, considered the problem of enumerating {\em
semi-magic} squares which are squares such that only the row and column
sums add to a magic sum (see \cite{stanley2}).

In this thesis, we have used methods from algebra, combinatorics, and
polyhedral geometry to construct and enumerate magic squares and these
methods are similar to the methods used by Stanley to enumerate
semi-magic squares \cite{stanley2}. Part of our work is based on
generalizing his ideas.  Polyhedral methods were also used by Halleck
in his 2000 Ph.D thesis (see \cite{halleck}) and Beck et al. in
\cite{becketal} to enumerate various kinds of magic squares. Our
techniques are different and are more algebraic in flavor.

\begin{figure}[h]
 \begin{center}
     \includegraphics[scale=0.7]{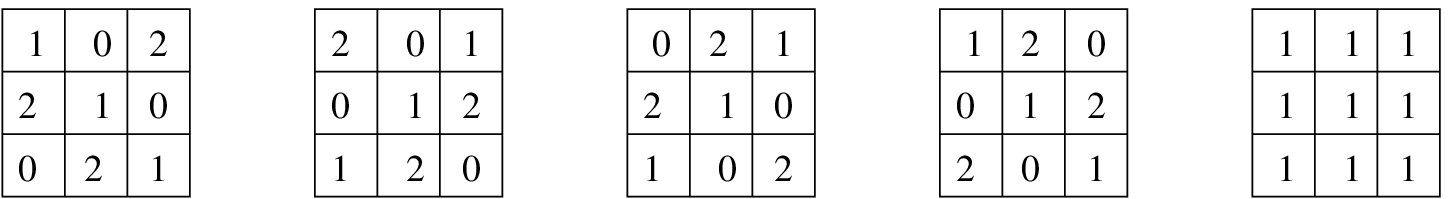}
\caption{The minimal Hilbert Basis of $3 \times 3$ Magic squares.} 
\label{3x3magichilbert}
 \end{center}
 \end{figure}

\begin{figure}[h]
 \begin{center}
     \includegraphics[scale=0.7]{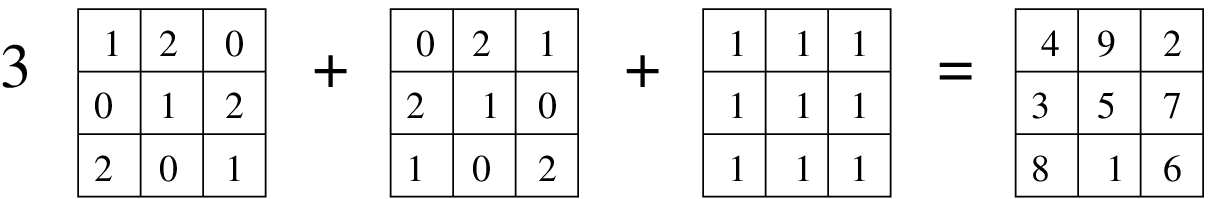}
\caption{A Hilbert basis construction of the Loh-Shu magic square.} 
\label{lohshuconstruct}
 \end{center}
 \end{figure}

The basic idea is to consider the defining linear equations of magic
squares. These equations, together with the nonnegativity requirement
on the entries, imply that the set of magic squares becomes the set of
integral points inside a {\em pointed polyhedral cone} (see \cite{adh}
or \cite{stanley}). The {\em minimal Hilbert basis} of this cone is
defined to be the smallest finite set $S$ of integral points with the
property that any integral point can be expressed as a linear
combination with nonnegative integer coefficients of the elements of
$S$ \cite{schrijver}. For example, the minimal Hilbert basis of the
$3 \times 3$ magic squares is given in Figure \ref{3x3magichilbert}
and a Hilbert basis construction of the Loh-shu magic square is given
in Figure \ref{lohshuconstruct}.

We map magic squares to monomials in a polynomial ring to enumerate
them and derive the formula for the number of $4 \times 4$ magic
squares of magic sum $s$ (also derived  simultaneously by Beck et
al. \cite{becketal}  using different techniques).

\begin{thm}  \label{4x4magicformula}

The number of $4 \times 4$ magic squares with magic sum $s$,

{\footnotesize
\[
M_4(s) = \left \{
\begin{array}{l}

\frac{1}{480}s^7+\frac{7}{240}s^6+\frac{89}{480}s^5+\frac{11}{16}s^4+\frac{779}{480}s^3
+\frac{593}{240}s^2+\frac{1051}{480}s+\frac{13}{16},
\\ \\ \hfill \mbox{when $s$ is odd}, \\ \\ \\ \frac{1}{480}s^7
+\frac{7}{240}s^6 +\frac{89}{480}s^5 +\frac{11}{16}s^4 +\frac{
49}{30}s^3 +\frac{38}{15}s^2 +\frac{71}{30}s + 1, \\ \\ \hfill \mbox{
when $s$ is even.}
\end{array}
\right .
\]
}
\end{thm}

We will describe this method in detail in the next section (also, see
\cite{adh} and \cite{ma}).

\section{Generating and enumerating  lattice points inside polyhedral
cones.} \label{method}

In this section we illustrate the algebraic techniques of generating
and enumerating lattice points inside polyhedral cones by applying the
method to construct and enumerate magic squares. As an example, we
walk through the details of the proof of Theorem \ref{4x4magicformula}.

For these purposes we regard $n \times n$ magic squares as either $n
\times n$ matrices or vectors in $\reals^{n^2}$ and apply the normal
algebraic operations to them. We also consider the entries of an $n
\times n$ magic square as variables $y_{ij} \hspace{0.05in} (1 \leq
i,j \leq n)$. If we set the first row sum equal to all other mandatory
sums, then magic squares become nonnegative integral solutions to a
system of linear equations $Ay = 0$, where $A$ is an $(2n+1)
\times n^2$ matrix each of whose entries is 0, 1, or -1. 

For example, the equations defining $3 \times 3$ magic squares are:

\[
\begin{array}{c}
y_{11}+y_{12}+y_{13} = y_{21}+y_{22}+y_{23} \\
y_{11}+y_{12}+y_{13} = y_{31}+y_{32}+y_{33}  \\
y_{11}+y_{12}+y_{13} = y_{11}+y_{21}+y_{31} \\
y_{11}+y_{12}+y_{13} = y_{12}+y_{22}+y_{32}  \\
y_{11}+y_{12}+y_{13} = y_{13}+y_{23}+y_{33}  \\
y_{11}+y_{12}+y_{13} = y_{11}+y_{22}+y_{33} \\
y_{11}+y_{12}+y_{13} = y_{13}+y_{22}+y_{31} 
\end{array}
\]

Therefore, $3 \times 3$ magic squares are nonnegative integer solutions
to the system of equations $Ay=0$ where:

\[
\begin{array}{lll}
A = \left[
\begin{array}{rrrrrrrrr}
1 & 1 & 1 & -1 & -1 & -1 & 0 & 0 & 0 \\
1 & 1 & 1 & 0 & 0 & 0 & -1 & -1 & -1 \\
0 & 1 & 1 & -1 & 0 & 0 & -1 & 0 & 0 \\
1 & 0 & 1 & 0 &-1 & 0 & 0 & -1 & 0 \\
1 & 1 & 0 & 0 &0 & -1 & 0 & 0 & -1  \\
0 & 1 & 1 & 0 & -1 & 0 & 0 & 0 & -1  \\
1 & 1 & 0 & 0 & -1 & 0 & -1& 0 & 0  \\
\end{array}
\right] & \mbox{and} &
y = \left[
\begin{array}{r}
y_{11} \\ y_{12} \\ y_{13} \\
y_{21} \\ y_{22} \\ y_{23} \\
y_{31} \\ y_{32} \\ y_{33} 
\end{array}
\right] 
\end{array}
\]

A nonempty set $C$ of points in $\reals^{n^2}$ is a {\em cone} if $au
+ bv$ belongs to $C$ whenever $u$ and $v$ are elements of $C$ and $a$
and $b$ are nonnegative real numbers. A cone is {\em pointed} if the
origin is its only vertex (or minimal face; see \cite{schrijver}). A
cone $C$ is {\em polyhedral} if $C = \{y: Ay \leq 0 \}$ for some
matrix $A$, i.e, if $C$ is the intersection of finitely many
half-spaces. If, in addition, the entries of the matrix $A$ are
rational numbers, then $C$ is called a {\em rational} polyhedral cone.
A point $y$ in the cone $C$ is called an {\em integral point} if all
its coordinates are integers.

It is easy to verify that the sum of two magic squares is a magic
square and that nonnegative integer multiples of magic squares are
magic squares. Therefore, the set of magic squares is the set of all
integral points inside a polyhedral cone $C_{M_n}=\{y:Ay=0, y \geq
0\}$ in $\reals^{n^2}$, where $A$ is the coefficient matrix of the
defining linear system of equations. Observe that $C_{M_n}$ is a
pointed cone.

In the example of $4 \times 4$ magic squares, there are three linear
relations equating the first row sum to all other row sums and four
more equating the first row sum to column sums. Similarly, equating
the two diagonal sums to the first row sum generates two more linear
equations.  Thus, there are a total of 9 linear equations that define
the cone of $4 \times 4$ magic squares.  The coefficient matrix $A$
has rank 8 and therefore the cone $C_{M_4}$ of $4 \times 4$ magic
squares has dimension $16-8=8$ (see \cite{schrijver}).

In 1979, Giles and Pulleyblank introduced the notion of a {\em Hilbert
basis} of a cone \cite{gilespulley}. For a given cone $C$, its set
$S_C = C \cap {\integers}^n$ of integral points is called the {\em
semigroup of the cone} $C$. 

\begin{defn}
A Hilbert basis for a cone $C$ is a finite set of points $HB(C)$ in
its semigroup $S_C$ such that each element of $S_C$ is a linear
combination of elements from $HB(C)$ with nonnegative integer
coefficients.
\end{defn}

\begin{figure}[h] 
 \begin{center}
     \includegraphics[scale=0.4]{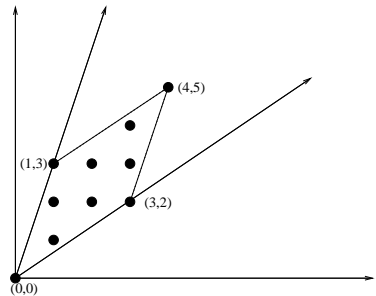}
\caption{A Hilbert Basis of a two dimensional cone.} \label{hilbexample}
 \end{center}
 \end{figure}

For example, the integral points inside and on the boundary of the
parallelepiped in $\reals^2$ with vertices $(0,0), (3,2), (1,3)$ and
$(4,5)$ in Figure \ref{hilbexample} form a Hilbert basis of the cone
generated by the vectors $(1,3)$ and $(3,2)$.

We recall an important fact about Hilbert bases \cite[Theorem
16.4]{schrijver}:
\begin{thm} \label{hilbtheorem}
Each rational polyhedral cone $C$ is generated by a Hilbert basis. If
$C$ is pointed, then there is a unique minimal integral Hilbert basis
generating $C$ (minimal relative to taking subsets).
\end{thm}

We present a proof of Theorem \ref{hilbtheorem} in Appendix
\ref{appendixa}. The minimal Hilbert basis of a pointed cone is unique
and henceforth, when we say the Hilbert basis,  we mean the minimal
Hilbert basis. An integral point of a cone $C$ is {\em irreducible} if
it is not a linear combination with integer coefficients of other
integral points.  All the elements of the minimal Hilbert basis are
irreducible \cite{raymond}, \cite{schrijver}. Since magic squares are
integral points inside a cone, Theorem \ref{hilbtheorem} implies that
every magic square is a nonnegative integer linear combination of
irreducible magic squares.

\begin{figure}[h]
 \begin{center}
     \includegraphics[scale=0.5]{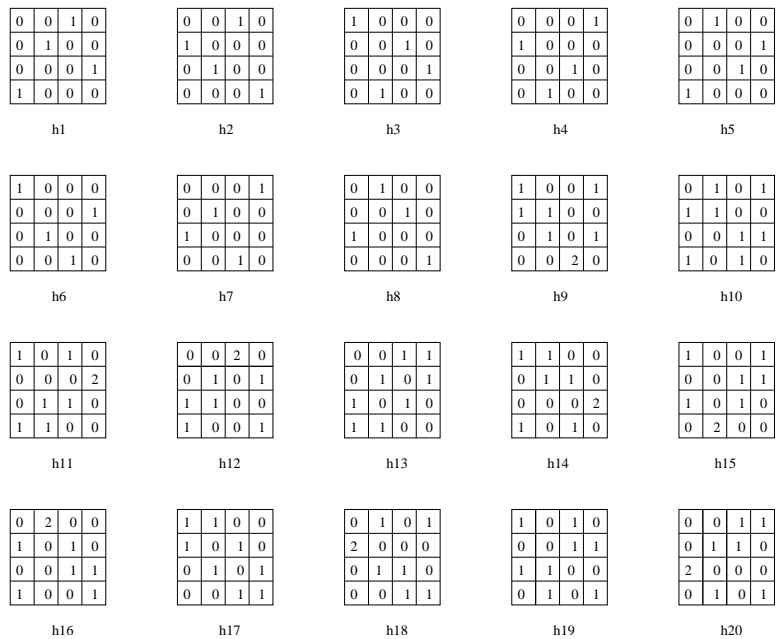}
\caption{The minimal Hilbert Basis of $4 \times 4$ Magic squares.} 
\label{4x4magichilbert}
 \end{center}
 \end{figure}

 The minimal Hilbert basis of the polyhedral cone of $4 \times 4$
magic squares is given in Figure \ref{4x4magichilbert} and was
computed using the software MLP (now called 4ti2) (see \cite{raymond};
software implementation 4ti2 is available from http://www.4ti2.de).
In fact, all Hilbert basis calculations in this thesis are done with
the software 4ti2. 

\begin{figure}[h]
 \begin{center}
     \includegraphics[scale=0.5]{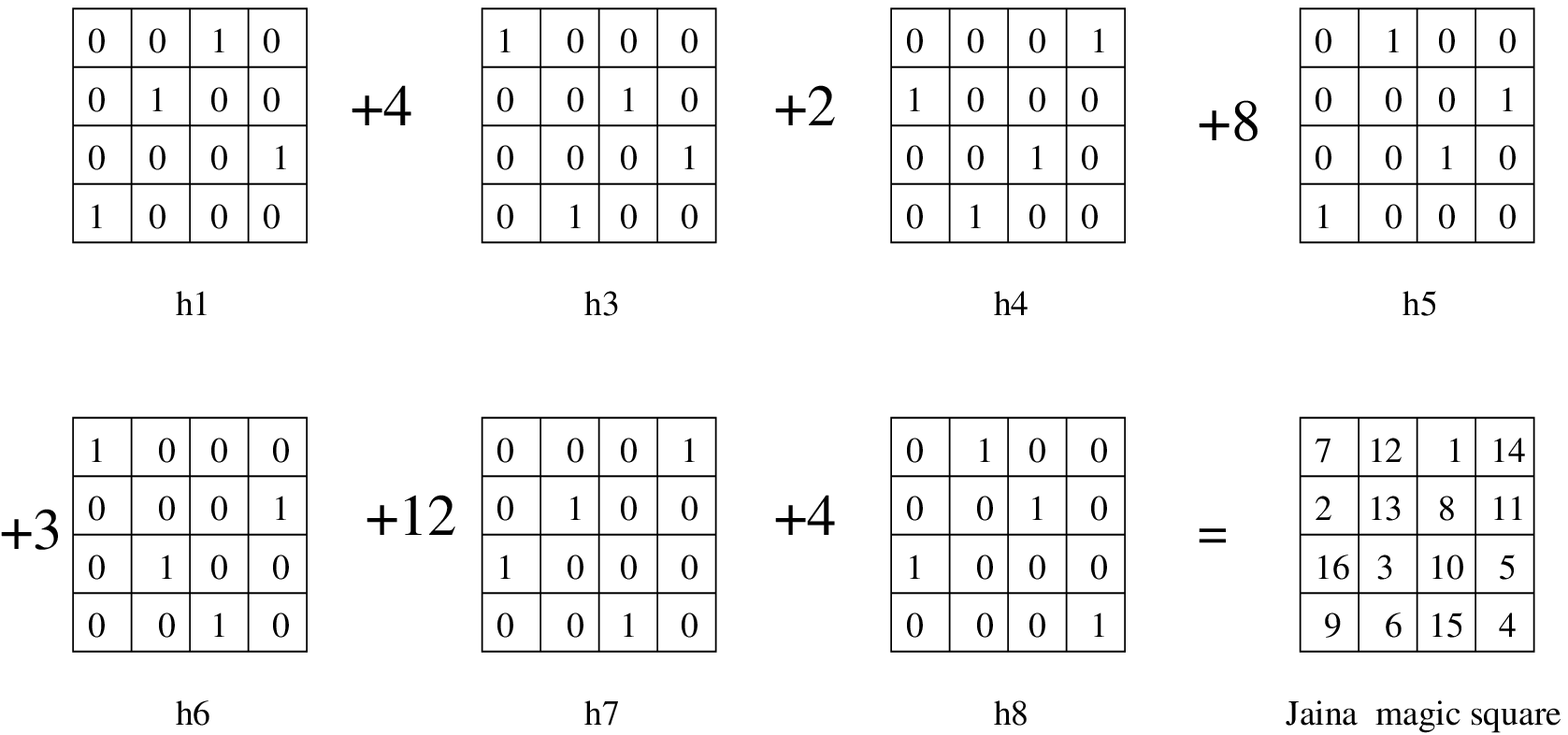}
\caption{A Hilbert basis construction of the Jaina  magic square.} 
\label{Jainaconstruct}
 \end{center}
 \end{figure}

\begin{figure}[h]
 \begin{center}
     \includegraphics[scale=0.5]{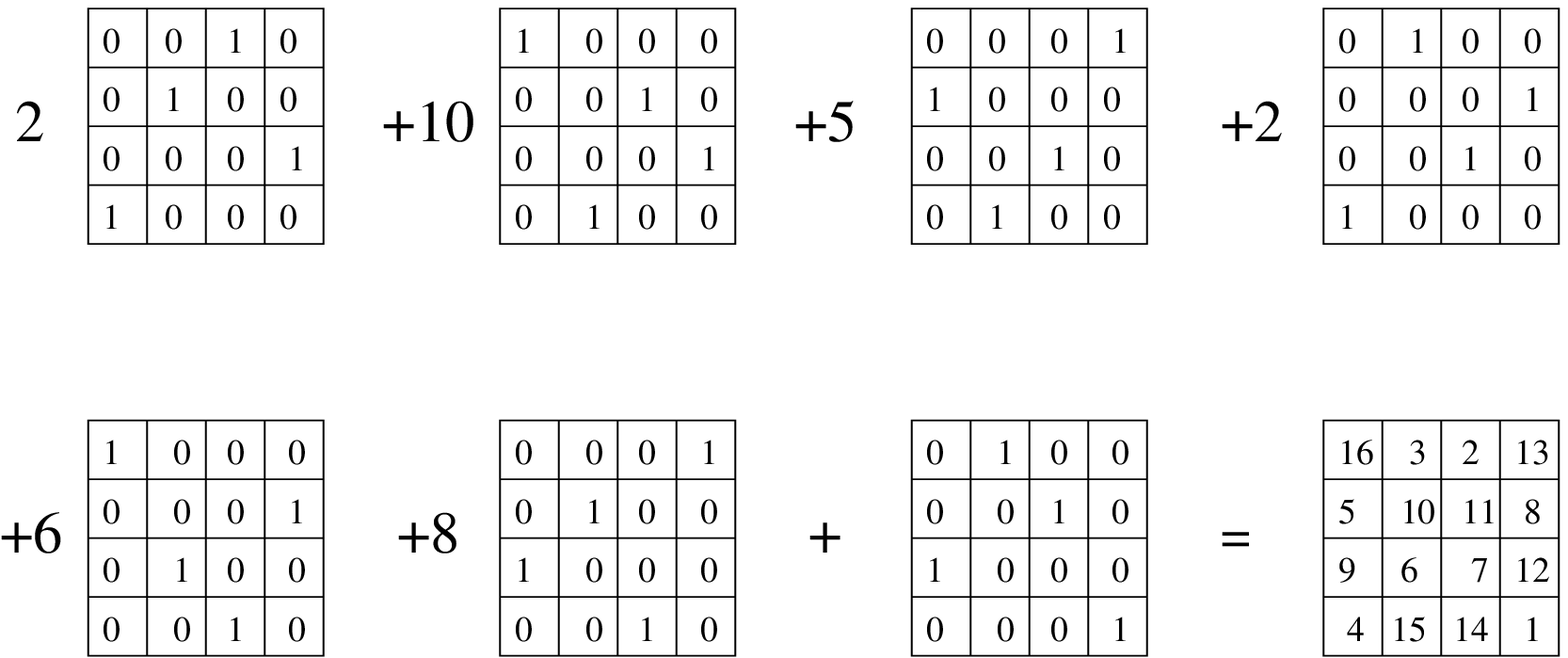}
\caption{A Hilbert basis construction of D\"urer's magic square.
} \label{durerconstruct}
 \end{center}
 \end{figure}

Hilbert basis constructions of the Jaina, and the D\"urer magic
squares are given in Figures \ref{Jainaconstruct} and
\ref{durerconstruct} respectively.

\begin{figure}[h]
 \begin{center}
     \includegraphics[scale=0.5]{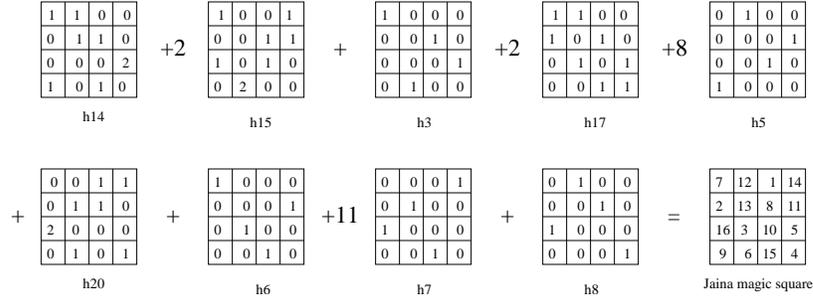}
\caption{Another Hilbert basis construction of the Jaina  magic square.} 
\label{anotherjainaconstruct}
 \end{center}
 \end{figure}

Different combinations of the elements of a Hilbert basis sometimes
produce the same magic square. Figures \ref{Jainaconstruct} and
\ref{anotherjainaconstruct} exhibit two different Hilbert basis
constructions of the Jaina magic square. This is due to algebraic
dependencies among the elements of the Hilbert basis.  Repetitions
have to be avoided when counting squares, a problem that we solve by
using a little bit of commutative algebra. Let $HB(C_{M_n}) = \{ h_1,h_2,
\dots h_r \}$ be a Hilbert basis for the cone of $n \times n$ magic
squares.  Denote the entries of the square $h_p$ by $y^p_{ij}$, and
let $k$ be any field. Let $\phi$ be the unique ring homomorphism
between the polynomial rings $k[x_1,x_2, \dots, x_r]$ and
$k[t_{11},t_{12}, \dots, t_{1n},t_{21},t_{22}, \dots t_{2n}, \dots,
t_{n1},t_{n2}, \dots, t_{nn}]$ such that $\phi(x_p) = t^{h_p}$, the
monomial defined by
\[
 t^{h_p} = \prod_{i,j=1, \dots, n} t_{ij}^{y^{p}_{ij}}.
\]

Monomials in $k[x_1,x_2, \dots, x_r]$ correspond to magic squares
under this map, and multiplication of monomials corresponds to
addition of magic squares.  For example, the monomial
$x_1^5x_{3}^{200}$ corresponds to the magic square $5h_1 +
200h_{3}$. Different combinations of Hilbert basis elements that give
rise to the same magic square can then be represented as polynomial
equations. Thus, from the two different Hilbert basis constructions of
the Jaina magic square represented in Figures \ref{Jainaconstruct}
and \ref{anotherjainaconstruct}, we learn that
\[
\begin{array}{l}
h1 + 4 \cdot h3 + 2 \cdot h4 + 8 \cdot h5 + 3 \cdot h6 + 12 \cdot h7 
+ 4 \cdot h8 = \\
h3  + 8 \cdot h5 + h6 + 11 \cdot h7 + h8 + h14 + 2 \cdot h15 + 
2 \cdot h17 + h20 
\end{array}
\]

\noindent In $k[x_1,x_2, \dots, x_r]$, this algebraic dependency 
of Hilbert basis elements translates to
\[
\begin{array}{l}
x_1x_3^{4}x_4^{2}x_{5}^{8}x_{6}^3x_7^{12}x_{8}^{4} -
x_{3}x_5^8 x_6x_7^{11}x_8x_{14}x_{15}^2x_{17}^{2}x_{20} = 0.
\end{array}
\]

Consider the set of all polynomials in $k[x_1,x_2, \dots, x_r]$ that
are mapped to the zero polynomial under $\phi$. This set, which
corresponds to all the algebraic dependencies of Hilbert basis
elements, forms an ideal in $k[x_1,x_2, \dots, x_r]$, an ideal known
as the {\em toric ideal} of $HB(C_{M_n})$ (see \cite{adh},
\cite{bigatti1}, or \cite{sturmfels} for details about toric
ideals). If we denote the toric ideal as $I_{HB(C_{M_n})}$, then the
monomials in the quotient ring $R_{C_{M_n}} = k[x_1,x_2,
\cdots,x_r]/I_{HB(C_{M_n})}$ are in one-to-one correspondence with
magic squares.

For example, in the case of $3 \times 3$ magic squares, there are 5 Hilbert
basis elements (see Figure \ref{3x3magichilbert}) and hence there are
5 variables $x_1, x_2, x_3, x_4, x_5$ which gets mapped by $\phi$ as follows:
{\small
\[
\begin{array}{cccc}
x_1 \mapsto &
\left [
\begin{array}{ccc}
1 & 0 & 2 \\
2 & 1 & 0 \\
0 & 2 & 1
\end{array}
\right ] 
\mapsto & t_{11}{t_{13}}^2 {t_{21}}^2 t_{22}{t_{32}}^2 t_{33}
\\ \\

x_2 \mapsto &

\left [
\begin{array}{ccc}
2 & 0 & 1 \\
0 & 1 & 2 \\
1 & 2 & 0
\end{array}
\right ] 

\mapsto & {t_{11}}^2 t_{13}t_{22}{t_{23}}^2 t_{31} {t_{32}}^2
\\ \\

x_3 \mapsto &

\left [
\begin{array}{ccc}
0 & 2 & 1 \\
 2 & 1 & 0 \\
 1 & 0 & 2
\end{array}
\right ] 
\mapsto & {t_{12}}^2 t_{13} {t_{21}}^2 t_{22}t_{31}{t_{33}}^2
\\ \\

x_4 \mapsto &
\left [
\begin{array}{ccc}
1 & 2 & 0 \\
0 & 1 & 2  \\
2 & 0 & 1
\end{array}
\right ] 
\mapsto & t_{11} {t_{12}}^2 t_{22} {t_{23}}^2 {t_{31}}^2 t_{33}
\\ \\

x_5 \mapsto &
\left [
\begin{array}{ccc}
1 & 1 & 1 \\
1 & 1 & 1 \\
1 & 1 & 1 
\end{array}
\right ]
\mapsto & t_{11} t_{12} t_{13} t_{21} t_{22} t_{23} t_{31} t_{32} t_{33}

\end{array}
\]
}

We use CoCoA to compute the toric ideal 

\[
I_{HB(C_{M_3})} = (x_1x_4 - x_5^2,x_2x_3 - x_1x_4).  
\]

Thus

\[ 
R_{C_{M_3}} = \frac{Q[x_1,x_2,x_3,x_4,x_5]}{(x_1x_4 - x_5^2,x_2x_3 - x_1x_4)}.
\]

Let $R_{C_{M_n}}(s)$ be the set of all homogeneous polynomials of
degree $s$ in the ring $R_{C_{M_n}}$.  Then $R_{C_{M_n}}(s)$ is a
$k$-vector space, and $R_{C_{M_n}}(0)=k$. The dimension
$\mbox{dim}_k(R_{C_{M_n}}(s))$ of $R_{C_{M_n}}(s)$ is precisely the
number of monomials of degree $s$ in $R_{C_{M_n}}$. 

Recall that if the variables $x_i$ of a polynomial ring $k[x_1,x_2,
\dots, x_r]$ are assigned nonnegative weights $w_i$, then the weighted
degree of a monomial $x_1^{\alpha_1} \cdots x_r^{\alpha_r}$ is
$\sum_{i=1}^r \alpha_i \cdot w_i$ (see \cite{atiyah}). Therefore, if
we take the weight of the variable $x_i$ to be the magic sum of the
corresponding Hilbert basis element $h_i$, then
$\mbox{dim}_k(R_{C_{M_n}}(s))$ is exactly the number of magic squares
of magic sum $s$. For example, in the case of $3 \times 3$ magic
squares, because all the elements of the Hilbert basis have sum 3, all
the variables are assigned degree 3.

Since $R_{C_{M_n}}$ is a graded $k$-algebra, it can be decomposed into
a direct sum of its graded components $R_{C_{M_n}}= \bigoplus
R_{C_{M_n}}(s)$ (see \cite{adh} or \cite{atiyah}). Consider a finitely
generated graded $k$-algebra $R_{C_{M_n}}= \bigoplus R_{C_{M_n}}(s)$.
The function $H(R_{C_{M_n}},s)=\mbox{dim}_k(R_{C_{M_n}}(s))$ is the
{\em Hilbert function} of $R_{C_{M_n}}$ and the {\em
Hilbert-Poincar\'e series} of $R_{C_{M_n}}$ is the formal power series
$$ H_{R_{C_{M_n}}}(t)=\sum_{s=0}^{\infty} H(R_{C_{M_n}},s)t^s.$$

We can now deduce the following lemma.
\begin{lemma}
Let the weight of a variable $x_i$ in the ring $R =
k[x_1,x_2,...,x_r]$ be the magic sum of the corresponding element of
the Hilbert basis $h_i$. With this grading of degrees on the monomials
of $R$, the number of distinct magic squares of magic sum $s$
is given by the value of the Hilbert function $H(R_{C_{M_n}},s)$.
\end{lemma}

We record here a version of the Hilbert-Serre Theorem. A proof of the
Hilbert-Serre theorem in all its generality is presented in Appendix
\ref{appendixa}.

\begin{thm}[Theorem 11.1 \cite{atiyah}] Let $k$ be a field and $R:=
k[x_1,x_2,...,x_r]$ be a graded Noetherian ring. let $x_1,x_2,...,x_r$
be homogeneous of degrees $>0$. Let $M$ be a finitely generated
$R$-module.  Then the Hilbert Poincar\'e series of $M$, $H_{M}(t)$ is a 
rational function of the form:
$$H_{M}(t)=\frac{p(t)}{\Pi_{i=1}^r(1-t^{deg x_i})},$$
where $p(t) \in {\integers}[t]$. 
\end{thm}

 By invoking the Hilbert-Serre theorem, we conclude that the
Hilbert-Poincar\'e series is a rational function of the form
$H_{R_{C_{M_n}}}(t)={p(t)}/{\Pi_{i=1}^r(1-t^{deg x_i})}$, where $p(t)$
belongs to ${\integers}[t]$.  We refer the reader to \cite{adh},
\cite{atiyah}, \cite{bigatti2}, or \cite{stanley} for information about
Hilbert-Poincar\'e series.

\begin{figure}[h]
 \begin{center}
     \includegraphics[scale=0.5]{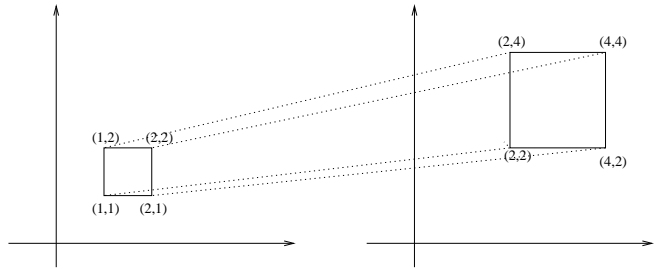}
\caption{Dilation of a polytope.} \label{dilation}
 \end{center}
 \end{figure}

We can also arrive at the conclusion that $H_{R_{C_{M_n}}}(t)$ is a
rational function by studying rational polytopes. A polytope {$\cal P$}
is called {\em rational} if each vertex of {$\cal P$} has rational
coordinates. The {\em dilation of a polytope ${\cal P}$ by an integer
$s$} is defined to be the polytope $s{\cal P} =
\{s{\alpha}:{\alpha}\in {\cal P}\}$ (see Figure \ref{dilation} for an
example).

Let $i({\cal P} ,s)$ denote the number of integer points inside the
polytope $s{\cal P}$. If ${\alpha \in {\rationals}^m}$, let den
${\alpha}$ be the least positive integer $q$ such that $q{\alpha} \in
{\integers}^m$. 
 
\begin{thm}[Theorem 4.6.25 \cite{stanley}] \label{Ehrhart}
Let $\cal P$ be a rational convex polytope of dimension $d$ in
${\reals}^m$ with vertex set $V$. Set $F({\cal P}, t) = 1 + {\sum}_{n
\geq 1}i({\cal P},s)t^s$. Then $F({\cal P},t)$ is a rational function,
which can be written with denominator ${\prod}_{{\alpha} \in V}(1 -
t^{\mbox{den } {\alpha}})$.
\end{thm}

To extract explicit  formulas from the generating function we need to define
the concept of {\em quasi-polynomials}.

\begin{defn}
A function $f:{\naturals} \mapsto {\complex}$
is a {\em quasi-polynomial} if there exists an integer $N > 0$ and
polynomials $f_0,f_1,...,f_d$ such that
\[
\begin{array}{ll}
 f(n) = f_i(n)  & if \  n \equiv i(\mbox{mod} N).
\end{array}
 \]
The integer $N$ is called a quasi-period of $f$.
\end{defn}

For example, the formula for the number of $3 \times 3$ magic squares
 of magic sum $s$ is a quasi-polynomial with quasi-period 3 (see
 Section \ref{magicsquares}). We now state some properties of
 quasi-polynomials.

\begin{prop}[Corollary 4.3.1 \cite{stanley}] \label{quasi}
The following conditions on a function $f:{\naturals} \mapsto {\complex}$ and 
integer $N > 0$ are equivalent:
\begin{enumerate}

\item $f$ is a quasi-polynomial of quasi-period $N$.

\item ${\sum}_{n \geq 0}f(n)x^n = \frac{P(x)}{Q(x)}$,

where $P(x)$ and  $Q(x) \in {\complex}[x]$, every zero $\alpha$ of $Q(x)$ satisfies
${\alpha}^N = 1$ (provided ${P(x)}/{Q(x)}$ has been reduced to lowest
terms) and deg $P < $  deg $Q$. 

\item For all $n \geq 0$,
\[
f(n) = {\sum}_{i=1}^k P_i(n){\gamma}_i^n
\]
where each $P_i$ is a polynomial function of $n$ and each ${\gamma}_i$ 
satisfies ${\gamma}_i^N = 1.$ The degree of  $P_i(n)$ is one less than
the multiplicity of the root ${\gamma}_i^{-1}$ in ${Q(x)}$ provided
${P(x)}/{Q(x)}$ has been reduced to lowest terms.
\end{enumerate}
\end{prop}

Theorem \ref{Ehrhart} together with Proposition \ref{quasi} imply that
$i({\cal P},s)$ is a quasi-polynomial and is generally called the {\em
Ehrhart quasi-polynomial of $\cal P$}. A polytope is called an {\em
integral polytope} when all its vertices have integral
coordinates. $i({\cal P},s)$ is a polynomial if $\cal P$ is an
integral polytope \cite{stanley}.

The convex hull of $n \times n$ matrices, with nonnegative real
entries, such that all the row sums, the column sums, and the diagonal
sums equal $1$, is called the {\em polytope of stochastic magic
squares}. Then, clearly $M_n(s)$ is the Ehrhart quasi-polynomial of
the polytope of stochastic magic squares. Therefore, by Theorem
\ref{Ehrhart}, we get, again, that $H_{R_{C_{M_n}}}(t)$ is a rational
function.

Coming back to the case of $4 \times 4$ magic squares, we used the
program CoCoA (see \cite{cocoa}; CoCoA software is available from
http://cocoa.dima.unige.it) to compute the Hilbert-Poincar\'e series
$\sum_{s=0}^{\infty}M_4(s)t^s$ and obtained

\[
\begin{array}{l}
\sum_{s=0}^{\infty} M_4(s)t^s= \frac{t^8+4t^7+18t^6+36t^5+50t^4+36t^3+18t^2+4t+1}{(1-t)^4(1-t^2)^4}= \\ \\
1+8t+48t^2+200t^3+675t^4+1904t^5+4736t^6+10608t^7+21925t^8+\ldots
\end{array}
\]

Recall that the coefficient of $t^s$ is the number of magic squares of
magic sum $s$. This information along with Proposition \ref{quasi}
enables us to recover the Hilbert function $M_4(s)$ in Theorem
\ref{4x4magicformula} from the Hilbert-Poincar\'e series by
interpolation. We provide basic algorithms to compute Hilbert bases,
toric ideals, and Hilbert-Poincar\'e series in Appendix \ref{appendixb}. 

\begin{figure}[hpt]
 \begin{center}
     \includegraphics[scale=.3]{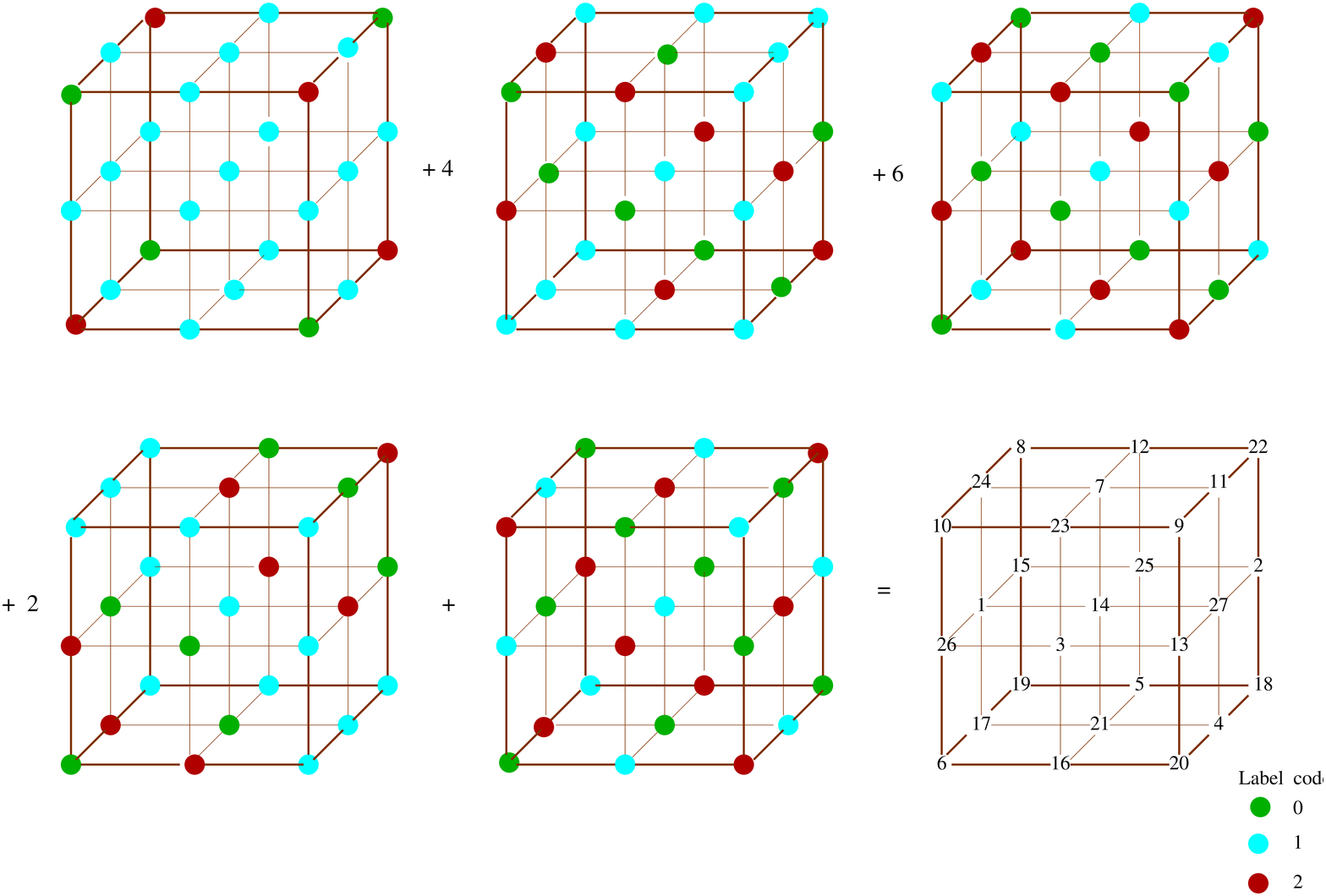}
 \caption{A Hilbert basis construction of a $3 \times 3 \times 3$ magic cube.}
 \label{magiccubeconstruct}
 \end{center}
 \end{figure}

We can also study magic figures of higher dimensions with our
methods. Similar enumerative results for magic cubes and further
properties of the associated cones are discussed in Chapter
\ref{magiccubechapter}. See Figure \ref{magiccubeconstruct} for a
Hilbert basis construction of a magic cube.

In the next section we apply these methods to construct and enumerate
pandiagonal magic squares.

%%%%%%%%%%%%%%%%%%%%%%%%%%%%%%%%%%%%%%%%%%%%%%%%%%%%%%%%%%%%%%%%
%%%%%%%%%%%%%%%%%%%%%%%%%%%%%%%%%%%%%%%%%%%%%%%%%%%%%%%%%%%%%%%%
%%%%%%%%%%%%%%%%%%%%%%%%%%%%%%%%%%%%%%%%%%%%%%%%%%%%%%%%%%%%%%%%
%%%%%%%%%%%%%%%%%%%%%%%%%%%%%%%%%%%%%%%%%%%%%%%%%%%%%%%%%%%%%%%%
%%%%%%%%%%%%%%%%%%%%%%%%%%%%%%%%%%%%%%%%%%%%%%%%%%%%%%%%%%%%%%%%
\section{Pandiagonal magic squares.}
We continue the the study of pandiagonal magic squares started in
\cite{alvin, halleck}. Here we investigate the Hilbert bases, and
recompute the formulas of Halleck (computed using a different method,
see \cite[Chapters 8 and 10]{halleck}).

The convex hull of $n \times n$ matrices, with nonnegative real
entries, such that all the row sums, the column sums, and the
pandiagonal sums equal $1$, is called the {\em polytope of
panstochastic magic squares}.  The integrality of the polytope of
panstochastic magic squares was fully solved in \cite{alvin}.

\begin{figure}[h]
 \begin{center}
     \includegraphics[scale=0.5]{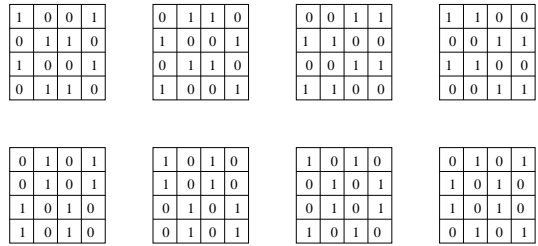}
\caption{Hilbert basis of the $4 \times 4$ Pandiagonal magic squares.}
\label{panmagic4hilb} \end{center}
\end{figure}

\begin{figure}[h]
 \begin{center}
     \includegraphics[scale=0.5]{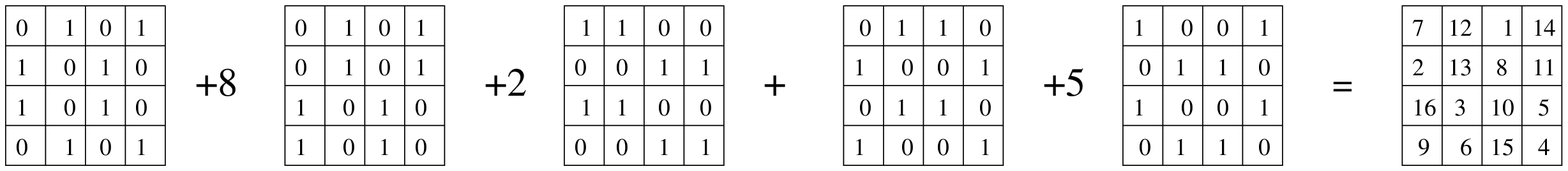}
\caption{A construction of the Jaina magic square with the Hilbert basis 
of  pandiagonal magic squares.} \label{panjainaconstruct}
 \end{center}
 \end{figure}

Let us denote by $MP_n(s)$ the number of $n \times n$ pandiagonal
magic squares with magic sum $s$. Halleck \cite{halleck} computed the
dimension of the cone to be $(n-2)^2$ for odd $n$ and $(n-2)^2+1$ for
even $n$ (degree of the quasipolynomial $MP_n(s)$ is one less than
these).  For the $4\times 4$ pandiagonal magic squares a fast
calculation corroborates that there are 8 elements in the Hilbert
basis (see Figure \ref{panmagic4hilb}). In his investigations, Halleck
\cite{halleck} identified a much larger generating set. Recall that
the Jaina magic square is also a pandiagonal magic square. A
pandiagonal Hilbert basis construction of the Jaina magic square is
given in Figure \ref{panjainaconstruct}.

%%%%%%%%%%%%%%%%%%%%%%%%%%%%%%%%%%%%%%%%%%%%%%%%%%%%%%%%%%%%%%%%

\begin{figure}[h]
 \begin{center}
     \includegraphics[scale=0.5]{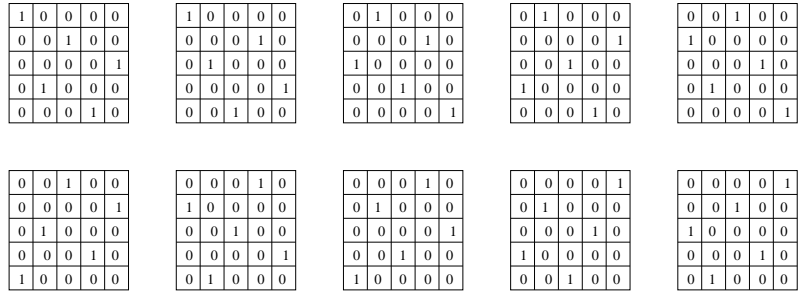}
\caption{Hilbert basis of the $5 \times 5$ Pandiagonal magic squares.}
\label{panmagic5hilb} \end{center}
\end{figure}

We verify that the $5\times 5$ pandiagonal magic squares have indeed a
polynomial counting formula.  This case requires in fact no
calculations thanks to earlier work by \cite{alvin} who proved that
for $n=5$ the only pandiagonal rays are precisely the pandiagonal
permutation matrices. It is easy to see that only 10 of the 120
permutation matrices of order 5 are pandiagonal.  A 4ti2 computation
shows that the set of pandiagonal permutation matrices is also the
Hilbert basis (see Figure \ref{panmagic5hilb}).

%%%%%%%%%%%%%%%%%%%%%%%%%%%%%%%%%%%%%%%%%%%%%%%%%%%%%%%%%%%%%%%%

We calculate the formulas stated in Theorem
\ref{thmpanmagic} using {\tt CoCoA}:

\begin{thm} \label{thmpanmagic}

\[
MP_4(s) = \left\{
\begin{array}{ll}
\frac{1}{48}(s^2+4s+12)(s+2)^2 & \mbox {if $2$ divides $s$,} \\
   0 & \mbox{otherwise.} \cr
\end{array}
\right.
\]

$$MP_5(s)= \frac{1}{8064} (s+4)(s+3)(s+2)(s+1)(s^2+5s+8)(s^2+5s+42).$$

\end{thm}

In the next section, we apply our methods to study Franklin squares
which are more complex than the squares we have seen so far.

%%%%%%%%%%%%%%%%%%%%%%%%%%%%%%%%%%%%%%%%%%%%%%%%%%%%%%%%%%%%%%%%
%%%%%%%%%%%%%%%%%%%%%%%%%%%%%%%%%%%%%%%%%%%%%%%%%%%%%%%%%%%%%%%%
%%%%%%%%%%%%%%%%%%%%%%%%%%%%%%%%%%%%%%%%%%%%%%%%%%%%%%%%%%%%%%%%
%%%%%%%%%%%%%%%%%%%%%%%%%%%%%%%%%%%%%%%%%%%%%%%%%%%%%%%%%%%%%%%%
%%%%%%%%%%%%%%%%%%%%%%%%%%%%%%%%%%%%%%%%%%%%%%%%%%%%%%%%%%%%%%%%

\section{Franklin Squares.}

%%%% Figure of Franklin squares %%%%%%%%%%%%%%%%%%%%%%%
\begin{figure}[h]
 \begin{center}
     \includegraphics[scale=0.5]{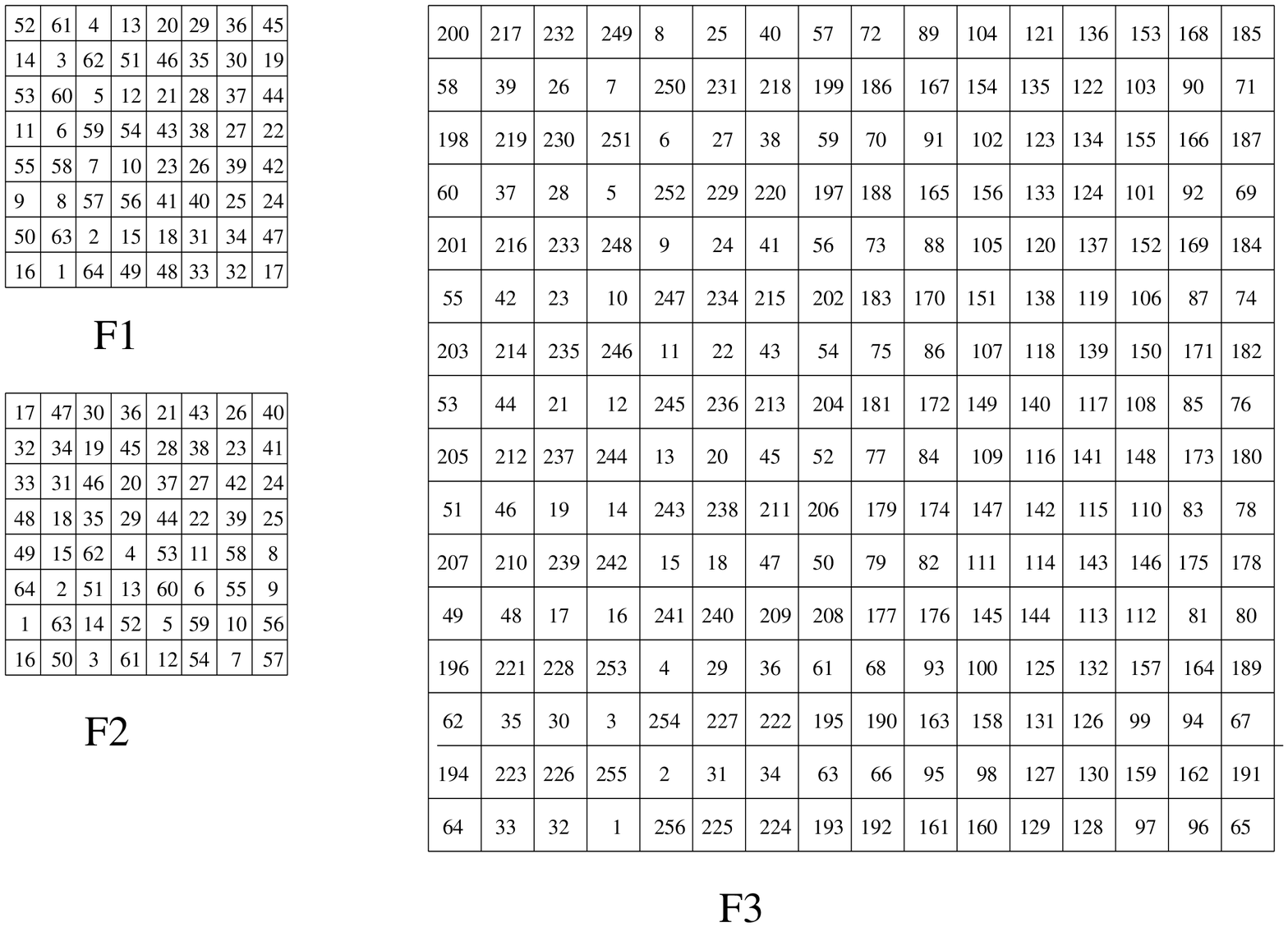}
\caption{Squares constructed by Benjamin Franklin.} \label{franklinsquares}
 \end{center}
 \end{figure}

\begin{figure}[h]
 \begin{center}
     \includegraphics[scale=0.3]{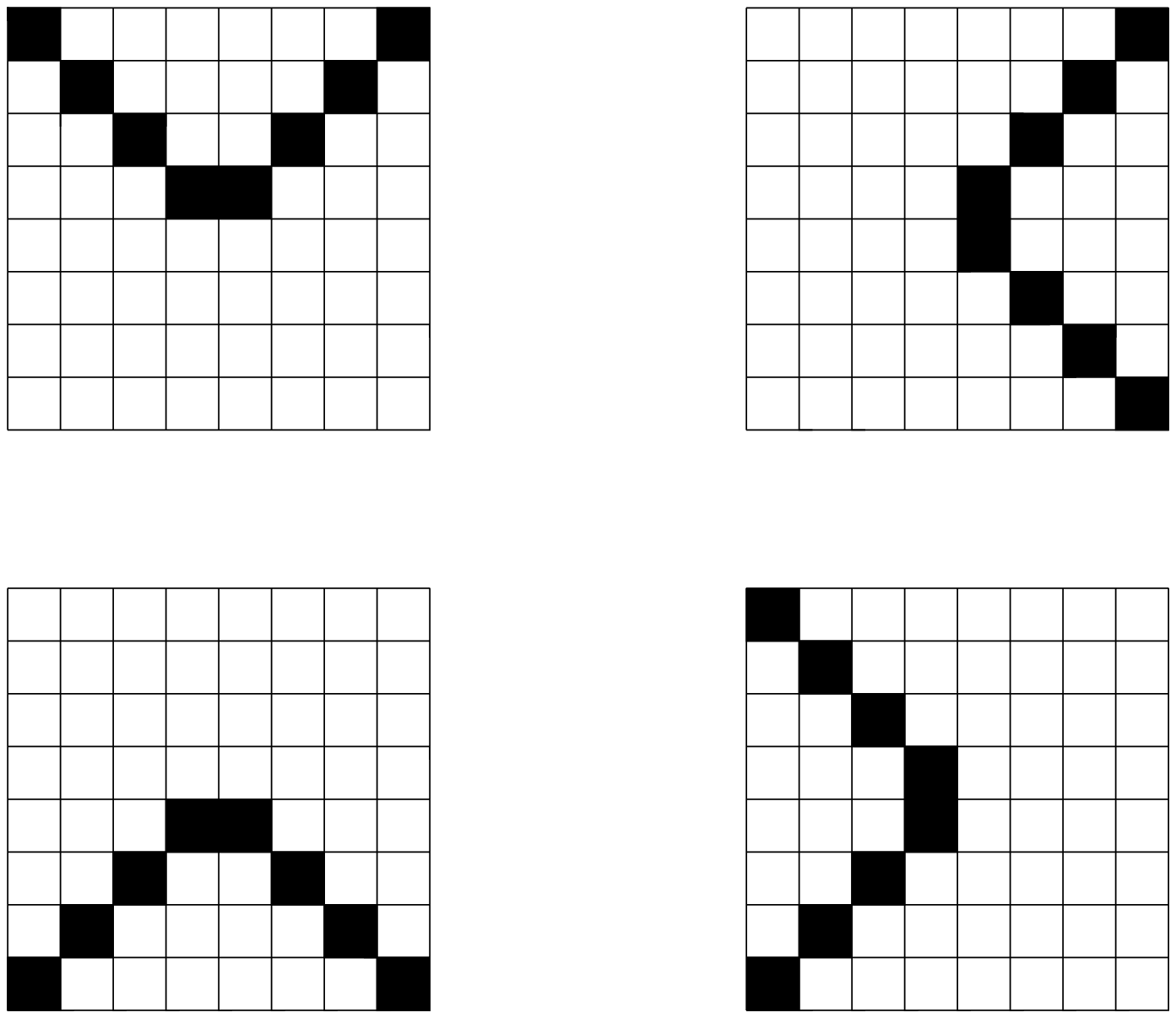}
\caption{The four main bent diagonals \cite{pasles}.} \label{bent_diags}
 \end{center}
 \end{figure}

\begin{figure}[h]
 \begin{center}
     \includegraphics[scale=0.5]{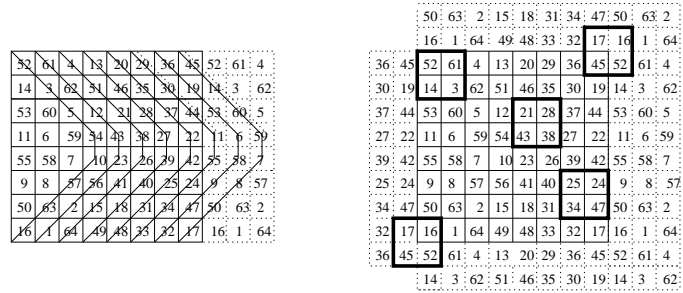}
\caption{Continuous properties of Franklin squares.} \label{parallelbents}
 \end{center}
 \end{figure}

The well-known squares  F1 and F3, as well as the less familiar 
F2, that appear in Figure \ref{franklinsquares} were
constructed by Benjamin Franklin (see \cite{andrews} and
\cite{pasles}). In a letter to Peter Collinson (see \cite{andrews}) he
describes the properties of the $8 \times 8$ square F1 as follows:
\begin{enumerate}
\item The entries of every row and column add to a common sum
called the {\em magic sum}. 
\item In every half-row and half-column the entries add to half the magic sum.
\item The entries of the main bent diagonals (see Figure
\ref{bent_diags}) and all the bent diagonals parallel to it (see
Figure \ref{parallelbents}) add to the magic sum.
\item The four corner entries together with the four middle entries
add to the magic sum.
\end{enumerate}

%%% all properties of 8x8 squares %%%%%%%%
\begin{figure}[h]
 \begin{center}
     \includegraphics[scale=0.5]{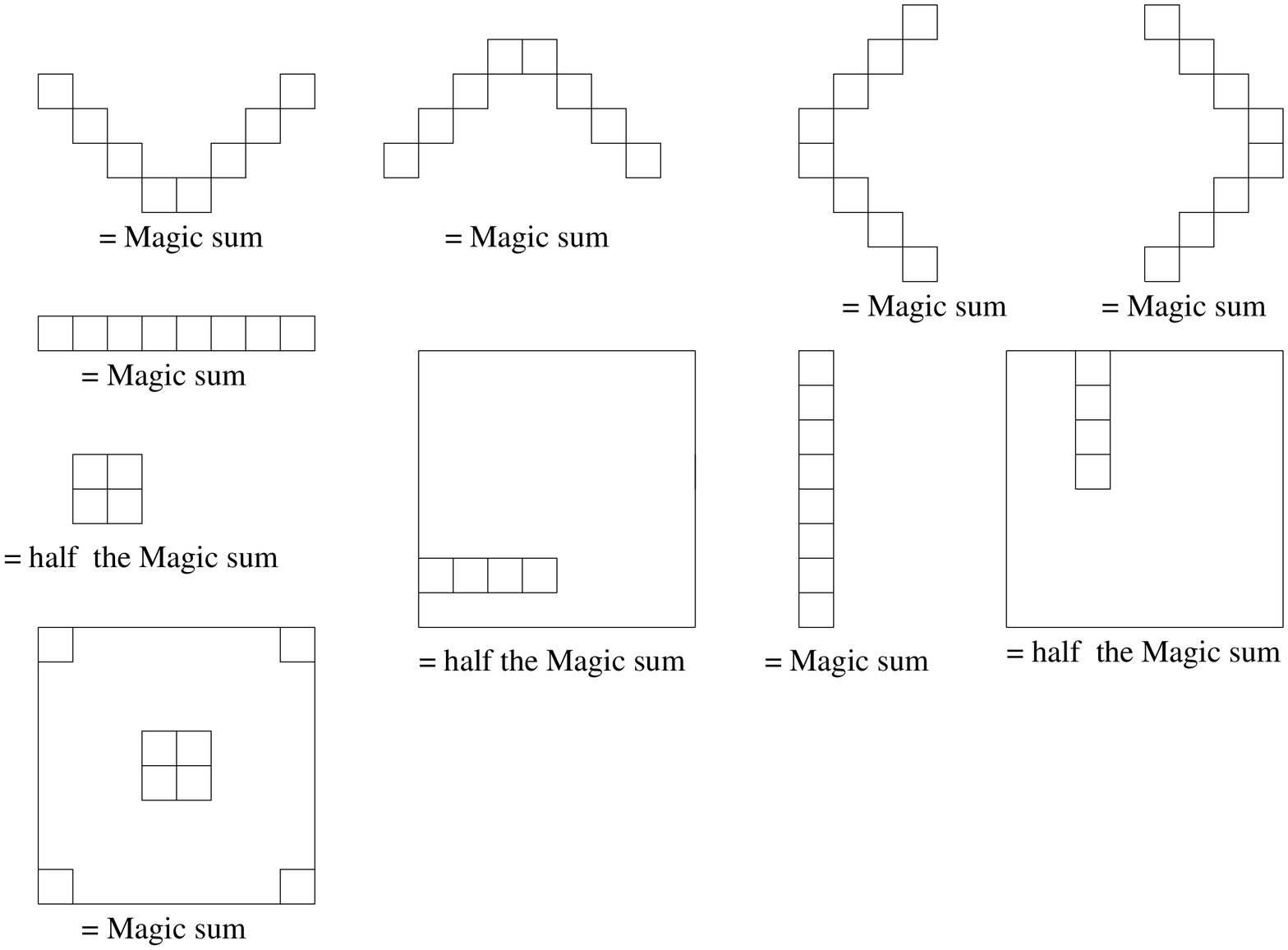}
\caption{Defining properties of the $8 \times 8$ Franklin squares \cite{andrews}.}
\label{8x8} \end{center} \end{figure}
%%%%%%%%%%%%%%%%%%%%%%%%%%%%%%%%%%%%%%%%%%%%%%%%%%%%%%%%%%%%%%%%

%%%% Figure of 4 main bent diagonals %%%%%%%%%%%%%%%%%%%%%%%

Henceforth, when we say row sum, column sum, bent diagonal sum, and
so forth, we mean that we are adding the entries in the corresponding
configurations.  Franklin mentions that the square F1 has five other
curious properties but fails to list them. He also says, in the same
letter, that the $16 \times 16$ square F3 has all the properties of
the $8 \times 8$ square, but that in addition, every $4 \times 4$
subsquare adds to the common magic sum. More is true about this square
F3.  Observe that every $2 \times 2$ subsquare in F3 adds to
one-fourth the magic sum. The $8 \times 8$ squares have magic sum 260
while the $16 \times 16$ square has magic sum 2056. For a detailed
study of these three ``Franklin'' squares, see \cite{andrews},
\cite{pasles}, or \cite{pasles2}.

We define {\em $8 \times 8$ Franklin squares} to be squares with
nonnegative integer entries that have the properties (1) - (4) listed
by Benjamin Franklin and the additional property that every $2 \times
2$ subsquare adds to one-half the magic sum (see Figure \ref{8x8}).
The $8 \times 8$ squares constructed by Franklin have this extra
property (this might be one of the unstated curious properties to
which Franklin was alluding in his letter). It is worth noticing that
the fourth property listed by Benjamin Franklin becomes redundant with
the assumption of this additional property.

Similarly, we define {\em $16 \times 16$ Franklin squares} to be $16
\times 16$ squares that have nonnegative integer entries with the
property that all rows, columns, and bent diagonals add to the magic
sum, the half-rows and half-columns add to one-half the magic sum, and
the $2 \times 2$ subsquares add to one-fourth the magic sum. The $2
\times 2$ subsquare property implies that every $4 \times 4$ subsquare
adds to the common magic sum.

The property of the $2 \times 2$ subsquares adding to a common sum and
the property of bent diagonals adding to the magic sum are
``continuous properties.'' By this we mean that, if we imagine the
square as the surface of a torus (i.e., if we glue opposite sides of
the square together), then the bent diagonals and the $2 \times 2$
subsquares can be translated without effect on the corresponding sums
(see Figure \ref{parallelbents}).

%%%%%%% parallel bent diagonals %%%%%%%%%%%%%%%%%
%%%%%%%%%%%%%%%%%%%%%%%%%%%%%%%%%%%%%%%%%%%%%%%%

When the entries of an $n \times n$ Franklin square ($n=8$ or $n=16$)
are $1,2,3,\dots,n^2$, it is called a {\em natural Franklin square}.
Observe that the squares in Figure \ref{franklinsquares} are natural
Franklin squares.  Very few such squares are known, for the simple
reason that squares of this type are difficult to construct. Many
authors have looked at the problem of constructing them (see, for
example, \cite{andrews}, \cite{pasles}, and \cite{patel} and the
references therein). 

Our method is a new way to construct F1, F2, and F3. We are also able
to construct new natural Franklin squares, not isomorphic to the ones
previously known. A permutation of the entries of a Franklin square is
a {\em symmetry operation} if it maps the set of all Franklin squares
to itself, and two squares are called {\em isomorphic} if it is
possible to transform one to the other by applying symmetry
operations. We start by describing several symmetries of the
 Franklin squares.

\begin{thm} \label{symmetries}
The following operations on $n \times n$ Franklin squares, where $n=8$
or $n= 16$, are symmetry operations: rotating the square by 90
degrees; reflecting the square across one of its edges; transposing
the square; interchanging alternate columns (respectively, rows) $i$
and $i+2$, where $1 \leq i \leq (n/2) - 2$ or $(n/2)+1 \leq i \leq
n-2$; interchanging the first $n/2$ columns (respectively, rows) of
the square with the last $n/2$ columns (respectively, rows)
simultaneously; interchanging all the adjacent columns (respectively,
rows) $i$ and $i+1$ $(i= 1,3,5, \dots n-1)$ of the square
simultaneously.

The following operations are additional symmetry operations on $16
\times 16$ Franklin squares: interchanging columns (respectively,
rows) $i$ and $i+4$, where $1 \leq i \leq 4$ or $9 \leq i \leq 12$.
\end{thm}

The proof of Theorem \ref{symmetries} is presented in Chapter
\ref{franklinchapter}.

%%% constructing 8x8 squares %%%%%%%%
\begin{figure}[h]
 \begin{center}
     \includegraphics[scale=0.5]{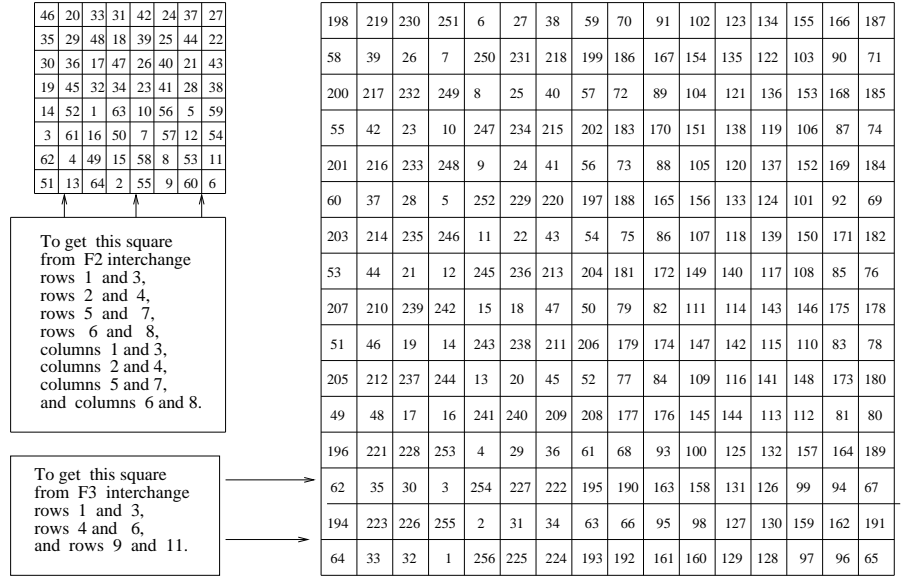}
\caption{Constructing Franklin squares by row and column exchanges of
Franklin squares.} \label{newexample} \end{center} \end{figure}
%%%%%%%%%%%%%%%%%%%%%%%%%%%%%%%%%%%%%%%%%%%%%%%%%%%%%%%%%%%%%%%%

%%% constructing new Franklin  squares %%%%%%%%
\begin{figure}[h]
 \begin{center}
     \includegraphics[scale=0.5]{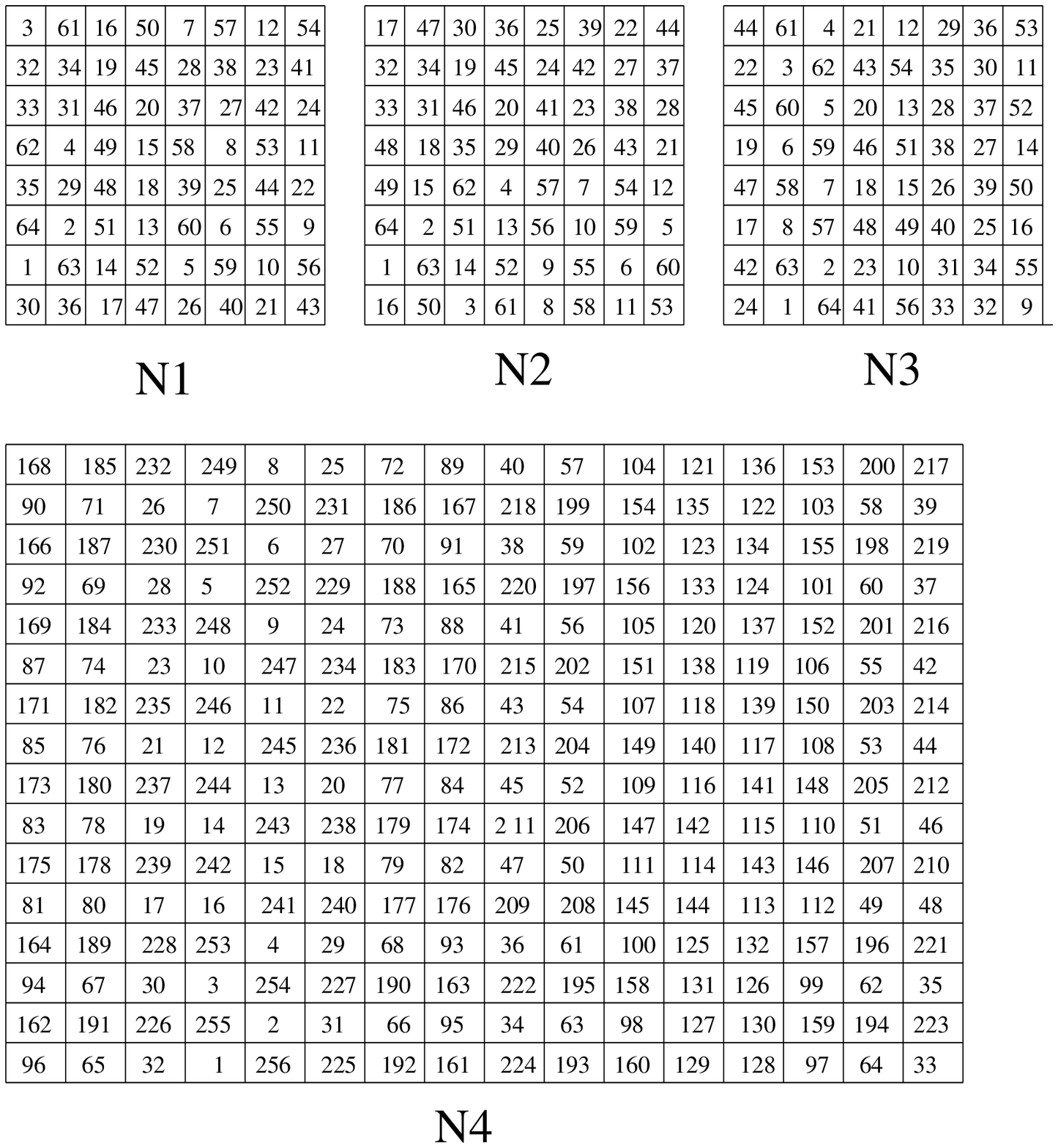}
\caption{New natural Franklin squares constructed using
Hilbert bases.} \label{newsquares} \end{center} \end{figure}
%%%%%%%%%%%%%%%%%%%%%%%%%%%%%%%%%%%%%%%%%%%%%%%%%%%%%%%%%%%%%%%%

The following theorem addresses itself to the squares in Figures 
\ref{franklinsquares} and \ref{newsquares}.

\begin{thm} \label{nonsymmetric}
The original Franklin squares F1 and F2 are not isomorphic.  The
squares N1, N2, and N3 are nonisomorphic natural $8 \times 8$ Franklin
squares that are not isomorphic to either F1 or F2.  The square N4 is
a natural $16 \times 16$ Franklin square that is not isomorphic to
square F3.
\end{thm}

We give a proof of Theorem \ref{nonsymmetric} in Chapter
\ref{franklinchapter}. We now enumerate Franklin squares with our
methods.

\begin{thm} \label{enumerate8x8}
Let $F_8(s)$ denote the number of $8 \times 8$ Franklin squares with
magic sum $s$. Then:

{ \footnotesize
\[
F_8(s)  = \left \{
\begin{array}{l}
\begin{array}{l}
{\frac {23}{627056640}}\,{s}^{9}+{\frac {23}{17418240}}\,{s}^{8}+{
\frac {167}{6531840}}\,{s}^{7}+{\frac {5}{15552}}\,{s}^{6} +{\frac {
2419}{933120}}\,{s}^{5}+{\frac {1013}{77760}}\,{s}^{4}+{\frac {701}{
22680}}\,{s}^{3} \\ \\  -{\frac {359}{10206}}\,{s}^{2}-{\frac {177967}{816480}
}\,s+{\frac {241}{17496}}
\end{array}  \\

\hfill \mbox{if s $\equiv 2$ (mod $12$) and $s \neq 2$, } \\ \\
 
\begin{array}{l}
{\frac {23}{627056640}}\,{s}^{9}+{\frac {23}{17418240}}\,{s}^{8}+{
\frac {167}{6531840}}\,{s}^{7}+{\frac {5}{15552}}\,{s}^{6} 
+{\frac {581 }{186624}}\,{s}^{5} +{\frac
{1823}{77760}}\,{s}^{4}+{\frac {6127}{45360 }}\,{s}^{3} \\ \\ +{\frac
{10741}{20412}}\,{s}^{2}+{\frac {113443}{102060}}\, s+{\frac
{3211}{2187}}
\end{array}
\\   \hfill \mbox{if $s \equiv 4$ (mod $12$),}  \\ \\

\begin{array}{l}
{\frac {23}{627056640}}\,{s}^{9}+{\frac {23}{17418240}}\,{s}^{8}+{
\frac {167}{6531840}}\,{s}^{7}+{\frac {5}{15552}}\,{s}^{6}
+{\frac { 2419}{933120}}\,{s}^{5} +{\frac
{1013}{77760}}\,{s}^{4}+{\frac {701}{ 22680}}\,{s}^{3} \\ \\ -{\frac
{5}{378}}\,{s}^{2}-{\frac {3967}{10080}}\,s-{ \frac {13}{8}}
\end{array}
 \\  \hfill \mbox{if $s \equiv 6$ (mod $12$),}  \\ \\

\begin{array}{l}
{\frac {23}{627056640}}\,{s}^{9}+{\frac {23}{17418240}}\,{s}^{8}+{
\frac {167}{6531840}}\,{s}^{7}+{\frac {5}{15552}}\,{s}^{6}+{\frac {581
}{186624}}\,{s}^{5}+{\frac {1823}{77760}}\,{s}^{4}+{\frac {6127}{45360
}}\,{s}^{3} \\ \\ +{\frac {11189}{20412}}\,{s}^{2}+{\frac {167203}{102060}}\,
s+{\frac {5771}{2187}}
\end{array}
 \\   \hfill \mbox{if $s \equiv 8$ (mod $12$),} \\ \\

\begin{array}{l}
{\frac {23}{627056640}}\,{s}^{9}+{\frac {23}{17418240}}\,{s}^{8}+{\frac {167}{6531840}}\,{s}^{7}+{\frac {5}{15552}}\,{s}^{6} +{\frac {
2419}{933120}}\,{s}^{5} +{\frac {1013}{77760}}\,{s}^{4}+{\frac {701}{
22680}}\,{s}^{3} \\ \\  -{\frac {583}{10206}}\,{s}^{2}-{\frac {608047}{816480}
}\,s-{\frac {20239}{17496}}
\end{array}
\\  \hfill   \mbox{if $s \equiv 10$ (mod $12$),}\\ \\

\begin{array}{l}
{\frac {23}{627056640}}\,{s}^{9}+{\frac {23}{17418240}}\,{s}^{8}+{
\frac {167}{6531840}}\,{s}^{7}+{\frac {5}{15552}}\,{s}^{6} +{\frac {581
}{186624}}\,{s}^{5} +{\frac {1823}{77760}}\,{s}^{4}+{\frac {6127}{45360
}}\,{s}^{3} \\ \\ +{\frac {431}{756}}\,{s}^{2}+{\frac {1843}{1260}}\,s+1
\end{array}
 \\  \hfill \mbox{if $s \equiv 0$ (mod $12$),} \\ \\
\begin{array}{l}
0 \end{array} \\ \hfill \mbox{otherwise.}
\end{array}
\right .
\]
}
\end{thm}

The {\em $8 \times 8$ pandiagonal Franklin squares} are the $8 \times
8$ Franklin squares that have all the pandiagonals adding to the
common magic sum (see Figure \ref{pans}). Our techniques enable us to
construct and count $8 \times 8$ pandiagonal Franklin squares as well:

\begin{thm} \label{enumeratepan8x8}
Let $PF_8(s)$ denote the number of $8 \times 8$ pandiagonal Franklin
squares with magic sum $s$. Then:

\[
PF_8(s)  = 
\left \{ \begin{array}{l}
{\frac {1}{2293760}}{s}^{8}+{\frac {1}{71680}}{s}^{7}+{\frac {1}{3840}}{s}^{6}+{\frac {1}{320}}{s}^{5}+{\frac{1}{40}}{s}^{4}+{\frac{2}{15}}{s}^{3}+{\frac {197}{420}}{s}^{2}+{\frac {106}{105}}s+1  \\
\hfill  \mbox{if  $s \equiv 0$ (mod $4$),} \\  
0 \\ \hfill  \mbox{otherwise.}
\end{array}
\right .
\]

\end{thm}
         
We introduce magic graphs and explore its connections to symmetric
magic squares in the next section.

\section{Magic labelings of graphs.}

\begin{figure}[h]
 \begin{center}
     \includegraphics[scale=0.5]{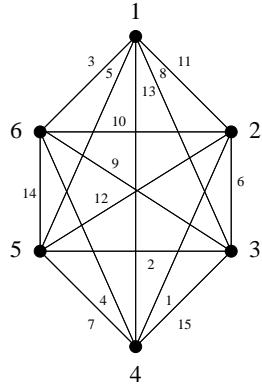}
\caption{A  magic labeling of the complete 
graph $K_6$ of magic sum $40$ \cite{stewartcomplete}.}  \label{mleg}
 \end{center}
 \end{figure}

Let $G$ be a finite graph. A {\em labeling} of $G$ is an assignment of
a nonnegative integer to each edge of $G$. A {\em magic labeling of
magic sum} $r$ of $G$ is a labeling such that for each vertex $v$ of
$G$ the sum of the labels of all edges incident to $v$ is the magic
sum $r$ (loops are counted as incident only once).  Graphs with a
magic labeling are also called {\em magic graphs} (see \ref{mleg} for
an example of a magic labeling of the complete graph $K_6$). Magic
graphs are also studied in great detail by Stanley and Stewart in
\cite{stanley3}, \cite{stanley4}, \cite{stewart}, and
\cite{stewartcomplete}.

%%%%complete6 decomp %%%%%%%%%%%%%%%%%%%%%%%

\begin{figure}[h]
 \begin{center}
\includegraphics[scale=0.5]{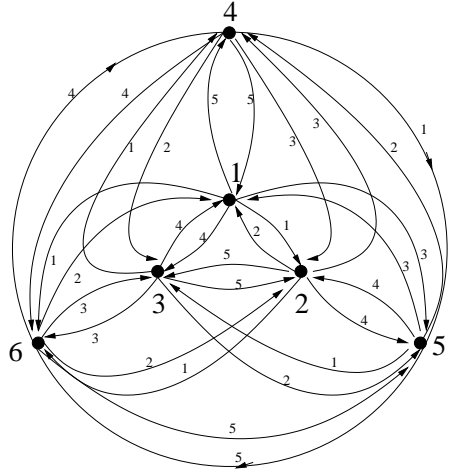}
\caption{Cayley digraph of the group $S_3$ \cite{konig}.}
 \label{cayleys3} \end{center} \end{figure}

We define a {\em magic labeling} of a digraph $D$ of {\em magic sum}
$r$ to be an assignment of a nonnegative integer to each edge of $D$,
such that for each vertex $v_i$ of $D$, the sum of the labels of all
edges with $v_i$ as the initial vertex is $r$, and the sum of the
labels of all edges with $v_i$ as the terminal vertex is also $r$.
Thus magic labelings of a digraph is a network flow, where the flow
into and out of every vertex, is the magic sum of the labeling (see
\cite{brualdi} for details about flows). Interesting examples of magic
digraphs are Cayley digraphs of finite groups. Let $G$ be a finite
group $\{g_1,g_2, \dots, g_n= I \}$. The {\em Cayley group digraph} of
$G$ is a graphical representation of $G$: every element $g_i$ of the
group $G$ corresponds to a vertex $v_i$ $(i=1,2, \dots, n)$ and every
pair of distinct vertices $v_i, v_j$ is joined by an edge labeled with
$\alpha$ where $g_{\alpha} = g_j g_i^{-1}$ (see \cite{hilton} or
\cite{konig}). For example, the Cayley digraph for the permutation
group
\[
S_3 = \{ g_1 = (123), g_2 = (132), g_3 = (23), g_4 = (12), g_5 = (13), g_6 = I \}
\]
 is given in Figure \ref{cayleys3}. 

\begin{prop} \label{cayley}
The Cayley digraph of a group of order $n$ is a magic digraph with
magic sum $\frac{n(n-1)}{2}$.
\end{prop}

\noindent {\em Proof.} Let $e_{ij}$ denote an edge between the vertex
$v_i$ and the vertex $v_j$ of the Cayley digraph such that $v_i$ is the initial
vertex and $v_j$ is the terminal vertex.  Let $v_l$ be a vertex of the
Cayley digraph, and let $\alpha$ be an integer in the set $\{ 1, 2,
\dots , n-1 \}$.  Let $g_p = g_{\alpha}g_l$ and let $g_q = g_l
g_{\alpha}$. Then, the edges $e_{lp}$ and $e_{ql}$ are labeled by
$\alpha$. Also, $g_jg_i^{-1} = g_n = I$ if and only if $i=j$.  Hence,
a Cayley group digraph is a magic digraph with magic sum $1 + 2 +
\cdots + (n-1) = \frac{n(n-1)}{2}$ (see also Chapter 8, Section 5 in
\cite{konig}). $\square$

%%%% Example %%%%%%%%%%%%%%%%%%%%%%%

A digraph is called {\em Eulerian} if for each vertex $v$ the indegree
and the outdegree of $v$ is the same. Therefore, Eulerian digraphs can
also be studied as magic digraphs where all the edges are labeled by 1
(see \cite{alantarsi} for the applications of Eulerian digraphs to
digraph colorings).

\begin{figure}[h]
 \begin{center}
     \includegraphics[scale=0.4]{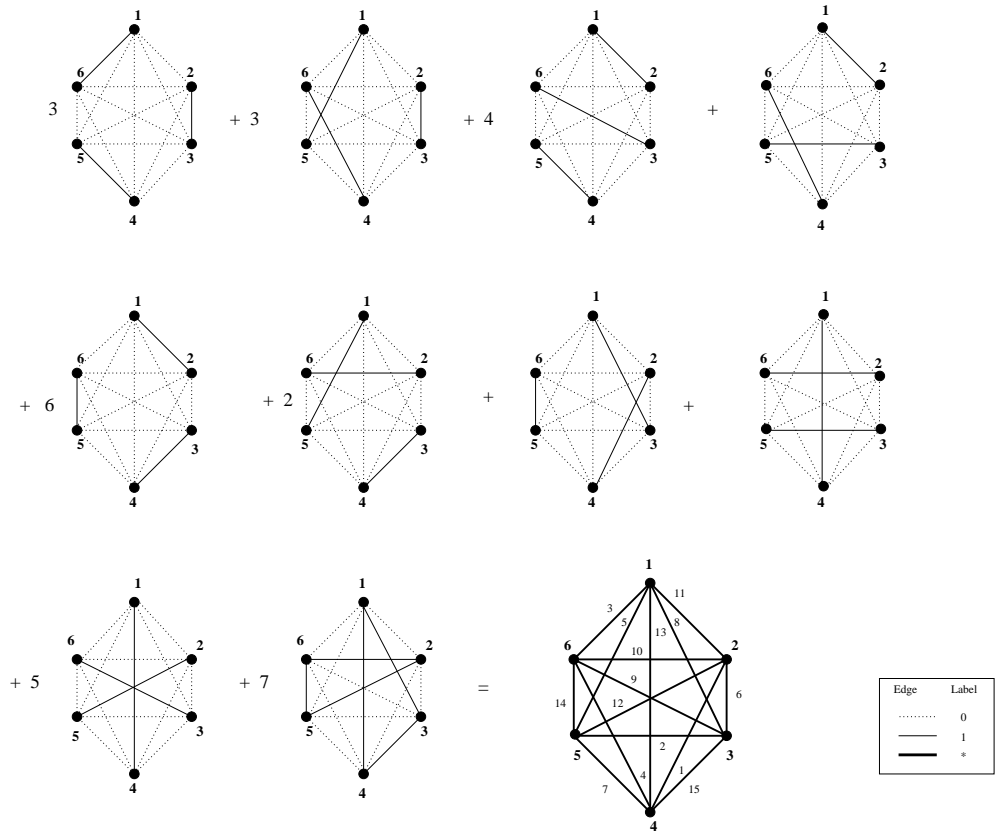}
\caption{A Hilbert basis construction of a magic labeling of the complete 
graph $K_6$.}  \label{hilbdecompeg}
 \end{center}
 \end{figure}

If we consider the labels of the edges of a magic graph $G$ as
variables, then the defining magic sum conditions are linear
equations. Thus, like before, the set of magic labelings of $G$
becomes the set of integral points inside a pointed polyhedral cone
$C_G$ (see Section \ref{method}). Therefore, with our methods we can
construct and enumerate magic labelings of graphs. For example, a
Hilbert basis construction of a magic labeling of the complete graph
$K_6$ is given in Figure \ref{hilbdecompeg}.

%%%%complete6 decomp %%%%%%%%%%%%%%%%%%%%%%%

Let $H_G(r)$ denote the number of magic labelings of $G$ of magic sum
$r$.  The generating functions of $H_G(r)$ in this thesis were
computed using the software LattE (see \cite{rudy}; software
implementation LattE is available from
http://www.math.ucdavis.edu/~latte). LattE was able to handle
computations that CoCoA was not able to perform. LattE uses Barvinok's
algorithm to compute generating functions which is different from
the methods used by CoCoA (see \cite{rudy}).

For example, let $\Gamma_n$ denote the {\em complete general graph} on
$n$ vertices (i.e, the complete graph with one loop at every
vertex). The formulas $H_{\Gamma_n}(r)$ for $n=3$ and $n=4$ were
computed by Carlitz \cite{carlitz}, and $H_{\Gamma_5}(r)$ was computed
by Stanley \cite{stanley4}. We use LattE to derive $H_{\Gamma_6}(r)$:

\[
H_{\Gamma_6}(r) = \left \{
\begin{array}{l}

{\frac {243653}{1992646656000}}{r}^{15}+{\frac {243653}{44281036800}
}{r}^{14}+{\frac {91173671}{797058662400}}{r}^{13}+{\frac {5954623
}{4087480320}}{r}^{12} \\ \\ +{\frac
{3895930519}{306561024000}}{r}^{11}+ {\frac
{21348281}{265420800}}{r}^{10}+{\frac {1063362673}{2786918400
}}{r}^{9}+{\frac {7132193}{5160960}}{r}^{8}+{\frac {479710409}{
124416000}}{r}^{7} \\ \\ +{\frac {963567863}{116121600}}{r}^{6}+{\frac {
26240714351}{1916006400}}{r}^{5} +{\frac
{39000163}{2280960}}{r}^{4 }+{\frac
{1514268697}{96096000}}{r}^{3} \\ \\ +{\frac {74169463}{7207200}}
{r}^{2}+{\frac {176711}{40040}}r+1 \\

\hfill  \mbox{if 2 divides $r$,} \\ \\  

{\frac {243653}{1992646656000}}{r}^{15}+{\frac {243653}{44281036800}
}{r}^{14}+{\frac {91173671}{797058662400}}{r}^{13}+{\frac {5954623
}{4087480320}}{r}^{12} \\ \\ +{\frac {3895930519}{306561024000}}{r}^{11} 
+
{\frac {21348281}{265420800}}{r}^{10}+{\frac {1063362673}{2786918400
}}{r}^{9}+{\frac {7132193}{5160960}}{r}^{8}+{\frac {479710409}{
124416000}}{r}^{7} \\ \\ +{\frac {963567863}{116121600}}{r}^{6}+{\frac {
839695842607}{61312204800}}{r}^{5} 
+{\frac {9983039353}{583925760}}
{r}^{4}+{\frac {774706849739}{49201152000}}{r}^{3} \\ \\ +{\frac {
302389338073}{29520691200}}{r}^{2}+{\frac {353330563}{82001920}}r+
{\frac {58885}{65536}} \\

\hfill \mbox{ otherwise.}
\end{array}
\right .
\]

An {\em $n$-matching} of $G$ is a magic labeling of $G$ with magic sum
at most $n$ and the labels are from the set $\{0, 1,\dots,n \}$ (see
\cite{lovasz}, chapter 6).  A {\em perfect matching} of $G$ is a
1-matching of $G$ with magic sum 1.

\begin{prop} \label{perfectmatchings} 
The perfect matchings of $G$ are the minimal Hilbert basis elements of
$C_G$ of magic sum 1 and the number of perfect matchings of $G$ is
$H_G(1)$. 
\end{prop}

\noindent{\em Proof.}  Magic labelings of magic sum 1 always belong to
the minimal Hilbert basis because they are irreducible.  Therefore,
perfect matchings belong to the minimal Hilbert basis because they
have magic sum 1. Conversely, every magic labeling of magic sum 1 is a
perfect matching. So we conclude that the perfect matchings of $G$ are
the minimal Hilbert basis elements of $C_G$ of magic sum 1. The fact
that the number of perfect matchings of $G$ is $H_G(1)$ follows by the
definition of $H_G(1)$.  $\square$

The perfect matchings of $G$ can be found by computing a
truncated Hilbert basis of magic sum 1 using 4ti2 (see \cite{raymond}).

\begin{figure}[h]
 \begin{center} 
 \includegraphics[scale=0.4]{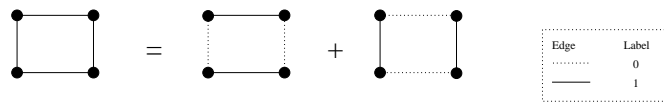}
\caption{Graph Factorization.} \label{factor}
\end{center}
\end{figure}

Hilbert bases can also be used to study factorizations of labeled
graphs. We define {\em Factors} of a graph $G$ with a labeling $L$ to
be labelings $L_i, i = 1, \dots, r$ of $G$ such that $L(G) =
\sum_{i=1}^r L_i(G)$, and if $L_i(e_k) \neq 0$ for some edge $e_k$ of
$G$, then $L_j (e_k) = 0$ for all $ j \neq i$. A decomposition of $L$
into factors is called a {\em factorization} of $G$.  An example of a
graph factorization is given in Figure \ref{factor}. See Chapters 11
and 12 of \cite{konig} for a detailed study of graph factorizations.

%%%% Example %%%%%%%%%%%%%%%%%%%%%%%

With our methods we can also construct and enumerate magic labelings
of digraphs. Let $H_D(r)$ denote the number of magic labelings of a
digraph $D$ of magic sum $r$. We now connect the magic labelings of
digraphs to magic labelings of bipartite graphs.

\begin{lemma} For every digraph $D$, there is a bipartite graph $G_D$
such that the magic labelings of $D$ are in one-to-one correspondence
with the magic labelings of $G_D$. Moreover, the magic sums of the
corresponding magic labelings of $D$ and $G_D$ are also the
same. \end{lemma}

{\em Proof.}  Denote a directed edge of a digraph $D$ with $v_i$ as
the initial vertex, and $v_j$ as the terminal vertex, by $e_{ij}$. Let
$L$ be a magic labeling of $D$ of magic sum $r$. Consider a bipartite
graph $G_D$ in $2n$ vertices, where the vertices are partitioned into
two sets $A = \{a_1, \dots, a_n \}$ and $B = \{b_1, \dots, b_n \}$,
such that there is an edge between $a_i$ and $b_j$, if and only if,
there is an edge $e_{ij}$ in $D$. Consider a labeling $L_{G_D}$ of
$G_D$ such that the edge between the vertices $a_i$ and $b_j$ is
labeled with $L(e_{ij})$.  Observe that the sum of the labels of the
edges incident to $a_{i}$ is the same as the sum of the labels of
incoming edges at the vertex $v_i$ of $D$. Also, the sum of the labels
of edges at a vertex $b_j$ is the sum of the labels of outgoing edges
at the vertex $v_j$ of $D$. Since $L$ is a magic labeling, it follows
that $L_{G_D}$ is a magic labeling of $G_D$ with magic sum $r$.  Going
back-wards, consider a magic labeling $L^{\prime} $ of $G_D$. We label
every edge $e_{ij}$ of $D$ with the label of the edge between $a_i$
and $b_j$ of $G_D$ to get a magic labeling $L_D$ of $D$. Observe that
 $L^{\prime} $ and $L_D$ have the same magic sum. Hence, there is a
one-to-one correspondence between the magic labelings of $D$ and the
magic labelings of $G_D$. $\square$

\begin{figure}[h]
 \begin{center} 
 \includegraphics[scale=0.5]{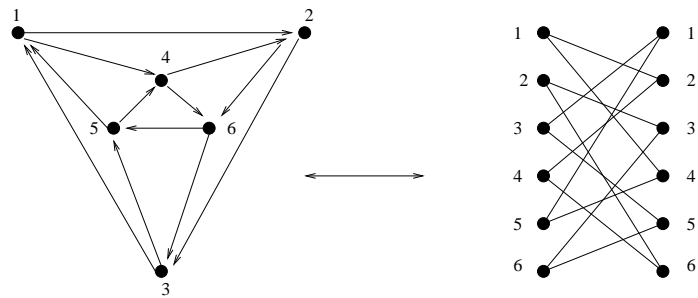}
\caption{The octahedral digraph with a given orientation $D_O$ and its
corresponding bipartite graph $G_{D_O}$.} \label{octexample}
\end{center} \end{figure}

For example, the magic labelings of the Octahedral digraph with the
given orientation $D_O$ in Figure \ref{octexample} are in one-to-one
correspondence with the magic labelings of the bipartite graph
$G_{D_O}$.

%%%% Example %%%%%%%%%%%%%%%%%%%%%%%

Similarly, for a bipartite graph $B$ we can get a digraph $B_D$ such
that the magic labelings of $B$ are in one-to-one correspondence with
the magic labelings of $B_D$: let $B$ be such that the vertices are
partitioned into sets $A = \{a_1, \dots, a_n \}$ and $B = \{b_1,
\dots, b_m \}$. Without loss of generality assume $n> m$.  Then $B_D$
is the digraph with $n$ vertices such that there is an edge $e_{ij}$
in $B_D$ if and only if there is an edge between the vertices $a_i$
and $b_j$ in $B$. This correspondence enables us to generate and
enumerate perfect matchings of bipartite graphs.

\begin{prop} \label{digraphperfectmatchings}
There is a one-to-one correspondence between the perfect matchings of
a bipartite graph $B$ and the elements of the Hilbert basis of
$C_{B_D}$. The number of perfect matchings of $B$ is $H_{B_D}(1)$.
\end{prop}

We present a proof of Proposition \ref{digraphperfectmatchings} in
Chapter \ref{magicgraphchapter}. A graph $G$ is called a {\em positive
graph} if for any edge $e$ of $G$ there is a magic labeling $L$ of $G$
for which $L(e) > 0$ \cite{stanley3}.  Since edges of $G$ that are
always labeled zero for any magic labeling of $G$ may be ignored to
study magic labelings, we will concentrate on positive graphs in
general. We use the following results by Stanley from \cite{stanley3}
and \cite{stanley4} to prove Theorems \ref{graphpolytopethm} and
\ref{digraphpolytopethm} and Corollary \ref{digraphenumerate}
\begin{thm} [Theorem 1.1, \cite{stanley4}] \label{quasipoly}
Let $G$ be a finite positive graph. Then either $H_G(r)$ is the
Kronecker delta $\delta_{0r}$ or else there exist polynomials $I_G(r)$
and $J_G(r)$ such that $H_G(r) = I_G(r) + (-1)^r J_G(r)$ for all $r
\in \naturals$.
\end{thm}

\begin{thm} [Theorem 1.2, \cite{stanley4}]  \label{degreeofG}
Let $G$ be a finite positive graph with at least one edge. The degree
of $H_G(r)$ is $q-n+b$, where $q$ is the number of edges of $G$, $n$
is the number of vertices, and $b$ is the number of connected
components of $G$ which are bipartite.
\end{thm}

\begin{thm}[Theorem 1.2, \cite{stanley3}] \label{bippoly}
Let $G$ be a finite positive bipartite graph with at least one edge,
then $H_G(r)$ is a polynomial.
\end{thm}

We now conclude that $H_{D}(r)$ is a polynomial for every digraph $D$.

\begin{corollary}  \label{digraphenumerate}
If $D$ is a digraph, then $H_{D}(r)$ is a polynomial of degree
$q-2n+b$, where $q$ is the number of edges of $D$, $n$ is the number
of vertices, and $b$ is the number of connected components of the
bipartite graph $G_D$.
\end{corollary}
\noindent {\em Proof.} The one-to one correspondence between the magic
labelings of $D$ and the magic labelings of $G_D$, implies by Theorem
\ref{bippoly} that $H_D(r)$ is a polynomial, and by Theorem
\ref{degreeofG} that the degree of $H_D(r)$ is $q-2n+b$, where $b$ is
the number of connected components of $G_D$ that are
bipartite. $\square$

Consider the polytope ${\cal P} := \{x|Ax \leq b \}$.  Let $c$ be a
nonzero vector, and let $\delta = $ max $\{cx|Ax \leq b\}$.  The
affine hyperplane $\{x|cx = \delta \}$ is called a {\em supporting
hyperplane} of $\cal P$. A subset $F$ of $\cal P$ is called a {\em
face} of $\cal P$ if $F = \cal P$ or if $F$ is the intersection of
$\cal P$ with a supporting hyperplane of $\cal P$. Alternatively, 
$F$ is a face of $\cal P$ if and only if $F$ is nonempty and 
\[
F = \{ x \in {\cal P} | A^{\prime} x = b^{\prime} \}
\]
for some subsystem $ A^{\prime} x \leq b^{\prime}$ of $Ax \leq b$. See
\cite{schrijver} for basic definitions with regards to polytopes.

Let $v_1, v_2, \dots, v_n$ denote the vertices of a graph $G$ and let
$e_{i_1}, e_{i_2}, \dots, e_{i_{m_i}}$ denote the edges of $G$ that
are incident to the vertex $v_i$ of $G$.  Consider the polytope
\[
{\cal P}_G = \{ L \in C_G \subseteq \reals^{q}, \hspace{0.05 in}
\sum_{j=1}^{m_i}   L(e_{i_j}) = 1 ;  i=1, \dots, n
\}.
\]

We will refer to ${\cal P}_G$ as the polytope of magic labelings of
$G$. Then, $H_G(r)$ is the Ehrhart quasi-polynomial of ${\cal P}_G$
(see Section \ref{method}). A face of ${\cal P}_G$ is a polytope of
the form
\[
 \{ L \in {\cal P}_G,   L(e_{i_k}) = 0 ;  e_{i_k} \in E_0 \}, 
\]
where $E_0 = \{ e_{i_1}, \dots, e_{i_r} \}$ is a subset of the set of
edges of $G$.

\begin{thm} \label{graphpolytopethm}
Let $G$ be a finite positive graph with at least one edge. Then the
polytope of magic labelings of $G$, ${\cal P}_G$ is a rational
polytope with dimension $q-n+b$, where $q$ is the number of edges of
$G$, $n$ is the number of vertices, and $b$ is the number of connected
components of $G$ that are bipartite. The $d$-dimensional faces of
${\cal P}_G$ are the $d$-dimensional polytopes of magic labelings of
positive subgraphs of $G$ with $n$ vertices and at most $n-b+d$ edges.
\end{thm}

\begin{figure}[h]
 \begin{center}
     \includegraphics[scale=0.4]{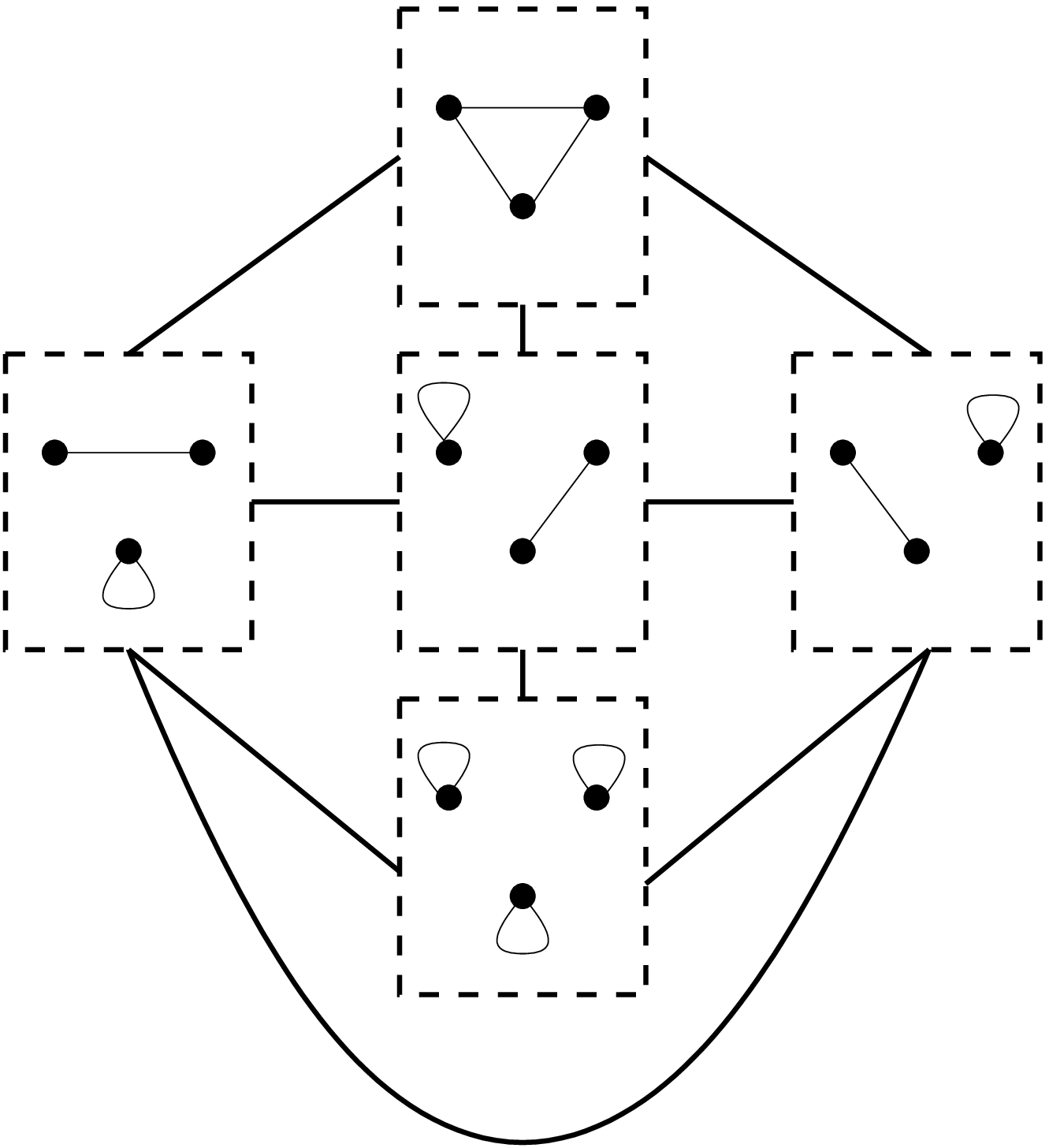}
\caption{The edge graph of  ${\cal P}_{\Gamma_3}$.}  \label{edgegraph}
 \end{center}
 \end{figure}

\begin{figure}[h]
 \begin{center}
     \includegraphics[scale=0.3]{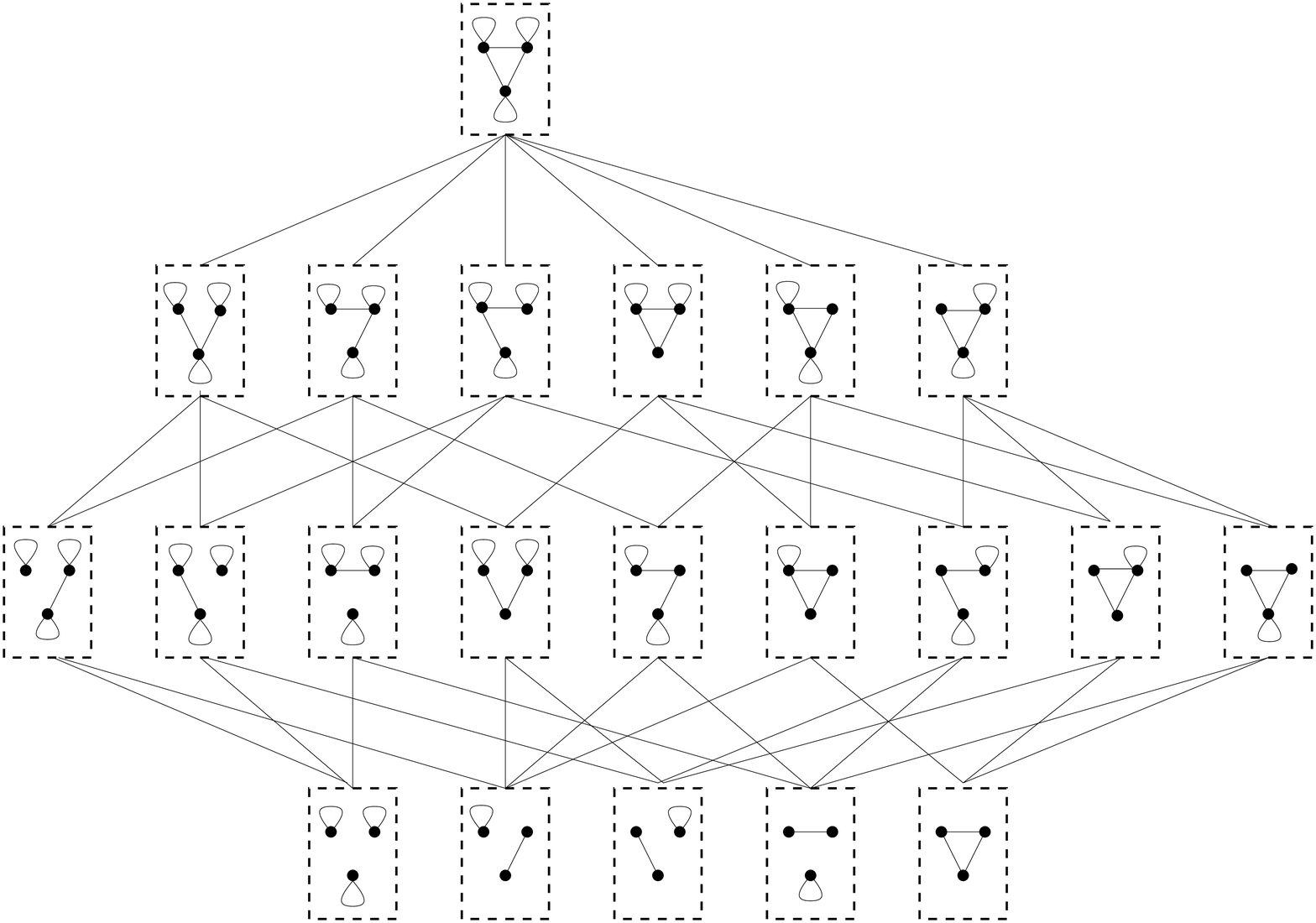}
\caption{The face poset of  ${\cal P}_{\Gamma_3}$.}  \label{faceposet}
 \end{center}
 \end{figure}

We prove Theorem \ref{graphpolytopethm} in Chapter
\ref{magicgraphchapter}.  Observe from Theorem \ref{graphpolytopethm}
that there is an edge between two vertices $v_i$ and $v_j$ of ${\cal
P}_G$ if and only if there is a graph with at most $n-b+1$ edges, with
magic labelings $v_i$ and $v_j$.  The edge graph of ${\cal
P}_{\Gamma_3}$ is given in Figure \ref{edgegraph}. Similarly, we can
draw the face poset of ${\cal P}_{G}$ (see Figure \ref{faceposet} for
the face poset of ${\cal P}_{\Gamma_3}$).

%%%% Figure edge graph %%%%%%%%%%%%%%%%%%%%%%%

%%%% Figure face poset %%%%%%%%%%%%%%%%%%%%%%%

An $n \times n$ {\em semi-magic square} of magic sum $r$ is an $n
\times n$ matrix with nonnegative integer entries such that the
entries of every row and column add to $r$. {\em Doubly stochastic
matrices} are $n \times n$ matrices in $\reals^{n^2}$ such that their
rows and columns add to 1. The set of all $n \times n$ doubly
stochastic matrices form a polytope $B_n$, called the {\em Birkhoff
polytope}. See \cite{billerasarang}, \cite{brualdigibson}, or
\cite{schrijver} for a detailed study of the Birkhoff polytope.

A {\em symmetric magic square} is a semi-magic square that is also a
 symmetric matrix. Let $H_n(r)$ denote the number of symmetric magic
 squares of magic sum $r$ (see \cite{carlitz}, \cite{gupta}, and
 \cite{stanley4} for the enumeration of symmetric magic squares).  We
 define the polytope ${\cal S}_n$ of $n \times n$ symmetric magic
 squares to be the convex hull of all real nonnegative $n \times n$
 symmetric matrices such that the entries of each row (and therefore
 column) add to one.

\begin{figure}[h]
 \begin{center}
     \includegraphics[scale=0.5]{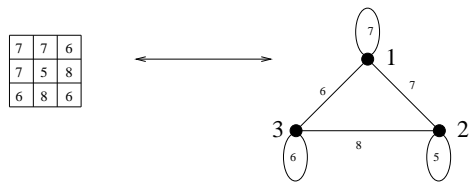}
\caption{A magic labeling of $\Gamma_3$ and its corresponding
symmetric magic square.} \label{gamma3eg}
 \end{center}
 \end{figure}

 A one-to-one correspondence between symmetric magic squares $M =
 [m_{ij}]$ of magic sum $r$, and magic labelings of the graph
 $\Gamma_n$ of the same magic sum $r$ was established in
 \cite{stanley4}: let $e_{ij}$ denote an edge between the vertex $v_i$
 and the vertex $v_j$ of $\Gamma_n$. Label the edge $e_{ij}$ of
 $\Gamma_n$ with $m_{ij}$, then this labeling is a magic labeling of
 $\Gamma_n$ with magic sum $r$. See Figure \ref{gamma3eg} for an
 example. Therefore, we get ${\cal P}_{\Gamma_n}$ = ${\cal S}_{n}$ and
 $H_{\Gamma_n}(r) = H_n(r)$.

%%%% Figure of symmetric magic correspondence %%%%%%%%%%%%%%%%%%%%%%%

\begin{corollary} \label{mayathm}
The polytope of magic labelings of the complete general graph ${\cal
P}_{\Gamma_n}$ is an $n(n-1)/2$ dimensional rational polytope with the
following description
\[\begin{array}{llll}
{\cal P}_{\Gamma_n} &=& \{ L = (L(e_{ij}) \in {\reals}^{\frac{n(n+1)}{2}}; 
&  L(e_{ij}) \geq 0; 1 \leq i,j \leq n, i \leq j,   \\ 
&&& \sum_{j=1}^i L(e_{ji}) + \sum_{j=i+1}^n L(e_{ij}) = 1 
\mbox{ for } i=1, \dots, n 
\}.
\end{array}
\]
The $d$-dimensional faces of ${\cal P}_{\Gamma_n}$ are $d$-dimensional
polytopes of magic labelings of positive graphs with $n$ vertices and
at most $n+d$ edges. There are $2n-1 \choose n$ faces of ${\cal
P}_{\Gamma_{2n}}$ that are copies of the Birkhoff polytope $B_n$.
\end{corollary} 

See Chapter \ref{magicgraphchapter} for the proof of Corollary
\ref{mayathm}.  Again, as in the case of graphs, we define a polytope
${\cal P}_D$ of magic labelings of $D$: Let $e_{i_1}, e_{i_2}, \dots,
e_{i_{m_i}}$ denote the edges of $D$ that have the vertex $v_i$ as the
initial vertex and let $f_{i_1}, f_{i_2}, \dots, f_{i_{s_i}}$ denote
the edges of $D$ for which the vertex $v_i$ is the terminal vertex,
then:

\[
{\cal P}_D = \{ L \in C_D \subseteq \reals^{q}, \hspace{0.05 in}
\sum_{j=1}^{m_i} L(e_{i_j}) = \sum_{j=1}^{s_i} L(f_{i_j}) = 1 ; i=1,
\dots, n \}.
\]

We define a digraph $D$ to be a {\em positive digraph} if the corresponding 
bipartite graph $G_D$ is positive. 

\begin{thm} \label{digraphpolytopethm}
Let $D$ be a positive digraph with at least one edge. Then, ${\cal
P}_D$ is an integral polytope with dimension $q-2n+b$, where $q$ is
the number of edges of $D$, $n$ is the number of vertices, and $b$ is
the number of connected components of $G_D$ that are bipartite. The
$d$-dimensional faces of ${\cal P}_D$ are the $d$-dimensional
polytopes of magic labelings of positive subdigraphs of $D$ with $n$
vertices and at most $2n-b+d$ edges.
\end{thm}

See Chapter \ref{magicgraphchapter} for the proof of Theorem
\ref{digraphpolytopethm}.

\begin{figure}[h]
 \begin{center} 
 \includegraphics[scale=0.5]{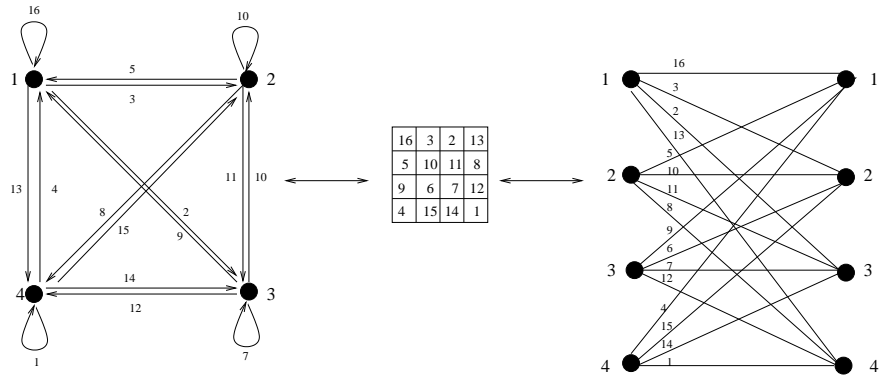}
\caption{Two different graph labelings associated to a semi-magic square.}
 \label{semi-magic} \end{center} \end{figure}

Let $\Pi_n$ denote the complete digraph with $n$ vertices, i.e, there
is an edge from each vertex to every other, including the vertex
itself (thereby creating a loop at every vertex), then $G_{\Pi_n}$ is
the the complete bipartite graph $K_{n,n}$. We get a one-to-one
correspondence between semi-magic squares $M = [m_{ij}]$ of magic sum
$r$ and magic labelings of $\Pi_n$ of the same magic sum $r$ by
labeling the edges $e_{ij}$ of $\Pi_n$ with $m_{ij}$.  This also
implies that there is a one-to-one correspondence between semi-magic
squares and magic labelings of $K_{n,n}$ (this correspondence is also
mentioned in \cite{stanley3} and \cite{stewart}).  See Figure
\ref{semi-magic} for an example.

%%%% Example %%%%%%%%%%%%%%%%%%%%%%%

A good description of the faces of Birkhoff polytope is not known \cite{pak}.
We can now give an explicit description of the faces of the Birkhoff polytope.

\begin{thm} \label{facesbirk}
${\cal P}_{\Pi_n}$ is the Birkhoff polytope $B_n$.  The
$d$-dimensional faces of $B_n$ are $d$-dimensional polytopes of magic
labelings of positive digraphs with $n$ vertices and at most $2n+d-1$
edges. The vertices of ${\cal P}_D$, where $D$ is a positive digraph,
are permutation matrices.
 
\end{thm}

\begin{figure}[h]
 \begin{center}
     \includegraphics[scale=0.5]{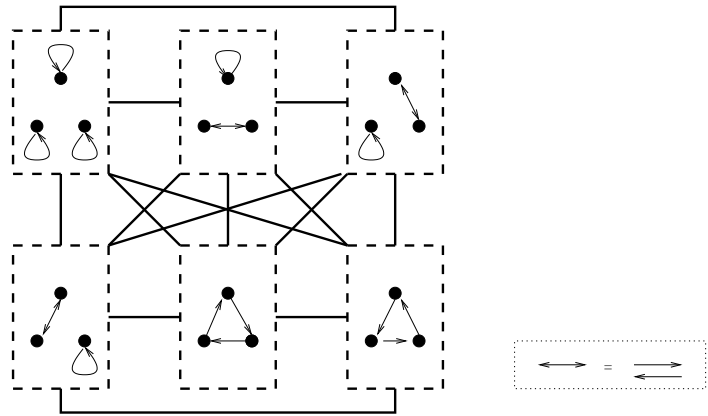}
\caption{The edge graph of the Birkhoff Polytope $B_3$.}  \label{b3edgegraph}
 \end{center}
 \end{figure}

\begin{figure}[h]
 \begin{center}
     \includegraphics[scale=0.5]{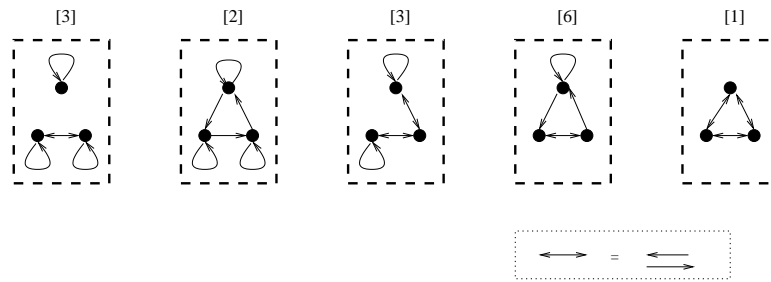}
\caption{The generators of the edges of the Birkhoff
Polytope $B_3$.}  \label{b3edges}
 \end{center}
 \end{figure}

\begin{figure}[h]
 \begin{center}
     \includegraphics[scale=0.5]{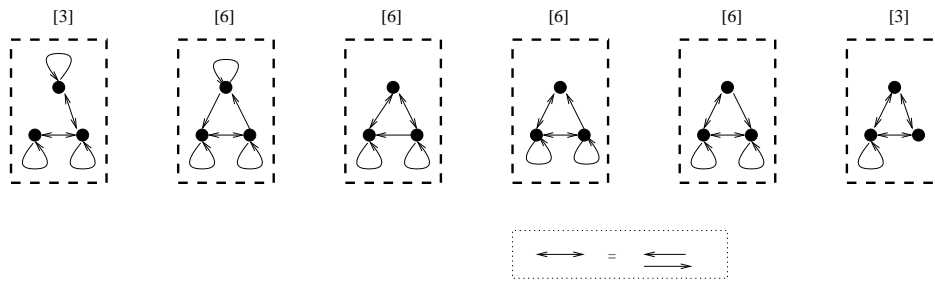}
\caption{The generators of the 2-dimensional faces of the Birkhoff
Polytope $B_3$.}  \label{b32dimfaces}
 \end{center}
 \end{figure}

\begin{figure}[h]
 \begin{center}
     \includegraphics[scale=0.5]{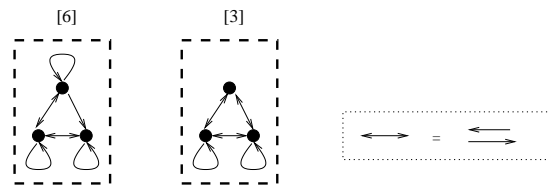}
\caption{The generators of the facets of the Birkhoff
Polytope $B_3$.}  \label{b3facets}
 \end{center}
 \end{figure}

\begin{figure}[h]
 \begin{center}
     \includegraphics[scale=0.5]{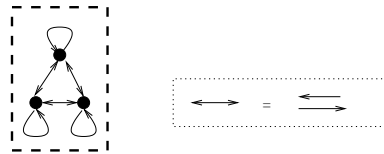}
\caption{The Birkhoff Polytope $B_3$.}  \label{b3}
 \end{center}
 \end{figure}

The proof of Theorem \ref{facesbirk} is presented in Chapter
 \ref{magicgraphchapter}. See Figure \ref{b3edgegraph} for the edge
 graph of $B_3$. Two faces of a polytope of magic labelings of a graph
 (or a digraph) are said to be {\em isomorphic faces} if the subgraphs
 (subdigraphs, respectively) defining the faces are isomorphic. A set
 of faces is said to be a {\em generating set of $d$-dimensional
 faces} if every $d$-dimensional face is isomorphic to one of the
 faces in the set. See Figures \ref{b3edges}, \ref{b32dimfaces},
 \ref{b3facets}, and \ref{b3} for the generators of the edges, the two
 dimensional faces, the facets, and the Birkhoff polytope $B_3$,
 respectively (the numbers in the square brackets indicate the number
 of faces in the isomorphism class of the given face).

%%%% Figure edge graph %%%%%%%%%%%%%%%%%%%%%%%

%%%% Figure edges  %%%%%%%%%%%%%%%%%%%%%%%

%%%% Figure 2 dim faces %%%%%%%%%%%%%%%%%%%%%%%

%%%% Figure  facets %%%%%%%%%%%%%%%%%%%%%%%

%%%% Figure  facets %%%%%%%%%%%%%%%%%%%%%%%

\newpage
\pagestyle{myheadings} 
\markright{  \rm \normalsize CHAPTER 2. \hspace{0.5cm}
 Magic Cubes }
\large 
\chapter{Magic Cubes} \label{magiccubechapter}
\thispagestyle{myheadings}
{\small
\begin{verse}
They flash upon that inward eye \\
Which is the bliss of solitude; \\
And then my heart with pleasure fills, \\
And dances with the daffodils.

-- William Wordsworth.
\end{verse}
}

\begin{figure}[hpt]
 \begin{center}
     \includegraphics[scale=.4]{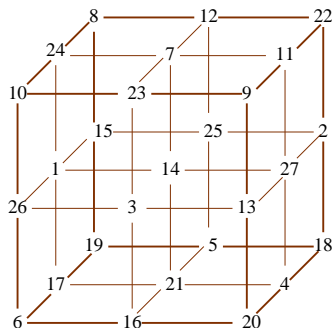}
 \caption{A magic cube.} \label{cubepicture}
 \end{center}
 \end{figure}

A {\em semi-magic hypercube} is a $d$-dimensional $n\times n \times
\dots \times n$ array of $n^d$ non-negative integers, which sum up to
the same number $s$ for any line parallel to some axis. A {\em magic
hypercube} is a semi-magic cube that has the additional property that
the sums of all the main diagonals, the $2^{d-1}$ copies of the
diagonal $x_{1,1,\dots,1},x_{2,2,\dots,2},\dots,x_{n,n,\dots,n}$ under
the symmetries of the $d$-cube, are also equal to the magic sum. For
example, in a $2 \times 2 \times 2$ cube there are 4 diagonals with
sums
$x_{1,1,1}+x_{2,2,2}=x_{2,1,1}+x_{1,2,2}=x_{1,1,2}+x_{2,2,1}=x_{1,2,1}+x_{2,1,2}.$

\begin{figure}[hpt]
 \begin{center}
     \includegraphics[scale=0.3]{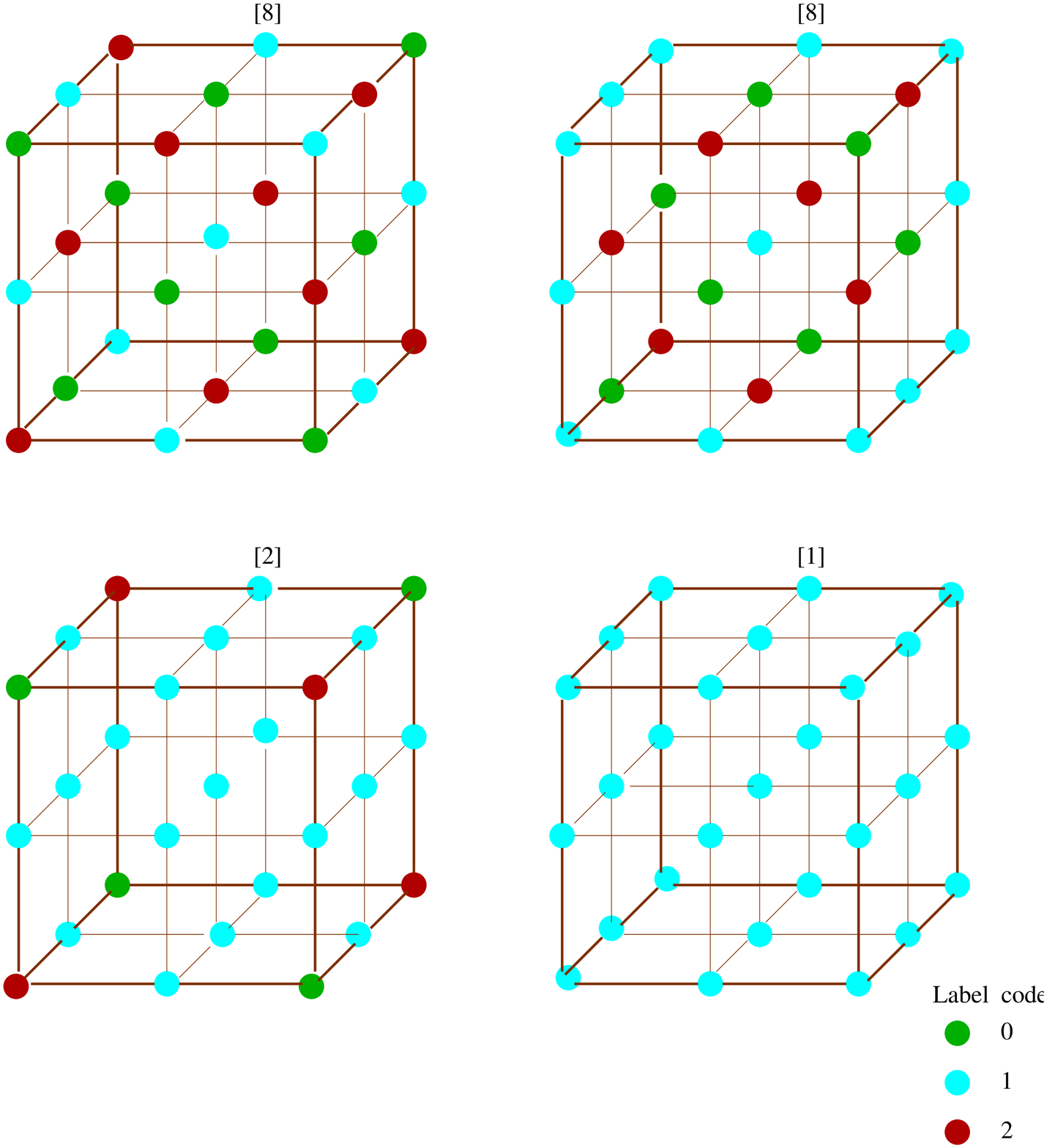}
     \caption{The generators of the Hilbert basis of the $3 \times 3
     \times 3$ magic cube.}
\label{cubehilbertbasis}
 \end{center}
 \end{figure}

An example of a  $3 \times 3 \times 3$ magic cube is given  in Figure
\ref{cubepicture}. If we consider the entries of a magic cube to be variables,
 the defining equations form a linear system of equations and thus
magic cubes are integral points inside a pointed polyhedral
cone. Therefore, we can use the methods described in Section
\ref{method} to construct and enumerate magic cubes. Let $G$ denote
the group of rotations of a cube \cite{dummit}. Two cubes are called
{\em isomorphic} if we can get one from the other by using a series of
rotations. A set of magic cubes are called {\em generators} of the
Hilbert basis if every element of the Hilbert basis is isomorphic to
one of the cubes in the set. The generators of the Hilbert basis of $3
\times 3 \times 3$ magic cubes are given in Figure
\ref{cubehilbertbasis} (the numbers in square brackets indicate the
number of elements in the orbit of a generator under the action of
$G$).  There are 19 elements in the Hilbert basis and all of them have
magic sum value of 3. An example of a Hilbert basis construction of a
magic cube is given in Figure
\ref{magiccubeconstruct}.

Let $MC_n(s)$ denote the number of $n \times n \times n$ magic cubes
of magic sum $s$. We use the algorithm presented in Section
\ref{method} to compute the generating function for the number of $3
\times 3 \times 3$ magic cubes:
\[
\begin{array}{l}
\sum_{s=0}^{\infty} MC_3(s)t^s={\frac
{{t}^{12}+14\,{t}^{9}+36\,{t}^{6}+14\,{t}^{3}+1}{\left (1-{t}^{
3}\right )^{5}}} \\ \\
=1+19\,{t}^{3}+121\,{t}^{6}+439\,{t}^{9}+1171\,{t}^{12}+2581\,{t}^{15}
+4999\,{t}^{18}+\dots 
\end{array}
\]

and we derive

\begin{thm}
The number of $3 \times 3 \times 3$ magic cubes
\[
MC_3(s) = \left\{
\begin{array}{ll}
\frac{11}{324}s^4+\frac{11}{54}s^3+\frac{25}{36}s^2+\frac{7}{6}s+1 
& \mbox {if $3 \vert s$,} \\ 
   0 & \mbox{otherwise.} 
\end{array}
\right.
\]
\end{thm}

\begin{thm} \label{thmcubomagia}
The number of $n \times n \times n$ magic cubes of magic sum $s$,
$MC_n(s)$ is a quasipolynomial of degree $(n-1)^3-4$ for $n \geq 3, n
\neq 4$. For $n=4$ it has degree $(4-1)^3-3=24$.
\end{thm}

\noindent {\em Proof.} 
The function that counts magic cubes is a quasipolynomial whose degree
is the same as the dimension of the cone of magic cubes minus one. For
small values (e.g $n=3,4$) we can directly compute this. We present an
argument for its value for $n>4$. Let $B$ be the $(3n^2+4) \times n^3$
matrix with $0,1$ entries determining axial and diagonal sums. In this
way we see that $n \times n \times n$ magic cubes of magic sum $s$ are
the integer solutions of $Bx=(s,s,\dots,s)^T, x \geq 0$. 

It is known that for semi-magic cubes the dimension is
$(n-1)^3$ \cite{becketal}, which means that the rank of the submatrix $B'$ of $B$
without the 4 rows that state diagonal sums is $n^3-(n-1)^3$.  It
remains to be shown that the addition of the $4$ sum constraints on
the main diagonals to the defining equations of the $n\times n\times
n$ semi-magic cube increases the rank of the defining matrix $B$ by
exactly $4$.

Let us denote the $n^3$ entries of the cube by
$x_{1,1,1},\ldots,x_{n,n,n}$ and consider the $(n-1)\times (n-1)\times
(n-1)$ sub-cube with entries $x_{1,1,1},\ldots,x_{n-1,n-1,n-1}$. For a
semi-magic cube we have complete freedom to choose these $(n-1)^3$
entries. The remaining entries of the $n\times n\times n$ magic cube
become known via the semi-magic cube equations, and all entries
together form a semi-magic cube. For example:

$x_{n,1,1}  =  -\sum_{i=1}^{n-1} x_{i,1,1},\
x_{1,n,n}  =  \sum_{i=1}^{n-1}\sum_{j=1}^{n-1} x_{i,j,1},\
x_{n,n,n}  =  -\sum_{i=1}^{n-1}\sum_{j=1}^{n-1}\sum_{k=1}^{n-1}
x_{i,j,k}.$

%\begin{eqnarray*}
%x_{n,1,1} & = & -\sum_{i=1}^{n-1} x_{i,1,1},\\
%x_{1,n,n} & = & \sum_{i=1}^{n-1}\sum_{j=1}^{n-1} x_{i,j,1},\\
%x_{n,n,n} & = & -\sum_{i=1}^{n-1}\sum_{j=1}^{n-1}\sum_{k=1}^{n-1}
%x_{i,j,k}.
%\end{eqnarray*}

However, for the magic cube, $4$ more conditions have to be
satisfied along the main diagonals. Employing the above semi-magic
cube equations, we can rewrite these $4$ equations for the main
diagonals such that they involve only the variables
$x_{1,1,1},\ldots,x_{n-1,n-1,n-1}$. Thus, as we will see, the
complete freedom of choosing values for the variables
$x_{1,1,1},\ldots,x_{n-1,n-1,n-1}$ is restricted by $4$
independent equations. Therefore the dimension of the kernel of $B$ is
reduced by $4$.

Let us consider the $3$ equations in
$x_{1,1,1},\ldots,x_{n-1,n-1,n-1}$ corresponding to the main diagonals
$x_{1,1,n},\ldots,x_{n,n,1}$, $x_{1,n,1},\ldots,x_{n,1,n}$, and
$x_{n,1,1},\ldots,x_{1,n,n}$.  They are linearly independent, since
the variables $x_{n-1,n-1,1}$, $x_{n-1,1,n-1}$, and $x_{1,n-1,n-1}$
appear in exactly one of these equations. The equation corresponding
to the diagonal $x_{1,1,1},\ldots,x_{n,n,n}$ is linearly independent
from the other $3$, because, when rewritten in terms of only variables
of the form $x_{i,j,k}$ with $1 \leq i,j,k <n$, it contains the
variable $x_{2,2,3}$, which for $n>4$ does not lie on a main diagonal
and is therefore not involved in one of the other $3$
equations. Therefore, for $n>4$ the kernel of the matrix $B$ has
dimension $(n-1)^3-4$.  This completes the proof. $\square$

\begin{figure}[hpt]
 \begin{center}
     \includegraphics[scale=0.25]{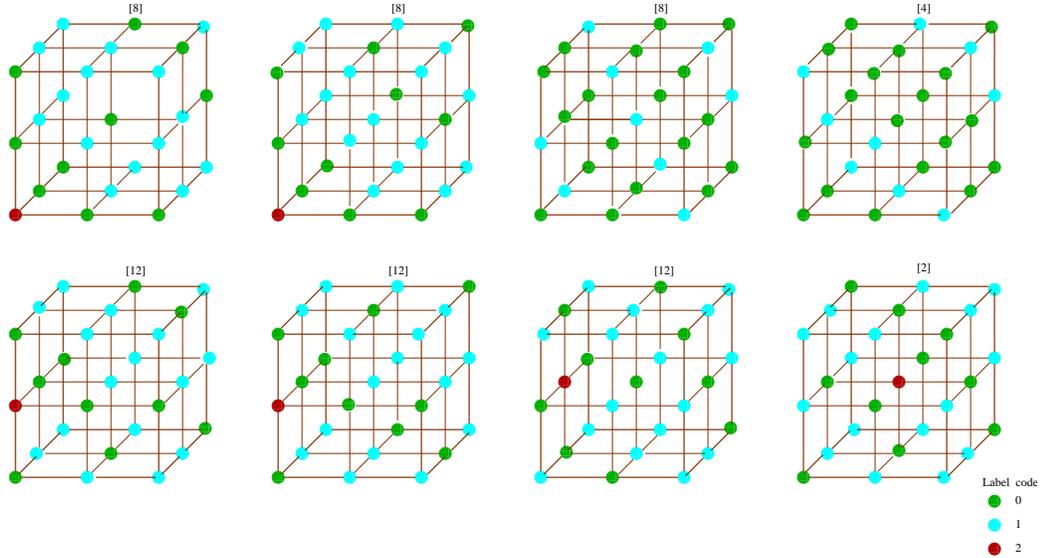}
     \caption{The generators of the minimal Hilbert basis of the $3
     \times 3 \times 3$ semi-magic cube.}
\label{semicubehb}
 \end{center}
 \end{figure}

Similarly, we can construct and enumerate semi-magic cubes. Bona
\cite{bona} had already observed that a Hilbert basis must contain
only elements of magic constant one and two. Here, we compute the 12
Hilbert basis elements of magic sum 1 and the 54 elements of magic sum
2 using 4ti2. The generators of the Hilbert basis of $3 \times 3
\times 3$ semi-magic cubes are given in Figure \ref{semicubehb} (the
definition of a generating set of the Hilbert basis of semi-magic
cubes is similar to the corresponding definition for magic cubes).

Denote by $SH^d_n(s)$ the number of semi-magic $d$-dimensional
hypercubes with $n^d$ entries. We use CoCoA to compute the generating
function $SH^d_3(s)$ :
\[
\begin{array}{l}
\sum_{s=0}^{\infty}
SH^3_3(s)t^s=\frac{t^8+5t^7+67t^6+130t^5+242t^4+130t^3+67t^2+5t+1}{(1-t)^9(1+t)^2} \\ \\
=1+12t+132t^2+847t^3+3921t^4+14286t^5+43687t^6+116757t^7+\ldots.
\end{array}
\]

In \cite{bona}, Bona presented a proof that the counting function of
 $3 \times 3 \times 3$ semi-magic cubes is a quasi-polynomial of
non-trivial period. We improve on his result by computing an explicit
formula.

\begin{thm}
The number of  $3 \times 3 \times 3$ semi-magic cubes of magic sum $s$,
\[
SH^3_3(s) = \left\{
\begin{array}{l}
{\frac {9}{2240}}\,{s}^{8}+{\frac {27}{560}}\,{s}^{7}+{\frac {87}{320}
}\,{s}^{6}+{\frac {297}{320}}\,{s}^{5}+{\frac {1341}{640}}\,{s}^{4}+{
\frac {513}{160}}\,{s}^{3}+{\frac {3653}{1120}}\,{s}^{2}+{\frac {627}{
280}}\,s+1
     \\ \hfill  \mbox {if $2 \vert s$,} \\ \\
    
{\frac {9}{2240}}\,{s}^{8}+{\frac {27}{560}}\,{s}^{7}+{\frac {87}{320}
}\,{s}^{6}+{\frac {297}{320}}\,{s}^{5}+{\frac {1341}{640}}\,{s}^{4}+{
\frac {513}{160}}\,{s}^{3}+{\frac {3653}{1120}}\,{s}^{2}+{\frac {4071}
{2240}}\,s+{\frac {47}{128}} \\ \hfill \mbox{otherwise.} 
\end{array}
\right.
\]
\end{thm}  
The convex hull of all real nonnegative semi-magic cubes (of given
size) all whose mandated sums equal 1 is called the {\em polytope of
stochastic semi-magic cubes}.  The polytope of  $3\times 3
\times 3$ stochastic semi-magic cubes is actually not equal to the convex hull of
integral semi-magic cubes. This is because the 54 elements of degree
two in the Hilbert basis, when appropriately normalized, give rational
stochastic matrices that are all vertices. In other words, the
Birkhoff-von Neumann theorem \cite[page 108]{schrijver} about
stochastic semi-magic matrices is false for $3 \times 3 \times 3$
stochastic semi-magic cubes. We prove the following result about the
number of vertices of stochastic semi-magic cubes.
 
\begin{thm} \label{thmsemicubomagic}
The number of vertices of the polytope of $n \times n \times n$
stochastic semi-magic cubes is bounded below by $(n!)^{2n}/n^{n^2}$.
\end{thm}

\noindent {\em Proof.}
We exhibit a bijection between integral stochastic semi-magic cubes
and $n\times n$ latin squares: Each $2$-dimensional layer or slice of
the integral stochastic cubes are permutation matrices (by
Birkhoff-Von Neumann theorem), the different slices or layers cannot
have overlapping entries else that would violate the fact that along a
line the sum of the entries equals one. Thus make the permutation
coming from the first slice be the first row of the latin square, the
second slice permutation gives the second row of the latin square,
etc. From well-known bounds for latin squares we obtain the lower
bound  (see Theorem 17.2 in \cite{lintwilson}). $\square$

\newpage
\pagestyle{myheadings} 
\markright{  \rm \normalsize CHAPTER 3. \hspace{0.5cm}
 Franklin Squares }
\large 
\chapter{Franklin Squares} \label{franklinchapter}
\thispagestyle{myheadings}
%CHAPTER --- Franklin Squares.
{\small
\begin{verse}
If back we look on ancient Sages Schemes, \\
They seem ridiculous as Childrens Dreams

-- Benjamin Franklin.
% in Poor Richahrd's Almanack, June. 
\end{verse}
}
\section{All about $8 \times 8$ Franklin squares.} \label{8x8Franklin}
Like in the case of magic squares, we consider the entries of an $n
\times n$ Franklin square as variables $y_{ij} \hspace{0.05in} (1 \leq
i,j \leq n)$ and set the first row sum equal to all other mandatory
sums. Thus, Franklin squares become nonnegative integral solutions to
a system of linear equations $Ay = 0$, where $A$ is an $(n^2+8n-1)
\times n^2$ matrix each of whose entries is 0, 1, or -1.

In the case of the $8 \times 8$ Franklin squares, there are seven
linear relations equating the first row sum to all other row sums and
eight more equating the first row sum to column sums. Similarly,
equating the eight half-row sums and the eight half-column sums to the
first row sum generates sixteen linear equations. Equating the four
sets of parallel bent diagonal sums to the first row sum produces
another thirty-two equations. We obtain a further sixty-four equations
by setting all the $2 \times 2$ subsquare sums equal to the first row
sum. Thus, there are a total of 127 linear equations that define the
cone of $8 \times 8$ Franklin squares.  The coefficient matrix $A$ has
rank 54 and therefore the cone $C$ of $8 \times 8$ Franklin squares
has dimension 10.

Let $(r_i,r_j)$ denote the operator that acts on the space of $n
\times n$ matrices by interchanging rows $i$ and $j$ of each matrix,
and let $(c_i,c_j)$ signify the corresponding operator on columns.

Consider the group $G$ of symmetry operations of $8 \times 8$ Franklin
squares (see Lemma \ref{G}): $G$ is generated by the set
\[\{(c_1,c_3),(c_5,c_7),(c_2,c_4),(c_6,c_8),
(r_1,r_3),(r_5,r_7), (r_2,r_4),(r_6,r_8)\}.\]

\begin{figure}[h]
 \begin{center}
     \includegraphics[scale=0.5]{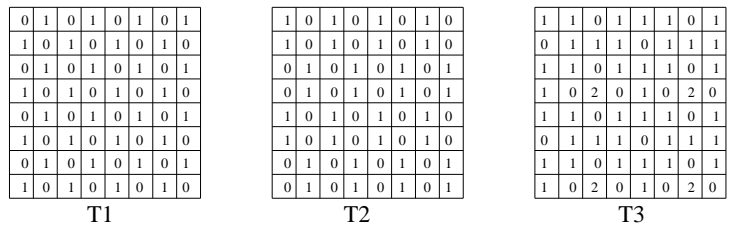}
\caption{Generators of the Hilbert basis of $8 \times 8$ Franklin squares.}
\label{hilb8x8} \end{center} \end{figure}

The Hilbert basis of the polyhedral cone of $8 \times 8$ Franklin
squares is generated by the action of the group $G$ on the three
squares T1, T2, and T3 in Figure \ref{hilb8x8} and their
counterclockwise rotations through 90 degree angles.  Not all squares
generated by these operations are distinct. Let $R$ denote the
operation of rotating a square 90 degrees in the counterclockwise
direction. Observe that $R^2 \cdot$T1 is the same as T1 and $R^3
\cdot$T1 coincides with $R \cdot$T1.  Similarly, $R^2 \cdot$T2 is just
T2, and $R^3 \cdot$T2 is the same as $R \cdot$T2. Also T1 and $R
\cdot$T1 are invariant under the action of the group $G$.  Therefore
the Hilbert basis of the polyhedral cone of $8 \times 8$ Franklin
squares consists of the ninety-eight Franklin squares: T1 and $R \cdot$T1;
the thirty-two squares generated by the action of $G$ on T2 and $R \cdot$T2;
the sixty-four squares generated by the action of $G$ on T3 and its
three rotations $R \cdot$T3, $R^2 \cdot$T3, and $R^3 \cdot$T3.

%%%%%%%  Hilbert basis  %%%%%%%%%%%%%%%%%
%%%%%%%%%%%%%%%%%%%%%%%%%%%%%%%%%%%%%%%%%%%%%%%%

\begin{figure}[h]
 \begin{center}
     \includegraphics[scale=0.4]{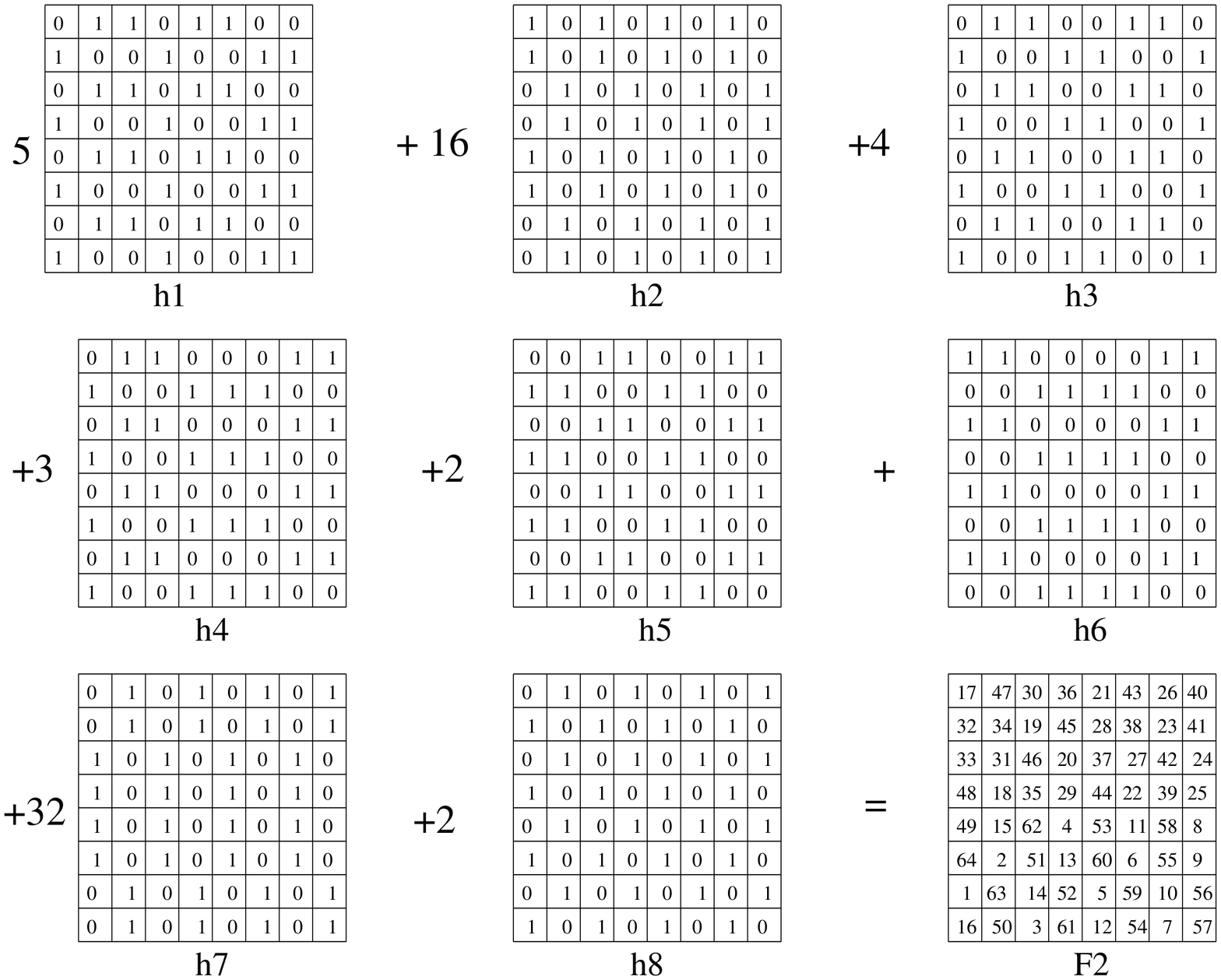}
\caption{Constructing Benjamin Franklin's $8 \times 8$ square F2.} \label{constructf2}
 \end{center}
 \end{figure}

\begin{figure}[h]
 \begin{center}
     \includegraphics[scale=0.4]{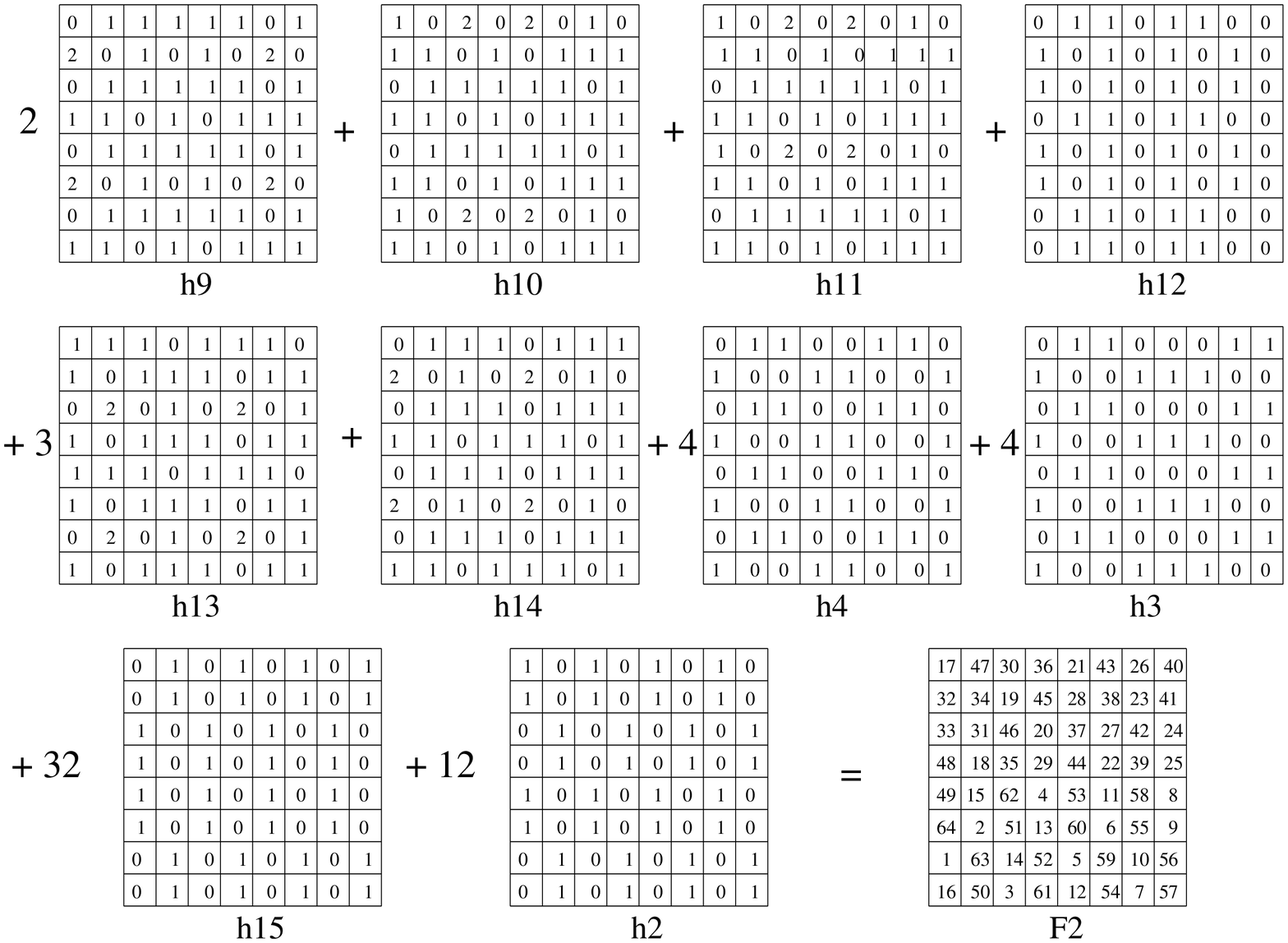}
\caption{Another construction of Benjamin Franklin's $8 \times 8$ square F2.} 
\label{anotherconstructf2}
 \end{center}
 \end{figure}

The Hilbert basis constructions of the Franklin squares F2, N1, N2,
F1, and N3 read as follows (see Figures \ref{constructf2} and
\ref{anotherconstructf2} for clarification of the notation):
\[
\begin{array}{lll}
\mbox{F2} & = &5 \cdot \mbox{h}1 + 16 \cdot \mbox{h}2 + 4\mbox{h}3+ 3 \cdot \mbox{h}4 + 2 \cdot \mbox{h}5 + \mbox{h}6 + 
+ 32 \cdot \mbox{h}7 + 2 \cdot \mbox{h}8; \\ \\

\mbox{N1} &= &5 \cdot \mbox{h}1 + 2 \cdot \mbox{h}2 + 4 \cdot \mbox{h}3+ 3 \cdot \mbox{h}4 + 2 \cdot \mbox{h}5 + \mbox{h}6  
+ 32 \cdot \mbox{h}7 + 16 \cdot \mbox{h}8; \\ \\

\mbox{N2}& = &5 \cdot \mbox{h}1 + 16 \cdot \mbox{h}2 + 4 \cdot
(c_5,c_7) \cdot (c_6,c_8) \cdot \mbox{h}3+ 3 \cdot \mbox{h}4 + 2 \cdot
\mbox{h}5 + \mbox{h}6 + 32 \cdot \mbox{h}7 \\ && + 2 \cdot \mbox{h}8 ; \\
\\

\mbox{F1} &=& 2 \cdot \mbox{h}1 + 14 \cdot (c1,c3) \cdot \mbox{h}1 +
\mbox{h}2 + (r6,r8) \cdot \mbox{h}2 + 3 \cdot (r1,r3) \cdot (r6,r8)
\cdot \mbox{h}2 \\ && + 30 \cdot \mbox{h}6 + 2 \cdot (c5,c7) \cdot
\mbox{h}6 + 6 \cdot \mbox{h}8 + 3 \cdot (c5,c7) \cdot \mbox{T3} +
(r2,r4) \cdot (c5,c7) \cdot \mbox{T3}; \\ \\

\mbox{N3} &=& 2 \cdot \mbox{h}1 + 6 \cdot (c1,c3) \cdot \mbox{h}1 +
\mbox{h}2 + (r6,r8) \cdot \mbox{h}2 + 3 \cdot (r1,r3) \cdot (r6,r8)
\cdot \mbox{h}2 \\ && + 30 \cdot \mbox{h}6 + 2 \cdot (c5,c7) \cdot
\mbox{h}6 + 14 \cdot \mbox{h}8 + 3 \cdot (c5,c7) \cdot \mbox{T3} +
(r2,r4) \cdot (c5,c7) \cdot \mbox{T3}. \\ \\
\end{array}
\]

%%% constructing 8x8 squares %%%%%%%%
%%%%%%%%%%%%%%%%%%%%%%%%%%%%%%%%%%%%%%%%%%%%%%%%%%%%%%%%%%%%%%%%

Similarly, let 
\[
\begin{array}{l}
\mbox{g}1 = (r_3,r_5)\cdot (r_4,r_6)\cdot
(r_{11},r_{13})\cdot (r_{12},r_{14}) \cdot \mbox{S2}, \hspace{0.05in}  \mbox{g}2 =
\mbox{S1}, \hspace{0.05in}  \mbox{g}3 = R \cdot \mbox{g}1, \\
\mbox{g}4 =(r_1,r_5)\cdot(r_4,r_8)
\cdot(r_{9},r_{13})\cdot (r_{12},r_{16})
\cdot \mbox{S2}, \hspace{0.05in} 
\mbox{g}5 = \mbox{ transpose of } S2,  \hspace{0.05in}  \\
\mbox{g}6 =(c_9,c_{13})\cdot (c_{10},c_{14})\cdot
(c_{11},c_{15})\cdot (c_{12},c_{16}) \cdot \mbox{g}5, \\ 
 \mbox{g}7 =(r_2,r_6)\cdot (r_{10},r_{14})\cdot
(r_{12},r_{16})\cdot \mbox{S2}, \hspace{0.05in} \\
 \mbox{g}8 =\mbox{S3}, \hspace{0.05in}
 \mbox{g}9=(r_1,r_5)\cdot (r_3,r_7)\cdot
(r_{10},r_{14})\cdot (r_{12},r_{16}) \cdot \mbox{S2}, \hspace{0.05in} \\
\mbox{g}10 =(r_2,r_6)\cdot(r_4,r_8)\cdot (r_{10},r_{14})
\cdot (r_{12},r_{16}) \cdot \mbox{S2}, \hspace{0.05in} \\
 \mbox{g}11 =(r_2,r_6)\cdot (r_3,r_7)\cdot
(r_{9},r_{13}) \cdot (r_{10},r_{14}) 
 \cdot(r_{12},r_{16}) \cdot
\mbox{S2}.
\end{array}
\]

These constructions, as we saw before in Section \ref{method} are not
unique.  A different construction of F2 is given in Figure
\ref{anotherconstructf2}.

%%% constructing 8x8 squares %%%%%%%%
%%%%%%%%%%%%%%%%%%%%%%%%%%%%%%%%%%%%%%%%%%%%%%%%%%%%%%%%%%%%%%%%

\begin{figure}[h]
 \begin{center} 
 \includegraphics[scale=0.5]{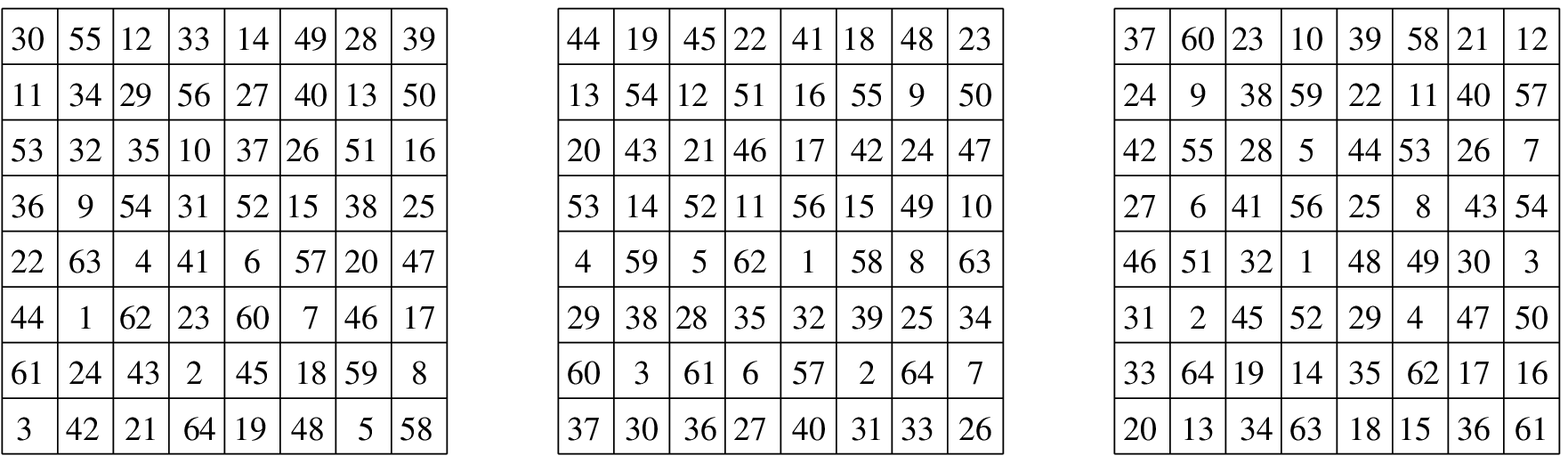}
\caption{Pandiagonal Franklin squares constructed by Ray Hagstorm \cite{hagstorm}.}\label{hagstorm}
 \end{center}
 \end{figure}

Interestingly, the Hilbert basis of $8 \times 8$ pandiagonal Franklin
squares is a subset of the Hilbert basis of $8 \times 8$ Franklin
squares. The thirty-two squares generated by the action of the group
$G$ on T2 and $R \cdot$T2 form the Hilbert basis of $8 \times 8$
pandiagonal Franklin squares. The pandiagonal Franklin squares in
Figure \ref{hagstorm} were constructed by Ray Hagstorm using the
minimal Hilbert basis of Pandiagonal Franklin squares \cite{hagstorm}.

%%%% Example %%%%%%%%%%%%%%%%%%%%%%%

Let $F_8(s)$ denote the number of $8 \times 8$ Franklin squares with
magic sum $s$.  We used the program CoCoA to compute the
Hilbert-Poincar\'e series $\sum_{s=0}^{\infty}F_8(s)t^s$ and obtained
\[
\begin{array}{l}

\sum_{s=0}^{\infty}F_8(s)t^s = \\ \\ \{
 ({t}^{36}-{t}^{34}+28\,{t}^{32}+33\,{t}^{30}+233\,{t}^{28}+390\,
 {t}^{26}+947\,{t}^{24}+1327\,{t}^{22}+1991\,{t}^{20} \\
 +1878\,{t}^{18} +1991\,{t}^{16} +1327\,{t}^{14}+947\,{t}^{12}
 +390\,{t}^{10}+233\,{t}^{8} +33\,{t}^{6}+28\,{t}^{4} \\ 
-{t}^{2}+1) \}/\{({t}^{2}-1)^{7}({t}^{6}-1)^{3}({t}^{2}+1)^{6} \} \\ \\

= 
1+34\,{t}^{4}+64\,{t}^{6}+483\,{t}^{8}+1152\,{t}^{10}+4228\,{t}^{12}+
9792\,{t}^{14}+25957\,{t}^{16}+ \cdots 
\end{array}
\]

We recover the Hilbert function $F_8(s)$ from the Hilbert-Poincar\'e
series by interpolation (see Section \ref{method}). The formulas for
the number of $8 \times 8$ pandiagonal Franklin squares in Theorem
\ref{enumeratepan8x8} are derived similarly.

Natural $8 \times 8$ Franklin squares always have magic sum 260. From
Theorem \ref{enumerate8x8} we find that $F_8(260)$ is
228,881,701,845,346. This number is an upper bound for the number of
natural $8 \times 8$ Franklin squares. The actual number of such
squares is still an open question.

\section{A few aspects of $16 \times 16$ Franklin squares.} \label{16x16section}
Finding the minimal Hilbert basis for the cone of $16 \times 16$
Franklin squares is computationally challenging and remains an
unresolved problem. However, we can provide a partial Hilbert basis
that enables us to construct Benjamin Franklin's square F3, as well as
the square N4. The following lemma proves that every $8 \times 8$
Franklin square corresponds to a $16 \times 16$ Franklin square.

\begin{lemma} \label{blockslemma}
Let M be an $8 \times 8$ Franklin square. Then the square T
constructed using M as blocks (as in  Figure \ref{8x8blocks}) is a $16
\times 16$ Franklin square.
%%% constructing 16x16 squares from 8x8 squares %%%%%%%%
\begin{figure}[h]
 \begin{center}
     \includegraphics[scale=0.6]{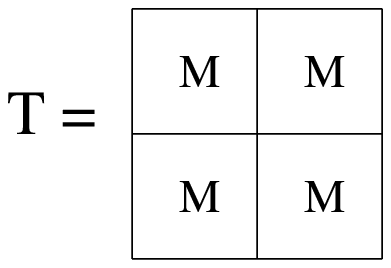}
\caption{Constructing a $16 \times 16$ Franklin square T using an $8
\times 8$ Franklin square M.} \label{8x8blocks} \end{center}
\end{figure}
%%%%%%%%%%%%%%%%%%%%%%%%%%%%%%%%%%%%%%%%%%%%%%%%%%%%%%%%%%%%%%%%
%---------------------------
\end{lemma}

\begin{figure}[h]
 \begin{center}
     \includegraphics[scale=0.4]{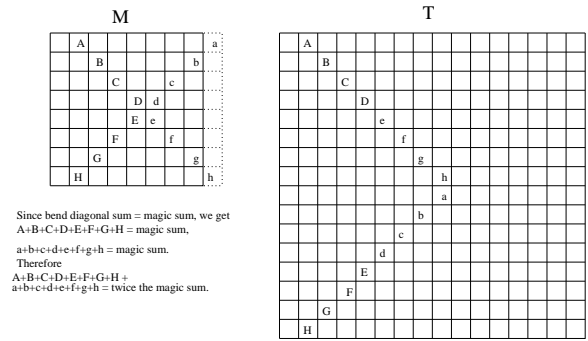}
\caption{Bent diagonals of T add to twice the magic sum of M.}
\label{explain8x8blocks} \end{center} \end{figure}

\noindent {\em Proof.} Let the magic sum of M be $s$.  The
half-columns and half-rows of T add up to $s$ since they are the
columns and rows of M, respectively. Also the columns and rows of T add
to $2s$.  The bent diagonals of T  sum to $2s$ (see Figure
\ref{explain8x8blocks} for an explanation). Since the $2 \times 2$
subsquares of M add to $s/2$, we infer that the $2 \times 2$
subsquares of T add to $s/2$. Thus T is a $16 \times 16$
Franklin square with magic sum $2s$.
$\square$

%%%%%%%%%%%%%%%%% Explanation %%%%%%%%%%%%%%%%%%%%%%%%%%%

%%%%%%%%%%%%%%%%%%%%%%%%%%%%%%%%%%%%%%%%%%%%%%%%%%%%%%%%

\begin{figure}[h]
 \begin{center}
     \includegraphics[scale=0.4]{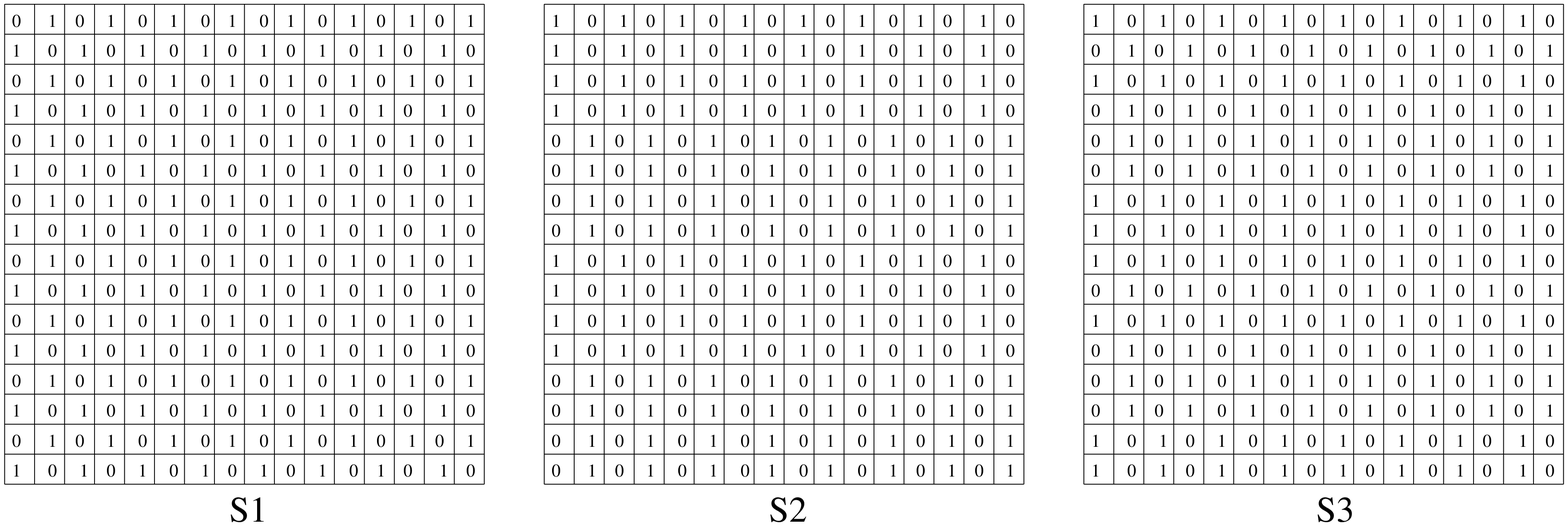}
\caption {Elements of a partial Hilbert basis of $16 \times 16$
Franklin squares.}
\label{hilbertbasis16x16} \end{center} \end{figure}

Consider the set $\cal{B}$ of $16 \times 16$ Franklin squares obtained
by applying the symmetry operations listed in Theorem \ref{symmetries}
to the squares constructed by applying Lemma \ref{blockslemma} to the
ninety-eight elements of the minimal Hilbert basis of $8 \times 8$
Franklin squares (for example, S1 in Figure \ref{hilbertbasis16x16} is
constructed from the $8 \times 8$ Franklin square T1 in Figure
\ref{hilb8x8}) and the $16 \times 16$ Franklin squares S2 and S3 in
Figure \ref{hilbertbasis16x16}. Observe that one-fourth the magic sum
of a $16 \times 16$ Franklin square is always an integer because its
$2 \times 2$ subsquares add to this number. This implies that the
squares in $\cal{B}$ are irreducible, for they have magic sums 8 or 12
(it is easy to verify that there are no $16 \times 16$ Franklin
squares of magic sum 4). Therefore, $\cal{B}$ is a subset of the
minimal Hilbert basis for the cone of $16 \times 16$ Franklin
squares. Thus, $\cal B$ forms a partial Hilbert basis.

%%% 16x16 Hilbert basis  %%%%%%%%
%%%%%%%%%%%%%%%%%%%%%%%%%%%%%%%%%%%%%%%%%%%%%%%%%%%%%%%%%%%%%%%%

We obtain F3 and N4 as follows:
\[
\begin{array}{lll}
\mbox{F3} &=& \mbox{g}1 + 17 \cdot \mbox{g}2 + 32 \cdot \mbox{g}3 + \mbox{g}4 
+ 64 \cdot \mbox{g}5 + 128 \cdot \mbox{g}6+ 2 \cdot \mbox{g}7  + 2
\cdot \mbox{g}8 + 7 \cdot \mbox{g}9 \\ && +  \mbox{g}10 + 2 \cdot \mbox{g}11; \\ \\
\mbox{N4} &=& \mbox{g}1 + 17 \cdot \mbox{g}2 + 64 \cdot \mbox{g}3 + \mbox{g}4 
+ 32 \cdot \mbox{g}5 + 128 \cdot \mbox{g}6+ 2 \cdot \mbox{g}7  + 2
\cdot \mbox{g}8 + 7 \cdot \mbox{g}9  \\ && +  \mbox{g}10 + 2 \cdot \mbox{g}11.

\end{array}
\]

 Since every $8 \times 8$ Franklin square corresponds to a $16 \times
16$ Franklin square by Lemma \ref{blockslemma}, the formulas for the
number of $8 \times 8$ Franklin squares of magic sum $s$ in Theorem
\ref{enumerate8x8} also yields a lower bound for the number of $16
\times 16$ Franklin squares of magic sum $2s$.

%%%%%%%  symmetries of square section %%%%%%%%%%%%%%%%%%%%%%%%%%%

\section{Symmetries of Franklin Squares.} \label{symmetry}

In this section we prove Theorem \ref{nonsymmetric}, which asserts
 that the new Franklin squares N1, N2, N3, and N4 are not derived from
 Benjamin Franklin's squares F1, F2, or F3 by symmetry operations. We
 first prove Theorem \ref{symmetries}.  

Rotation, reflection, and taking the transpose are plainly symmetry
 operations on Franklin squares. The proof of Theorem \ref{symmetries}
 follows from Lemmas \ref{G}, \ref{S}, \ref{fourrowslemma}, and
 \ref{switchrowslemma}.

Let $S_n$ denote the group of $n \times n$ permutation matrices acting
on $n \times n$ matrices.  As earlier, let $(r_i,r_j)$ denote the operation of
exchanging rows $i$ and $j$ of a square matrix, and let $(c_i,c_j)$
denote the analogous operation on columns.

\begin{lemma} \label{G}
Let $G$ be the subgroup of $S_8$ generated by
\[ 
\{(c_1,c_3),(c_5,c_7),(c_2,c_4),(c_6,c_8),(r_1,r_3),(r_5,r_7),
   (r_2,r_4),(r_6,r_8) \},
\] and let $H$ be the subgroup of $S_{16}$ generated by
\[ \begin{array}{l}
 \{(c_1,c_3),(c_2,c_4),(c_3,c_5),(c_4,c_6),(c_5,c_7) (c_6,c_8),
   (c_9,c_{11}),(c_{10},c_{12}), \\
   (c_{11},c_{13}),(c_{12},c_{14}),(c_{13},c_{15}) (c_{14},c_{16}),
   (r_1,r_3),(r_2,r_4),(r_3,r_5),(r_4,r_6), \\(r_5,r_7) (r_6,r_8),
   (r_9,r_{11}),(r_{10},r_{12}),(r_{11},r_{13}),(r_{12},r_{14}),
   (r_{13},r_{15}),(r_{14},r_{16}) \}.
\end{array}\]

\noindent The row and column permutations from the group $G$ map $8
\times 8$ Franklin squares to $8 \times 8$ Franklin squares, while the
row and column permutations from the group $H$ map $16 \times 16$
Franklin squares to $16 \times 16$ Franklin squares.
\end{lemma}

\begin{figure}[h]
 \begin{center}
     \includegraphics[scale=0.5]{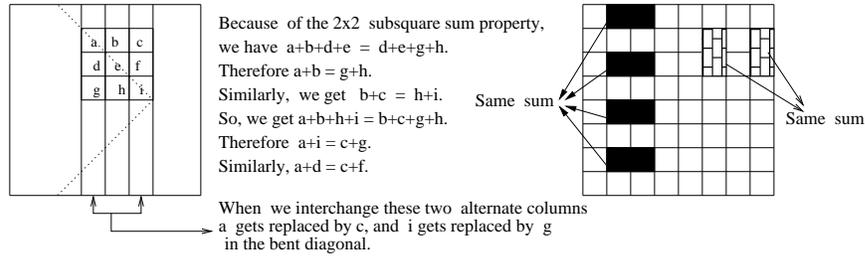}
\caption{Properties of Franklin squares. }
\label{explain} \end{center} \end{figure}

\noindent {\em Proof.} Clearly exchanging rows or columns of a
Franklin square preserves row and column sums.  Half-row and
half-column sums are preserved because the permutations of rows and
columns included here operate in some half of a Franklin square. That
$2 \times 2$ subsquare sums are preserved follows from the fact that
every alternate pair of entries in a pair of columns or rows add to
the same sum (see Figure \ref{explain} for an explanation).  For any
$3 \times 3$ subsquare of a Franklin square, the two sums of
diagonally opposite elements are equal (see Figure \ref{explain} for
details). This implies that, if we permute alternate rows or alternate
columns then the new entries preserve bent diagonal sums. Observe that for
the preservation of  bent diagonal sums, it is critical that the alternate row and
column permutations be restricted to act in one half of a Franklin
square (see Figure \ref{newexample} for examples).              
$\square$

%%%%%%%  explain  %%%%%%%%%%%%%%%%%
%%%%%%%%%%%%%%%%%%%%%%%%%%%%%%%%%%%%%%%%%%%%%%%%

\begin{lemma} \label{S}
Let $S$ be the subgroup of $S_{16}$ generated by the set 
\[
\begin{array}{l}
 \{(r_1,r_5), (r_2,r_6), (r_3,r_7), (r_4,r_8), (r_9,r_{13}), 
(r_{10},r_{14}),(r_{11},r_{15}),(r_{12},r_{16}), \\
(c_1,c_5),(c_2,c_6),(c_3,c_7),(c_4,c_8),(c_9,c_{13}),(c_{10},c_{14}) 
,(c_{11},c_{15}),(c_{12},c_{16})\}.
\end{array}
\]
The row and column permutations from the subgroup $S$ map $16 \times 16$
Franklin squares to $16 \times 16$ Franklin squares.

\end{lemma}

\begin{figure}[h]
 \begin{center}
     \includegraphics[scale=0.5]{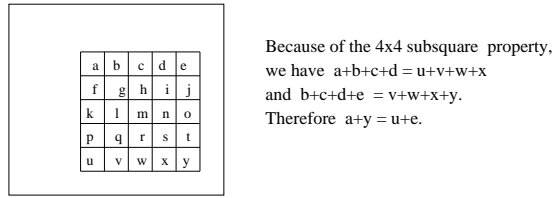}
\caption{Properties of $16 \times 16$ Franklin squares. }
\label{explain16x16} \end{center} \end{figure}

\noindent {\em Proof.} Half-column and half-row sums are preserved because the
specified row or column exchanges only affect some half of the $16
\times 16$ Franklin square.  Because of the $2 \times 2$ subsquare sum
property, we see that every $4 \times 4$ subsquare adds to the common
magic sum. Hence the two sums of diagonally opposite elements in a $5
\times 5$ subsquare are equal (see Figure \ref{explain16x16} for an
explanation). This implies that the bent diagonal sums are preserved,
again because they operate in only one-half of a Franklin square. It
is easy to verify that all other sums are preserved under the action
of elements of $S$. $\square$

%%%%%%%  explain  %%%%%%%%%%%%%%%%%
%%%%%%%%%%%%%%%%%%%%%%%%%%%%%%%%%%%%%%%%%%%%%%%%

Observe that the groups $G$, $H$, and $S$ are commutative, for each
nonidentity element in each of these groups has order 2. Therefore the
order of $G$ is $2^8$, the order of $H$ is $2^{24}$, and the order of
$S$ is $2^{16}$.

\begin{lemma} \label{fourrowslemma}
The operation of interchanging the first $n/2$ columns
(respectively, rows) and the last $n/2$ columns (respectively, rows) of an
 $n \times n$ Franklin square is a symmetry operation.
\end{lemma}

\begin{figure}[h]
 \begin{center}
     \includegraphics[scale=0.5]{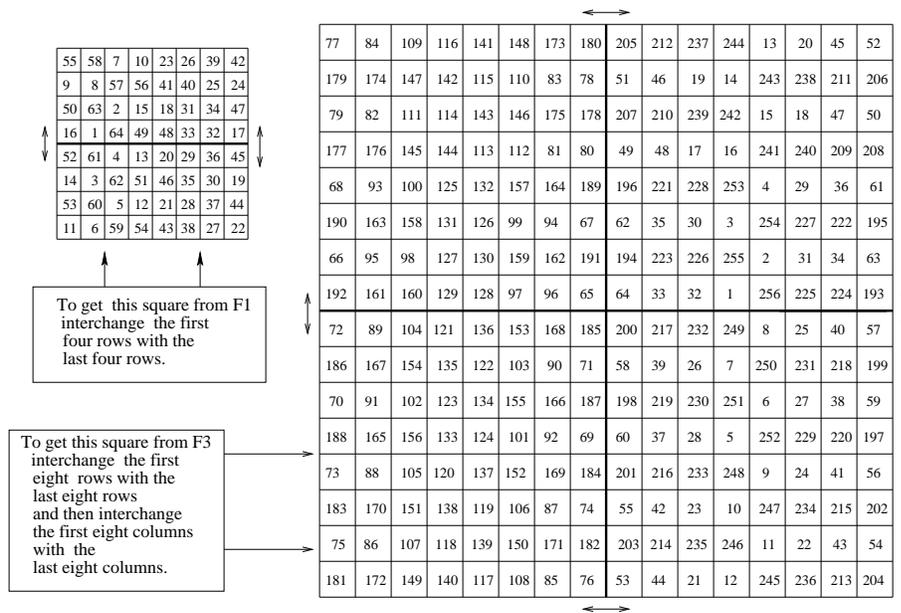}
\caption{Constructing Franklin squares by simultaneous row (or
column) exchanges of Franklin squares described in Lemma
\ref{fourrowslemma}.}
\label{fourrows} \end{center}
\end{figure}

\noindent{\em Proof.} These operations preserve half-row and half-column sums. The
row sums and column sums do not change. Bent diagonal and $2 \times 2$
subsquare sums are preserved because of continuity (see Figure
\ref{fourrows} for examples). $\square$

%%% constructing 8x8 squares by exchanging four rows%%%%%%%%
%%%%%%%%%%%%%%%%%%%%%%%%%%%%%%%%%%%%%%%%%%%%%%%%%%%%%%%%%%%%%%%%

\begin{lemma} \label{switchrowslemma}
Simultaneously interchanging all the adjacent columns (respectively,
rows) $i$ and $i+1$ ($i= 1,3,5, \dots, n-1$) of an $n \times n$ Franklin
square is a symmetry operation.  
\end{lemma}

\begin{figure}[h]
 \begin{center}
     \includegraphics[scale=0.5]{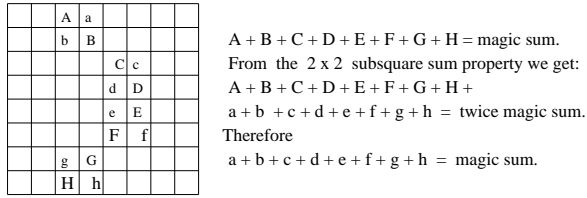}
\caption{Properties of $8 \times 8$ Franklin squares.} \label{explain2}
\end{center}
\end{figure}

\noindent {\em Proof.} It is clear that row, column, half-row, and
half-column sums are preserved by these operations. Moreover $2 \times
2$ subsquare sums are preserved because every alternate pair of
entries in a pair of columns or rows add to the same sum (see Figure
\ref{explain}). Bent diagonal sums are preserved because
of the $2 \times 2$ subsquare sum property. See Figure \ref{explain2}
for an explanation (the explanation for bent diagonal sums of $16
\times 16$ Franklin squares is similar).
$\square$

%%% Theorem explain2%%%%%%%%

Are the symmetry operations given in Theorem \ref{symmetries} all the
symmetry operations of a Franklin square? We do not know the answer to
this question, but the symmetries described in this section enable us
to interchange any two rows or any two columns within any half of a
Franklin square and get a Franklin square. Note that certain row or
column exchanges are not symmetry operations unless they are
accompanied by other simultaneous row or column exchanges. We now show
that N1, N2, N3, and N4 are not symmetric transformations of F1, F2, or
F3.

\begin{lemma} \label{F1F2case}
The squares F1 and F2 can be transformed by means of symmetry
operations neither to each other nor to any of the nonisomorphic
squares N1, N2, or N3.

\end{lemma}

\begin{figure}[h]
 \begin{center}
     \includegraphics[scale=0.5]{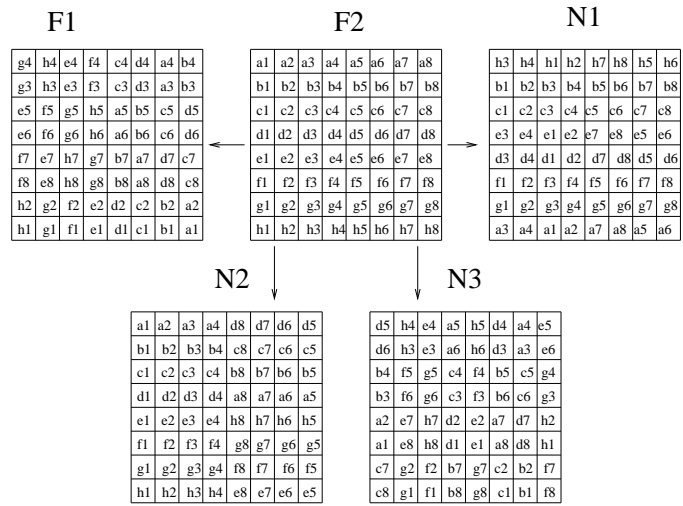}
\caption{The abstract permutation of F2 that gives F1, N1, N2, and N3.}
\label{f2tonewsquares} \end{center} \end{figure}

\noindent {\em Proof.} By definition, symmetry operations map a
Franklin square to another Franklin square. We can permute the entries
of the Franklin square F2 to get F1, N1, N2, and N3 (see Figure
\ref{f2tonewsquares}).  The permutation that maps F2 to N2 is not,
however, a symmetry operation: F1 does not transform to a Franklin
square under this permutation since bent diagonal sums are not
preserved. The other permutations of F2 in Figure \ref{f2tonewsquares}
likewise fail to be symmetries. Again, F1 does not map to a Franklin
square under these permutations for half-column sums are not
preserved.  Thus F2 cannot be transformed to F1, N1, N2, or N3 using
symmetry operations. Similarly, the permutations of the entries of F1
that map it to F2, N1, N2, and N3, respectively, are not symmetry
operations because F2 is not mapped to a Franklin square under any of
these permutations (in these instances half-row sums are not
preserved).  The permutations that map the square N1 to N2 and N3, and
the permutations that map N2 to N1 and N3 are not symmetry operations
because F1 is not mapped to a Franklin square under these permutations
(half-column sums of F1 are not preserved for all these
permutations). Similarly, the permutations that map N3 to N1 and N2
are not symmetry operations because F2 is not mapped to a Franklin
square under these operations. Therefore, the squares N1, N2, and N3
are not isomorphic to each other. $\square$

%%%%%%% F2 to F1 and new squares N1 and N2 %%%%%%%%%%%%%%%%%
%%%%%%%%%%%%%%%%%%%%%%%%%%%%%%%%%%%%%%%%%%%%%%%%

\begin{lemma} \label{F3case}
The square F3 cannot be transformed to  N4  using symmetry operations.
\end{lemma}
\noindent {\em Proof.} Permuting the entries of F3 to get N4 is
achieved by simultaneously interchanging columns 1 and 15, columns 2
and 16, columns 7 and 9, and columns 8 and 10 of F3. This permutation is
not a symmetry operation. To see this, note that the square A obtained
by transposing F3 is a Franklin square (Theorem \ref{symmetries}).
But $(c_1,c_{15})(c_{2},c_{16})(c_{7},c_9)(c_8,c_{10})\cdot$A is not a
Franklin square, since bent diagonal sums are not preserved.
$\square$

Lemmas \ref{F1F2case} and \ref{F3case}, in tandem, establish Theorem
\ref{nonsymmetric}.

\newpage
\pagestyle{myheadings} 
\markright{  \rm \normalsize CHAPTER 4. \hspace{0.5cm}
The  Magic Squares and Magic Graphs Connection}
\large 
\chapter{Symmetric Magic Squares and the Magic Graphs Connection} \label{magicgraphchapter}
\thispagestyle{myheadings}
{ \small
\begin{verse}
What immortal hand or eye \\
Dare frame thy fearful symmetry?

-- William Blake.
% in Tiger Tiger burning bright
\end{verse}
}
\section{Hilbert bases of polyhedral cones of magic labelings.}
In this section we derive some results about Hilbert bases of cones of
magic labelings of graphs and also prove Proposition
\ref{digraphperfectmatchings}.

\begin{lemma}
Let $G$ be a graph with $n$ vertices. A labeling $L$ of $G$ with magic
sum $s$ can be lifted to a magic labeling $L^{\prime}$ of $\Gamma_n$ with
magic sum $s$.
\end{lemma}

\noindent {\em Proof.}
Since $G$ is a subgraph of $\Gamma_n$, every
labeling $L$ of $G$ can be lifted to a labeling $L^{\prime}$ of
$\Gamma_n$, where
\[
L^{\prime} (e_{ij}) = \left \{ \begin{array}{ll} L(e_{ij}) & \mbox{ if
$e_{ij}$ is also an edge of $G$,} \\ 0 & \mbox{otherwise.}
\end{array}
\right .
\]

Since  the edges with nonzero labels are the same for both $L$  and 
$L^{\prime}$, it follows that the magic sums are also the same. $\square$ 

\begin{lemma} \label{hilbGprop}
Let $G$ be a graph with $n$ vertices. The minimal Hilbert basis of
$C_G$ can be lifted to a subset of the minimal Hilbert basis of
$C_{\Gamma_n}$.
\end{lemma}
\noindent {\em Proof.}  If $L$ is an irreducible magic labeling of
$G$, then clearly it lifts to an irreducible magic labeling $L^{\prime}$ of
$\Gamma_n$. Since the minimal Hilbert basis is the set of all
irreducible magic labelings, we get that the minimal Hilbert basis of
$C_G$ corresponds to a subset of the minimal Hilbert basis of
$C_{\Gamma_n}$.  $\square$

For example, the magic labelings O1 and O2 of the octahedral graph
in Figure \ref{octahilb} correspond to the magic labelings a1 and
a2, respectively, of $\Gamma_6$ (see Figure \ref{symm6hilb}). 

Similarly, we can prove:
\begin{lemma} \label{digraphhilb} 
For a digraph $D$ with $n$ vertices, a magic labeling $L$ with magic
sum $s$ can be lifted to a magic labeling $L^{\prime}$ of $\Pi_n$ with
the same magic sum $s$. The minimal Hilbert basis of $C_D$ can be
lifted to a subset of the minimal Hilbert basis of $C_{\Pi_n}$.
\end{lemma}

\begin{lemma} \label{magicsum1}
Let $D$ be a digraph with $n$ vertices. All the elements of the
minimal Hilbert basis of $C_D$ have magic sum 1.
\end{lemma}
\noindent {\em Proof.}  It is well-known that the minimal Hilbert
basis of semi-magic squares are the permutation matrices (see
\cite{schrijver}) and therefore have magic sum 1. The one-to-one
correspondence between magic labelings of $\Pi_n$ and semi-magic
squares implies that the minimal Hilbert basis elements of $C_{\Pi_n}$
have magic sum 1. It follows by Lemma \ref{digraphhilb} that all the
elements of the minimal Hilbert basis of $C_{D}$ have magic sum
1. $\square$

We now present the proof of the fact that perfect matchings of
bipartite graphs correspond to the elements of the minimal Hilbert
basis of its corresponding digraph.

\noindent {\em Proof of Proposition \ref{digraphperfectmatchings}.}  

Every element of the minimal Hilbert basis of $C_{B_D}$ corresponds to
a perfect matching of $B$ by Lemma \ref{magicsum1}. Moreover, all the
magic labelings of $B_D$ of magic sum 1 belong to the minimal Hilbert
basis of $C_{B_D}$ because they are irreducible. Since perfect
matchings of $B$ are in one-to-one correspondence with magic labelings
of $B_D$ of magic sum 1, we derive that perfect matchings of $B$ are in
one-to-one correspondence with the elements of the minimal Hilbert
basis of $C_{B_D}$. It follows that $H_{B_D}(1)$ is the number of
perfect matchings of $B$.  $\square$

\begin{figure}[h]
 \begin{center} 
 \includegraphics[scale=0.5]{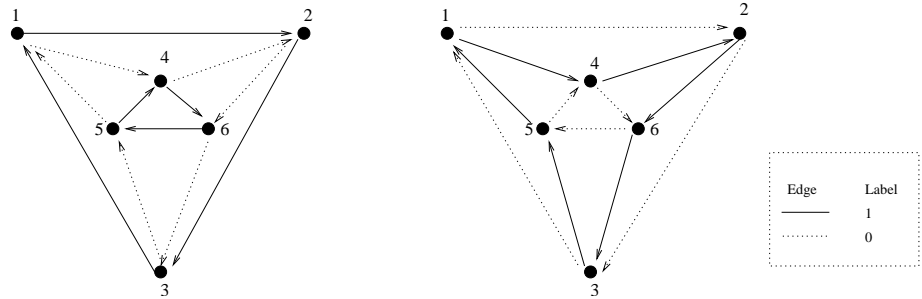}
\caption{The minimal Hilbert basis of the cone of magic labelings of
$D_O$.} \label{diroctexamp}
\end{center} \end{figure}

\begin{figure}[h]
 \begin{center} 
 \includegraphics[scale=0.5]{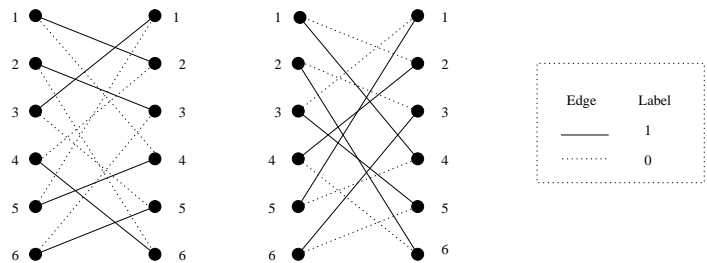}
\caption{Perfect matchings of $G_{D_O}$ corresponding to the minimal
Hilbert basis elements of $D_O$ (see Figure \ref{diroctexamp}).}
\label{bipperfect}
\end{center} \end{figure}

For example, consider the Octahedral graph with the given orientation
$D_O$ in Figure \ref{octexample}. The minimal Hilbert basis of $D_O$
is given in Figure \ref{diroctexamp}.  The perfect matchings of the
bipartite graph $G_{D_O}$ corresponding to the minimal Hilbert basis
elements of $C_{D_O}$ is given in Figure \ref{bipperfect}. We derive
$H_{D_O}(r) = r+1$ and thereby verify that the number of perfect
matchings of $G_{D_O}$ is indeed 2.

%%%% Example %%%%%%%%%%%%%%%%%%%%%%%

%%%% Example %%%%%%%%%%%%%%%%%%%%%%%

%%%%%%%%%%%%%%%%%%%%%%%%%%%%%%%%%%%%%%%%%%%%%%%%%%
%%%%%%%%%%%%%%%%%%%%%%%%%%%%%%%%%%%%%%%%%%%%%%%%%%
%%%%%%%%%%%%%%%%%%%%%%%%%%%%%%%%%%%%%%%%%%%%%%%%%

\section{Counting isomorphic simple labelings and Invariant rings.}
Let $S_n$ denote the group of permutations that acts on the vertex
set $\{v_1,v_2, \dots, v_n \}$ of $G$. Let $e_{ij}$ denote an edge
between the vertices $v_i$ and $v_j$. The action of $S_n$ on the
vertices of $G$ translates to an action on the labels of the edges of
$G$ by
\[
\sigma (L(e_{ij})) = L(e_{\sigma(i) \sigma(j)}) \mbox{ where } \sigma \in S_n.
\] 

%%%% Example %%%%%%%%%%%%%%%%%%%%%%%
\begin{figure}[h]
 \begin{center} 
 \includegraphics[scale=0.5]{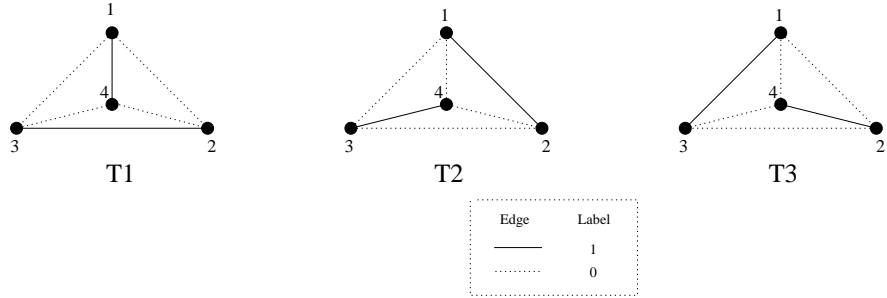}
\caption{The minimal Hilbert basis of magic labelings of the tetrahedral Graph.}
 \label{tetrahilb}
 \end{center} \end{figure}

Two labelings $L$ and $L^{\prime}$ of $G$ are {\em isomorphic} if
there exists a permutation $\sigma$ in $S_n$ such that
$L^{\prime}(e_{ij}) = L(e_{\sigma(i) \sigma(j)})$, i.e. $\sigma(L) =
L^{\prime}$. A set $S = \{g_1, g_2, \dots, g_r \}$ is said to {\em
generate the minimal Hilbert basis} of the cone of magic labelings of
$G$ if every element of the minimal Hilbert basis is isomorphic to
some $g_i$ in $S$.  For example, T1 generates the minimal Hilbert
basis of the cone of magic labelings of the tetrahedral graph (see
Figure
\ref{tetrahilb}). Observe that we get T2 by permuting the vertices
$v_2$ and $v_4$ of T1, and T3 by permuting the vertices $v_3$ and
$v_4$ of T1.

\begin{prop}
Let $L$ be a magic labeling in the minimal Hilbert basis of the cone
 of magic labelings of a graph (or a digraph), then all the labelings
 isomorphic to $L$ also belong to the minimal Hilbert basis.
\end{prop}
\noindent {\em Proof.} Let $L^{\prime}$ be a labeling isomorphic to
$L$ and $\sigma$ in $S_n$ be such that $\sigma(L) =
L^{\prime}$. Suppose $L^{\prime}$ does not belong to the Hilbert
basis. Then $L^{\prime}$ is reducible and can be written as sum of two
labelings: $L^{\prime} = L_1 + L_2$.  But, ${\sigma}^{-1}(L^{\prime})
= L$. Therefore, ${\sigma}^{-1}(L_1) +{\sigma}^{-1}(L_2) = L$. This
is not possible because $L$ is irreducible, since it belongs to the
minimal Hilbert basis.  Therefore, we conclude that $L^{\prime}$ must
also belong to the Hilbert basis.  $\square$

A labeling of $G$ is called a {\em simple labeling} if the labels are
0 or 1. Invariant theory \cite{sturmfels} provides an efficient
algebraic method of counting isomorphic simple labelings of a graph
$G$. Let $L$ be a simple labeling of $G$. Let $X^L$ denote the
monomial
\[
X^L= \prod_{i,j=1,\dots,n} x_{ij}^{L({e_{ij}})}.
\]

Consider the polynomial
\[
{\left ( X^L \right )}^{\circledast} = \sum_{\sigma \in S_n} X^{\sigma(L)}, \hspace{0.1in} 
\mbox { where }   X^{\sigma(L)} 
=  \prod_{i,j=1,\dots,n} x_{ij}^{L(e_{\sigma(i) \sigma(j)})}. 
\]

Observe that ${\left ( X^L \right )}^{\circledast}$ is an invariant
polynomial under the action of $S_n$ on the indices of the variables
$x_{ij}$. Let $k$ be any field. The set of polynomials invariant in
the polynomial ring $k[x_{ij}]$ under the action of the group $S_n$ is
called the {\em invariant ring} of $S_n$ and is denoted by
$k[x_{ij}]^{S_n}$. See \cite{sturmfels} for an introduction to
invariant rings.

Consider the simple labeling $L_G$ of $\Gamma_n$ associated
to $G$:
\[
L_G(e_{ij}) = 
\left \{
\begin{array}{ll}
1 & \mbox{if $e_{ij}$ is an edge of $G$,} \\
0  & \mbox {otherwise.}
\end{array}
\right .
\]
Then the polynomial ${\left ( X^L \right )}^{\circledast}$ evaluated
at $L_G$ counts the number of labelings of $G$ that are isomorphic to
$L$.

%%%% Example %%%%%%%%%%%%%%%%%%%%%%%
\begin{figure}[h]
 \begin{center} 
\includegraphics[scale=0.5]{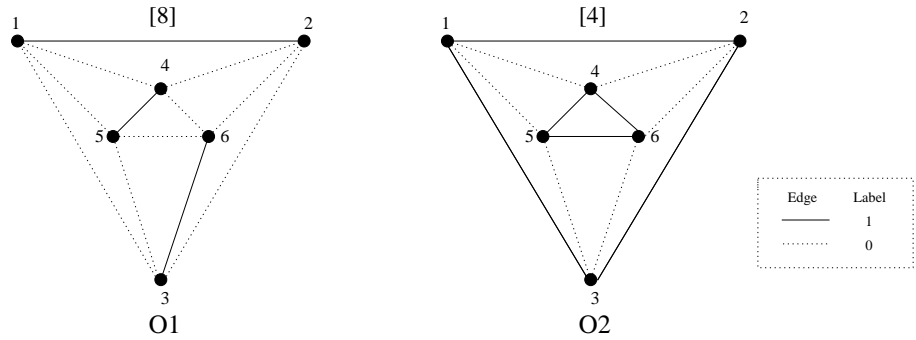}
\caption{Generators of the minimal Hilbert basis of magic labelings of the
Octahedral Graph.} \label{octahilb}  \end{center} \end{figure}

For example, consider the Octahedral graph ${\cal O}$ with the
labeling O1 in Figure \ref{octahilb}. Then,
\[\begin{array}{ll}
{ \left ( X^{O1} \right )}^{\circledast} = &
x_{14}x_{23}x_{56}+x_{16}x_{23}x_{45}+x_{13}x_{26}x_{45}+x_{12}x_{36}x_{45}+x_{15}x_{23}x_{46}+x_{13}x_{25}x_{46} +
\\  & x_{16}x_{25}x_{34}+x_{15}x_{26}x_{34}+x_{12}x_{34}x_{56}+x_{12}x_{35}x_{46}+x_{13}x_{24}x_{56}+x_{16}x_{24}x_{35}+ \\ & 
x_{14}x_{26}x_{35}+x_{15}x_{24}x_{36}+x_{14}x_{25}x_{36}.
\end{array}
\]

Substituting $ x_{12}= x_{13}= x_{14} = x_{15} = x_{23} = x_{24} =
x_{26} = x_{35} = x_{36} = x_{x45} = x_{46} = x_{56} = 1 $ and $x_{11}
= x_{22} = x_{33} = x_{44} = x_{55} = x_{66}= x_{16} = x_{25} = x_{34}
= 0$ in $ {\left(X^{O1} \right )}^{\circledast}$, we get that
\[ 
{ \left ( X^{O1} \right )}^{\circledast}(L_{\cal O} ) = 8.
\]

Therefore, there are 8 magic labelings in the $S_6$ orbit of the magic
labeling O1 of the octahedral graph. Similarly, there are four magic
labelings in the orbit of O2. The generators of the Hilbert basis of
the Octahedral graph are given in Figure \ref{octahilb}. The numbers
in square brackets in the figures indicate the number of elements in
the orbit class of each generator throughout the article.

For digraphs, we assign a variable $x_{ij}$ to every directed edge
$e_{ij}$, and use the corresponding invariant ring to count isomorphic
simple labelings. See \cite{thiery} for more aspects of labeled graph
isomorphisms and invariant rings.

\begin{figure}[h]
 \begin{center} 
 \includegraphics[scale=0.4]{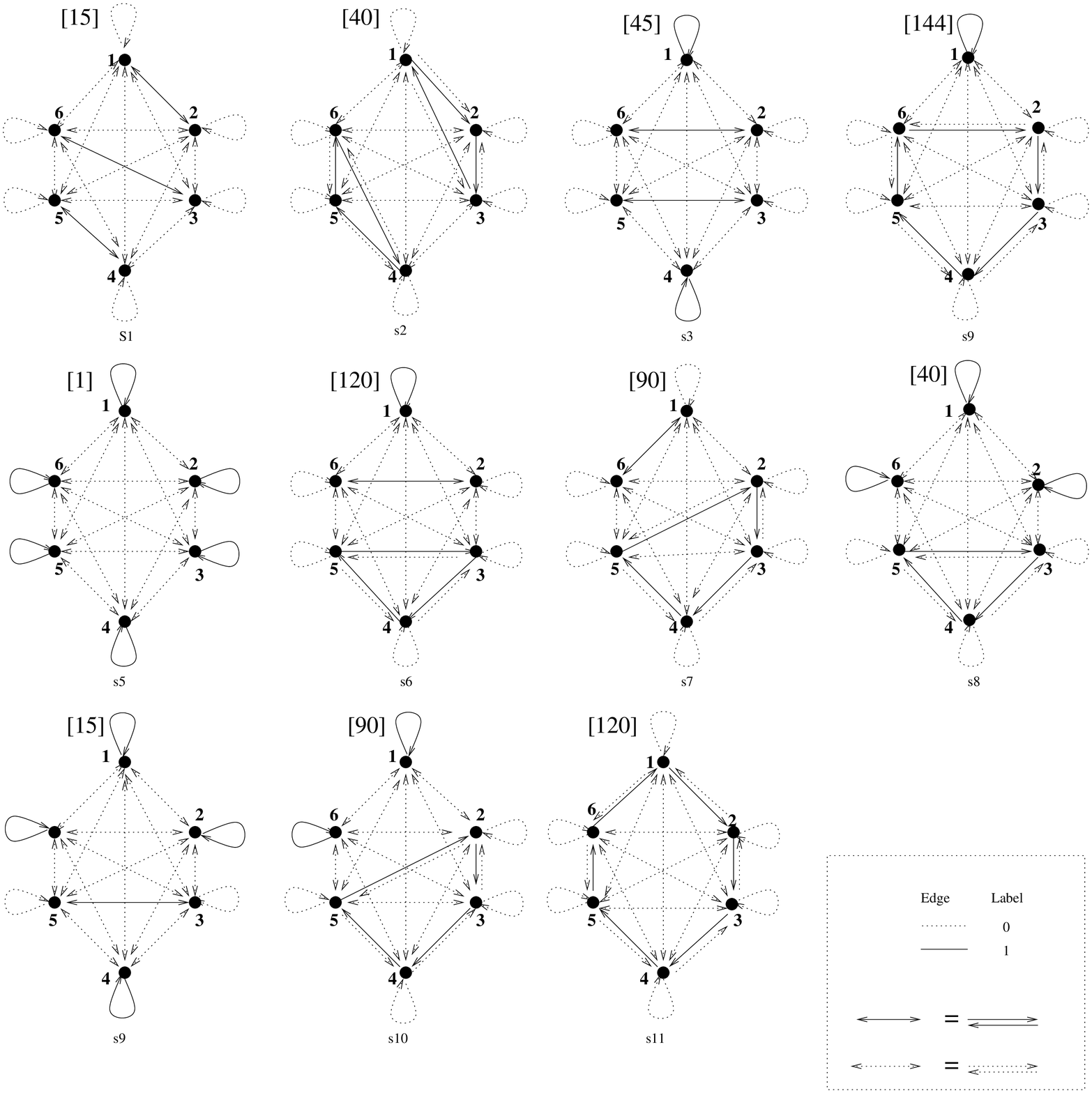}
\caption{Generators of the minimal Hilbert basis of $6 \times 6$
semi-magic squares.} \label{dir6hilb} \end{center}
\end{figure}

Since all the elements of the minimal Hilbert basis of $C_D$, where
$D$ is a digraph, are simple labelings, we can use invariant theory
effectively.  Since the number of elements in the minimal Hilbert
basis of semi-magic squares (and hence $C_{\Pi_n}$) is $n!$ (see
\cite{schrijver}), we can list all the generators of the minimal
Hilbert basis of $C_{\Pi_n}$. The generators of the Hilbert basis of
$C_{\Pi_6}$ are given in Figure \ref{dir6hilb}. 

\begin{algo} 
Computing minimal Hilbert basis of a finite digraph $D$ with $n$ vertices.
\end{algo}

Input: A digraph $D$ with $n$ vertices and the set of $n \times n$
permutation matrices.  

Output: The minimal Hilbert basis of the finite digraph $D$.

Step 0. List a set of generators of the minimal Hilbert basis of
$C_{\Pi_n}$.

Step 1. Choose all the elements $h_i$ among the generators of
the minimal Hilbert basis of $C_{\Pi_n}$ which have the edges not in
$D$ labeled 0. Delete the edges in $h_i$ that are not in $D$ to get a
magic labeling $g_i$ of $D$.

Step 2. $g_i$ and the magic labelings isomorphic to $g_i$ form the
minimal Hilbert basis of the cone of magic labelings of $D$.

For example, consider the digraph $D$ that have all the edges of
$\Pi_6$ except the loops.  Then the minimal Hilbert basis of $D$ are the 265
labelings corresponding to the labelings s1, s2, s7, and s11 in Figure
\ref{dir6hilb} and their isomorphic magic labelings.
 
%%%%%%%%%%%%%%%%%%%%%%%%%%%%%%%%%%%%%%%%%%%%%%%%%%%%%%%%%%%%%%
%%%%%%%%%%%%%%%%%%%%%%%%%%%%%%%%%%%%%%%%%%%%%%%%%%%%%%%%%%%%%
%%%%%%%%%%%%%%%%%%%%%%%%%%%%%%%%%%%%%%%%%%%%%%%%%%%%%%%%%%%%
\section{Polytopes of magic labelings.}

The proofs of Theorems \ref{graphpolytopethm},
\ref{digraphpolytopethm}, and \ref{facesbirk}, and Corollary \ref{mayathm}
 are presented in this section.
 
Let $G$ be a positive graph.  An element $\beta$ in the semigroup
$S_{C_G}$ is said to be {\em completely fundamental}, if for any
positive integer $n$ and $\alpha, {\alpha}^{\prime} \in S_{C_G}$, $n
\beta = \alpha + {\alpha}^{\prime}$ implies $\alpha = i \beta$ and
${\alpha}^{\prime} = (n-i) \beta$, for some positive integer $i$, such
that $0 \leq i \leq n$ (see \cite{stanley2}).

\begin{lemma} \label{polytopeverticeslemma} ${\cal P}_G$ is a rational 
polytope.
\end{lemma}
\noindent {\em Proof.} Proposition $4.6.10$ of Chapter 4 in
\cite{stanley2} states that the set of extreme rays of a cone and the
set of completely fundamental solutions are identical. Proposition 2.7
in \cite{stanley3} states that every completely fundamental magic
labeling of a graph $G$ has magic sum 1 or 2. Thus, the extreme rays
of the cone of magic labelings of a graph $G$ are irreducible
2-matchings of $G$. We get a vertex of ${\cal P}_{G}$ by dividing the
entries of a extreme ray by its magic sum. Thus, ${\cal P}_G$ is a
rational polytope.  $\square$

\begin{lemma} \label{polytopedimensionlemma}
The dimension of ${\cal P}_G$ is $q-n+b$, where $q$ is the number of
edges of $G$, $n$ is the number of vertices, and $b$ is the number of
connected components that are bipartite.
\end{lemma}

\noindent {\em Proof.}  Ehrhart's theorem states that the degree of
$H_G(r)$ is the dimension of ${\cal P}_G$ \cite{becketal}. The degree
of $H_G(r)$ is $q-n+b$ by Theorem \ref{degreeofG}. Therefore, the
dimension of ${\cal P}_G$ is $q-n+b$.  $\square$

\begin{lemma} \label{polytopefaceslemma}
The $d$-dimensional faces of ${\cal P}_G$ are the $d$-dimensional
polytopes of magic labelings of positive subgraphs of $G$ with $n$
vertices and at most $n-b+d$ edges.
\end{lemma}
\noindent {\em Proof.} An edge $e$ labeled with a zero in a magic
labeling $L$ of $G$ does not contribute to the magic sum, therefore,
we can consider $L$ as a magic labeling of a subgraph of $G$ with the
edge $e $ deleted. Since a face of ${\cal P}_{G}$ is the set of magic
labelings of $G$ where some edges are always labeled zero, it follows
that the face is also the set of all the magic labelings of a subgraph
of $G$ with these edges deleted.  Similarly, every magic labeling of a
subgraph $H$ with $n$ vertices corresponds to a magic labeling of $G$,
where the missing edges of $G$ in $H$ are labeled with 0. Now, let $H$
be a subgraph such that the edges $e_{r1}, \dots, e_{rm}$ are labeled
zero for every magic labeling of $H$. Then the face defined by $H$ is
same as the face defined by the positive graph we get from $H$ after
deleting the edges $e_{r1}, \dots, e_{rm}$. Therefore, the faces of
${\cal P}_G$ are polytopes of magic labelings of positive subgraphs.

By Lemma \ref{polytopedimensionlemma}, the dimension of ${\cal P}_G$
is $q-n+b$. Therefore, to get a $d$-dimensional polytope, we need to
label at least $q-n+b-d$ of $G$ edges always 0. This implies that the
$d$-dimensional face is the set of magic labelings of a positive
subgraph of $G$ with $n$ vertices and at most $n-b+d$ edges.
$\square$

The proof of Theorem \ref{graphpolytopethm} follows from Lemmas
\ref{polytopeverticeslemma}, \ref{polytopedimensionlemma}, and
\ref{polytopefaceslemma}. We can now prove Corollary \ref{mayathm}.

\noindent {\em Proof of Corollary \ref{mayathm}.}  

It is clear from the one-to-one correspondence between magic labelings
of $\Gamma_n$ and symmetric magic squares that ${\cal P}_{\Gamma_n}$
has the following description:

\[\begin{array}{lll}
{\cal P}_{\Gamma_n} &=& \{ L = (L(e_{ij}) \in {\reals}^{\frac{n(n+1)}{2}}; 
 L(e_{ij}) \geq 0; 1 \leq i,j \leq n, i \leq j,   \\  
&& \sum_{j=1}^i L(e_{ji}) + \sum_{j=i+1}^n L(e_{ij}) = 1 
\mbox{ for } i=1, \dots, n 
\}.
\end{array}
\]

Since the graph ${\Gamma_n}$ has $\frac{n(n+1)}{2}$ edges and $n$
vertices, and every graph is a subgraph of $\Gamma_n$, it follows from
Theorem \ref{graphpolytopethm} that the dimension of ${\cal
P}_{\Gamma_n}$ is $\frac{n(n-1)}{2}$; the $d$-dimensional faces of
${\cal P}_{\Gamma_n}$ are $d$-dimensional polytopes of magic labelings
of positive graphs with $n$ vertices and at most $n+d$ edges. 

We can partition the vertices of $\Gamma_{2n}$ into two equal sets $A$
and $B$ in ${2n-1} \choose n$ ways: Fix the vertex $v_1$ to be in the
set $A$, then we can choose the n vertices for the set $B$ in ${2n-1}
\choose n$ ways, and the remaining $n-1$ vertices will belong to the
set $A$. By adding the required edges, we get a complete bipartite
graph for every such partition of the vertices of $\Gamma_{2n}$. Thus,
the number of subgraphs of $\Gamma_{2n}$ that are isomorphic to
$K_{n,n}$ is ${2n-1} \choose n$. Therefore, there are $2n-1 \choose n$
faces of ${\cal P}_{\Gamma_{2n}}$ that are Birkhoff polytopes because
every isomorphic copy of $K_{n,n}$ contributes to a face of ${\cal
P}_{\Gamma_{2n}}$.  $\square$

We now prove our results about polytopes of magic digraphs.

\noindent {\em Proof of Theorem \ref{digraphpolytopethm}.}  By Lemma
\ref{magicsum1}, all the elements of the Hilbert basis of $C_D$ have
magic sum 1. Since the extreme rays are a subset of the Hilbert basis
elements, it follows that the vertices of ${\cal P}_D$ are
integral. Since ${\cal P}_D = {\cal P}_{G_D}$, it follows by Theorem
\ref{graphpolytopethm} that the dimension of ${\cal P}_D$ is $q-2n+b$;
the $d$-dimensional faces of ${\cal P}_{D}$ are the $d$-dimensional
polytopes of magic labelings of positive subdigraphs of $D$ with $n$
vertices and at most $2n-b+d$ edges.  $\square$

We derive our results about the faces of the Birkhoff polytope as a
consequence.

\noindent {\em Proof of Theorem \ref{facesbirk}.}  The one-to-one
correspondence between semi-magic squares and magic labelings of
${\Pi_n}$ gives us that ${\cal P}_{\Pi_n} = B_n$. Since every digraph
with $n$ vertices is a subdigraph of $\Pi_n$, by Theorem
\ref{digraphpolytopethm}, it follows that its $d$-dimensional faces
are $d$-dimensional polytopes of magic labelings of positive digraphs
with $n$ vertices and at most $2n-1+d$ edges. Since the vertex set of
a face of $B_n$ is a subset of the vertex set of $B_n$ it follows that
the vertices of ${\cal P}_D$, where $D$ is a positive digraph, are
permutation matrices. $\square$

Our results enable us to reprove some known facts about the Birkhoff
polytope as well.  For example, Theorem \ref{digraphpolytopethm} gives
us that the dimension of $B_n$ is ${(n-1)}^2$. The leading coefficient
of the Ehrhart polynomial of $B_n$ is the volume of $B_n$. This number
has been computed for $n=1,2, \dots, 9$ (see \cite{beck} and
\cite{chanrobbins}).

%%%%%%%%%%%%%%%%%%%%%%%%%%%%%%%%%%%%%%%%%%%%%%%%%%%%%%%%%%%%%%%%%%%%%%%
%%%%%%%%%%%%%%%%%%%%%%%%%%%%%%%%%%%%%%%%%%%%%%%%%%%%%%%%%%%%%%%%%%%%%%%
%%%%%%%%%%%%%%%%%%%%%%%%%%%%%%%%%%%%%%%%%%%%%%%%%%%%%%%%%%%%%%%%%%%%%
%%%%%%%%%%%%%%%%%%%%%%%%%%%%%%%%%%%%%%%%%%%%%%%%%%%%%%%%%%%%%%%%%%%%%%
\section{Computational results}

We will now list our computational results. The numbers in square
brackets in the figures represent the number of elements in the orbit
class of each generator.

\subsection{Symmetric magic squares.}

\begin{figure}[h]
 \begin{center} 
 \includegraphics[scale=0.4]{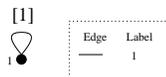}
\caption{The minimal Hilbert basis of $1 \times 1$ symmetric
magic squares.} \label{symm1hilb} \end{center}
\end{figure}

\begin{figure}[h]
 \begin{center} 
 \includegraphics[scale=0.4]{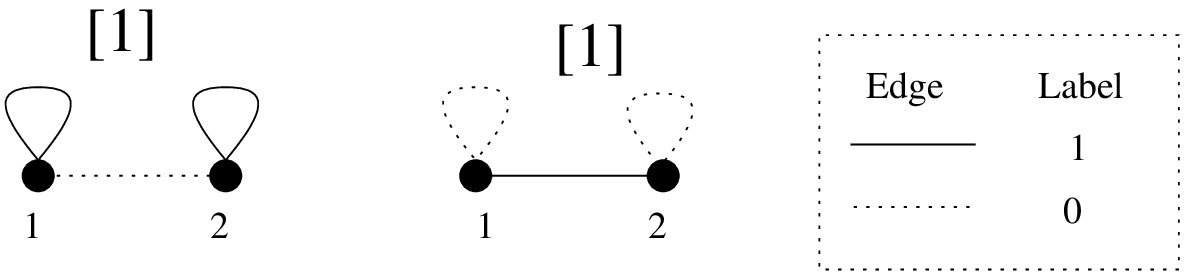}
\caption{Generators of the minimal Hilbert basis of $2 \times 2$  symmetric magic squares.} \label{symm2hilb} \end{center} 
\end{figure}

\begin{figure}[h]
 \begin{center} 
\includegraphics[scale=0.4]{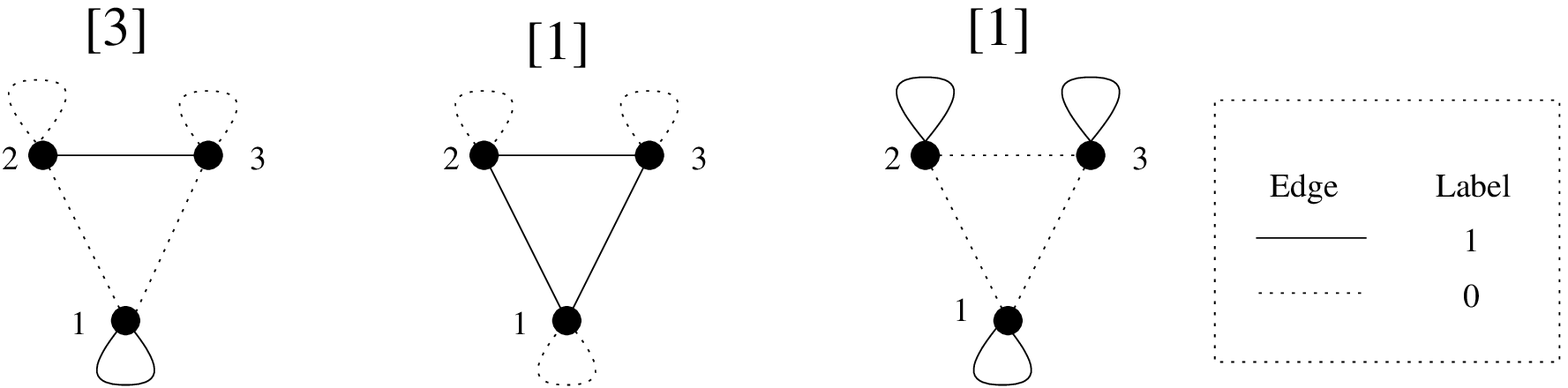}
\caption{Generators of the minimal Hilbert basis of $3 \times 3$ symmetric magic squares.} \label{symm3hilb}  \end{center} 
\end{figure}

\begin{figure}[h]
 \begin{center} 
 \includegraphics[scale=0.4]{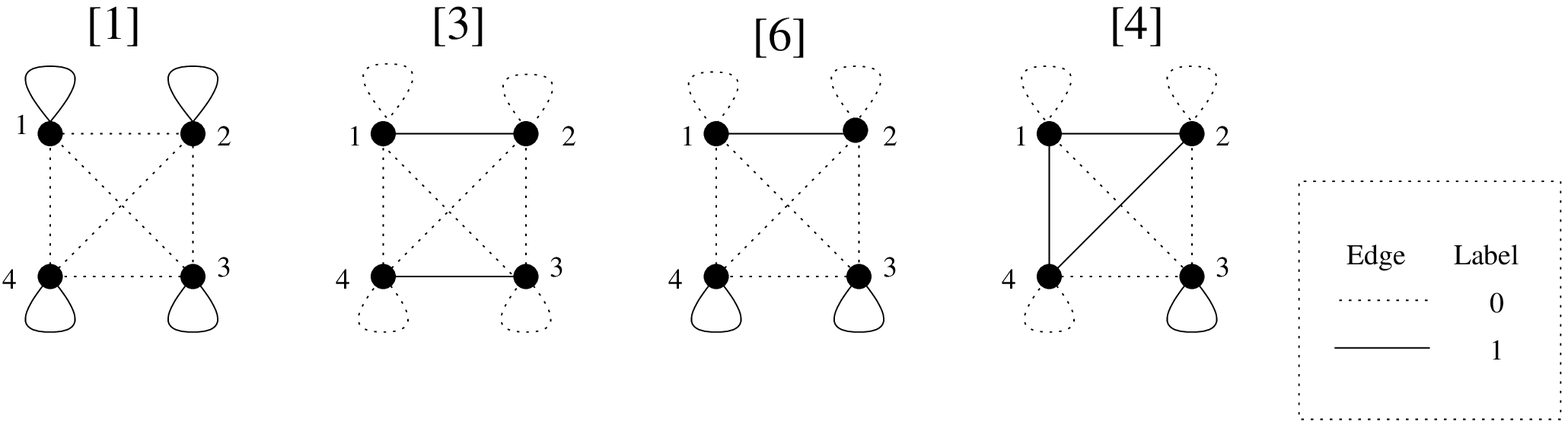}
\caption{Generators of the minimal Hilbert basis of $4 \times 4$ symmetric magic squares.} \label{symm4hilb}  \end{center} 
\end{figure}

\begin{figure}[h]
 \begin{center} 
 \includegraphics[scale=0.4]{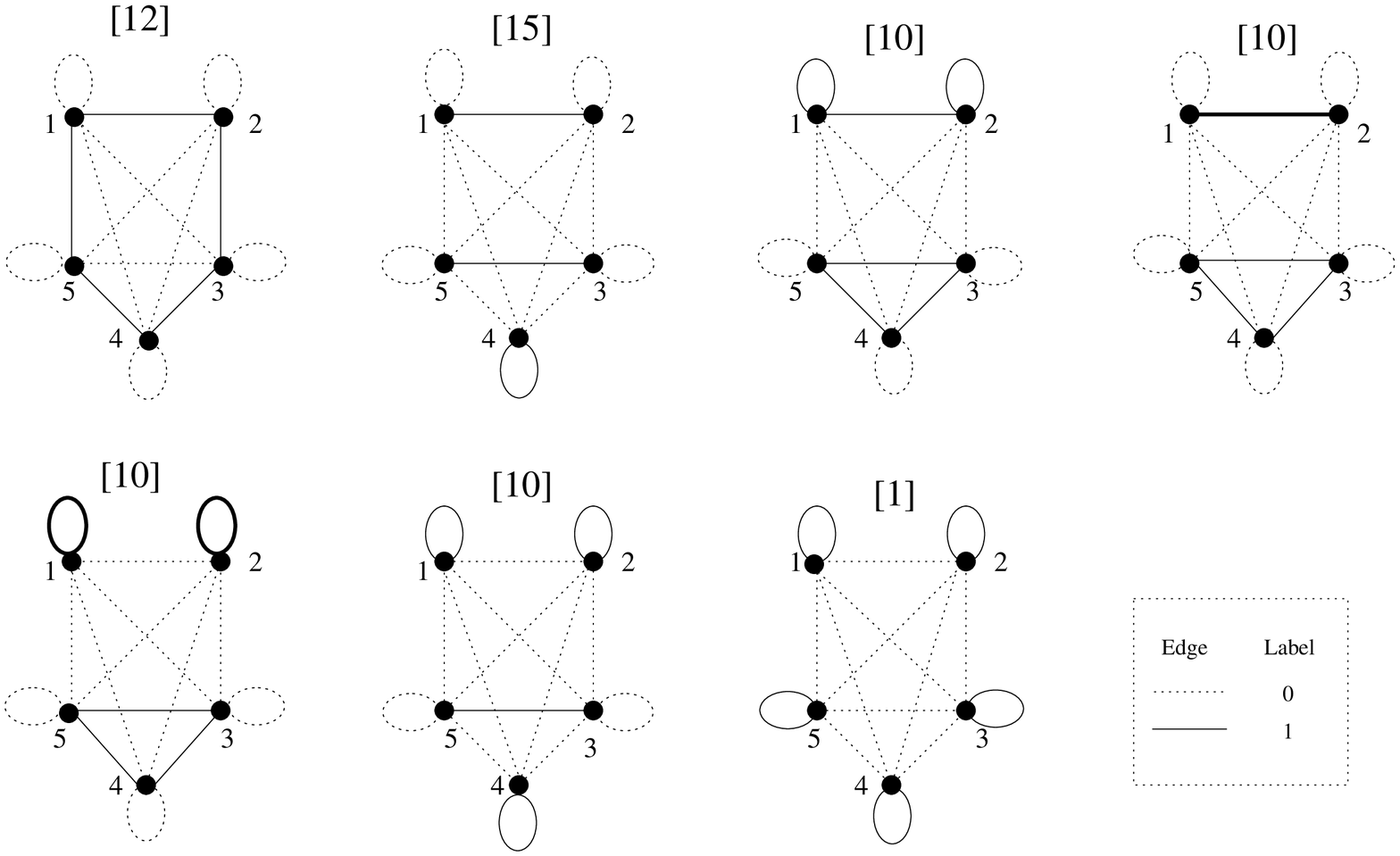}
\caption{Generators of the minimal Hilbert basis of $5 \times 5$ symmetric magic squares.}  \label{symm5hilb} \end{center} 
\end{figure}

\begin{figure}[h]
 \begin{center} 
 \includegraphics[scale=0.4]{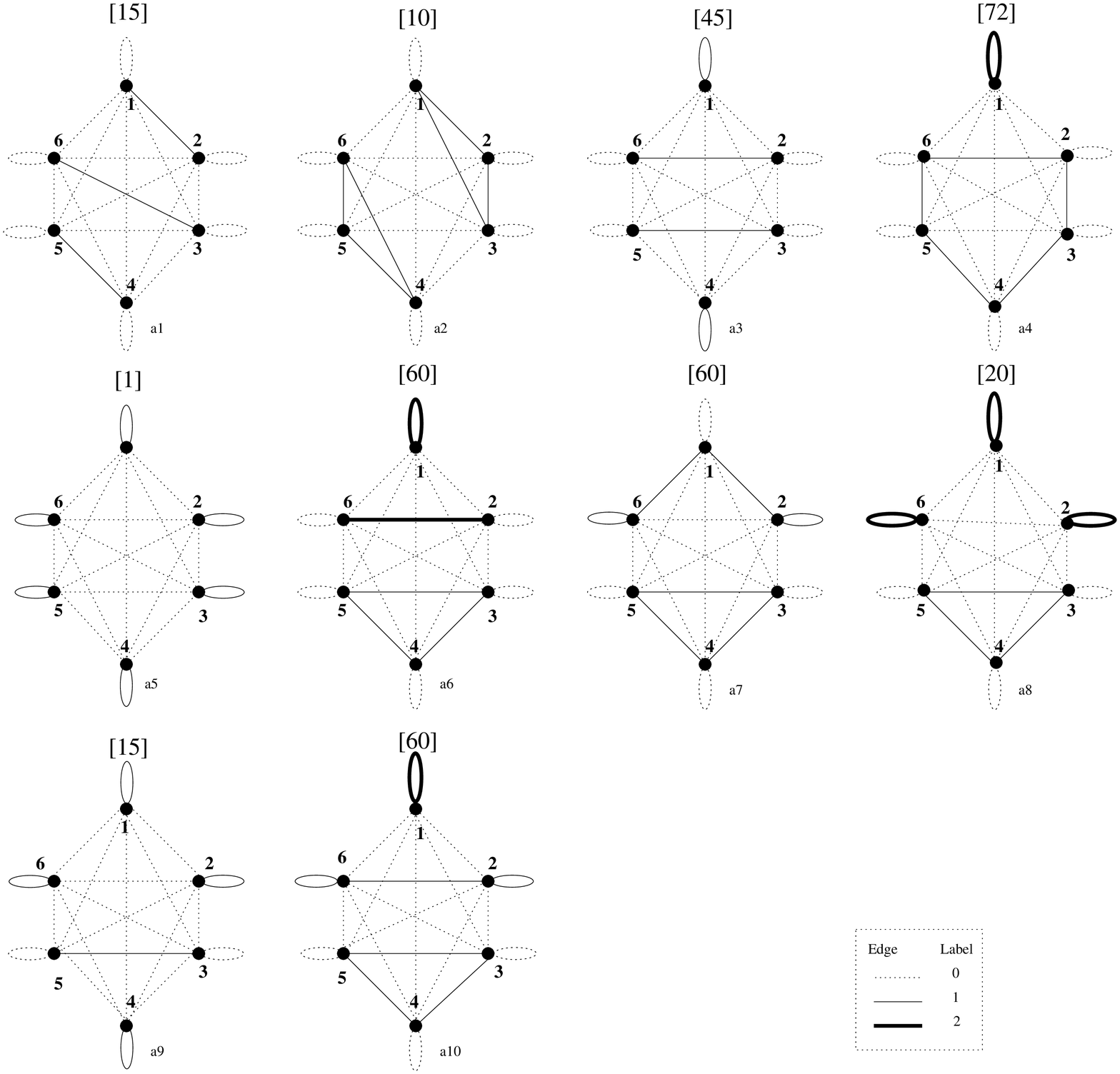}
\caption{Generators of the minimal Hilbert basis of $6 \times 6$ symmetric
magic squares.} \label{symm6hilb} \end{center}
\end{figure}

\begin{figure}[h]
 \begin{center} 
 \includegraphics[scale=0.4]{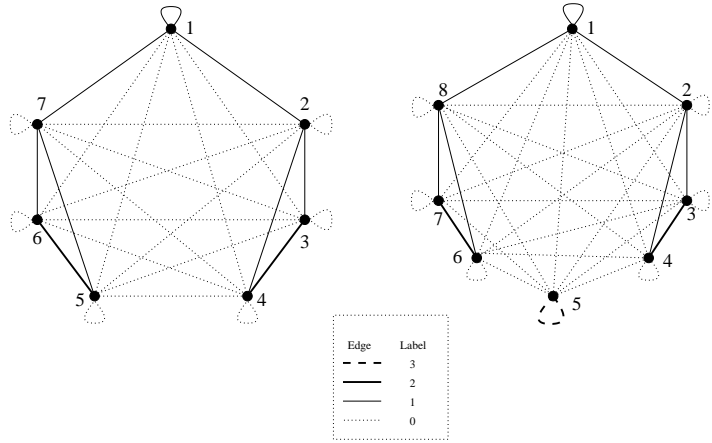}
\caption{Minimal Hilbert basis elements of $C_{\Gamma_7}$ and $C_{\Gamma_8}$
of magic sum 3.} \label{countersym} \end{center}
\end{figure}

The generators of the minimal Hilbert basis of $C_{\Gamma_n}$ for
$n=1,2,3,4,5$, and $6$ are given in Figures \ref{symm1hilb},
\ref{symm2hilb}, \ref{symm3hilb}, \ref{symm4hilb}, \ref{symm5hilb},
and \ref{symm6hilb}, respectively. It is interesting that all the
elements of the minimal Hilbert basis of $C_{\Gamma_n}$ for $n=1,
\dots, 6$ are 2-matchings. The minimal Hilbert basis elements are not,
in general, 2-matchings for all $n$ (see Figure \ref{countersym} for
examples of irreducible magic labelings of magic sum 3). In fact, it
follows from the results of chapter 11 in \cite{konig}, that there
exists a graph with an irreducible magic labeling of magic sum $r$ if
and only if $r$ is 2 or $r$ is odd.  By Proposition \ref{hilbGprop},
this implies, that there is a minimal Hilbert basis element of magic
sum $r$ of $C_{\Gamma_n}$, for some $n$, if and only if, $r$ is 2 or
$r$ is odd. The program 4ti2 was also able to compute the minimal
Hilbert bases of $C_{\Gamma_7}$ and $C_{\Gamma_8}$.

Recall that the number of $n \times n$ symmetric magic squares is the
same as $H_{\Gamma_n}(r)$ (the generating functions for
$H_{\Gamma_n}(r)$ for $n$ up to 5 are given in \cite{stanley4}). The
volume of $ {\cal P}_{\Gamma_n}$ is the leading coefficient of
$H_{\Gamma_n}(r)$.

%%%% Example %%%%%%%%%%%%%%%%%%%%%%%

%%%% Example %%%%%%%%%%%%%%%%%%%%%%%

%%%% Example %%%%%%%%%%%%%%%%%%%%%%%

%%%% Example %%%%%%%%%%%%%%%%%%%%%%%

%%%% Example %%%%%%%%%%%%%%%%%%%%%%%

%%%% Example %%%%%%%%%%%%%%%%%%%%%%%

%%%% Example %%%%%%%%%%%%%%%%%%%%%%%

\[
H_{\Gamma_1}(r)  = 1 \mbox{ for all $r \geq 0$.}
\]

\[
H_{\Gamma_2}(r)  = r+1  \mbox{ for all $r \geq 0$.}
\]

\[
H_{\Gamma_3}(r)  = \left \{
\begin{array}{ll}
\frac{1}{4}r^3+ \frac{9}{8}r^2+ \frac{7}{4}r+1
&  \mbox{if 2 divides $r$,} \\ \\ 
\frac{1}{4}r^3+ \frac{9}{8}r^2+ \frac{7}{4}r+ \frac{7}{8}
 & \mbox{ otherwise.}
\end{array}
\right .
\]

\[
H_{\Gamma_4}(r)  = \left \{
\begin{array}{ll}
{\frac {1}{72}}{r}^{6}+\frac{1}{6}{r}^{5}+{\frac {119}{144}}{r}^{4}+{
\frac {13}{6}}{r}^{3}+{\frac {29}{9}}{r}^{2}+\frac{8}{3}r+1
&  \mbox{if 2 divides $r$,} \\ \\ 
{\frac {1}{72}}{r}^{6}+\frac{1}{6}{r}^{5}+{\frac {119}{144}}{r}^{4}+{
\frac {13}{6}}{r}^{3}+{\frac {29}{9}}{r}^{2}+\frac{8}{3}r+{\frac {15}{16
}}
 & \mbox{ otherwise.}
\end{array}
\right .
\]

\[
H_{\Gamma_5}(r)   = \left \{
\begin{array}{l}
{\frac {365}{3096576}}{r}^{10}+{\frac {9125}{3096576}}{r}^{9}+{
\frac {22553}{688128}}{r}^{8}+{\frac {55085}{258048}}{r}^{7}+{
\frac {11083}{12288}}{r}^{6}+{\frac {7945}{3072}}{r}^{5}+{\frac {
1978913}{387072}}{r}^{4} \\ \\ +{\frac {335065}{48384}}{r}^{3} +{\frac
{ 50329}{8064}}{r}^{2}+{\frac {1177}{336}}r+1 \\
\hfill   \mbox{if 2 divides $r$,} \\ \\  \\

{\frac {365}{3096576}}{r}^{10}+{\frac {9125}{3096576}}{r}^{9}+{
\frac {22553}{688128}}{r}^{8}+{\frac {55085}{258048}}{r}^{7}+{
\frac {11083}{12288}}{r}^{6}+{\frac {63545}{24576}}{r}^{5}+{\frac 
{15807679}{3096576}}{r}^{4} \\ \\ +{\frac {5329855}{774144}}{r}^{3} 
+{
\frac {6327137}{1032192}}{r}^{2}+{\frac {1139917}{344064}}r+{
\frac {27213}{32768}} \\

 \hfill  \mbox{ otherwise.}
\end{array}
\right .
\]

\vspace{0.2in}

\subsection{Pandiagonal symmetric magic squares.}

\begin{figure}[h]
 \begin{center} 
 \includegraphics[scale=0.4]{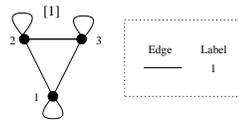}
\caption{Generators of the minimal Hilbert basis of $3 \times 3$
pandiagonal symmetric magic squares.}  \label{pansymm3hilb} \end{center}
\end{figure}

\begin{figure}[h]
 \begin{center} 
 \includegraphics[scale=0.4]{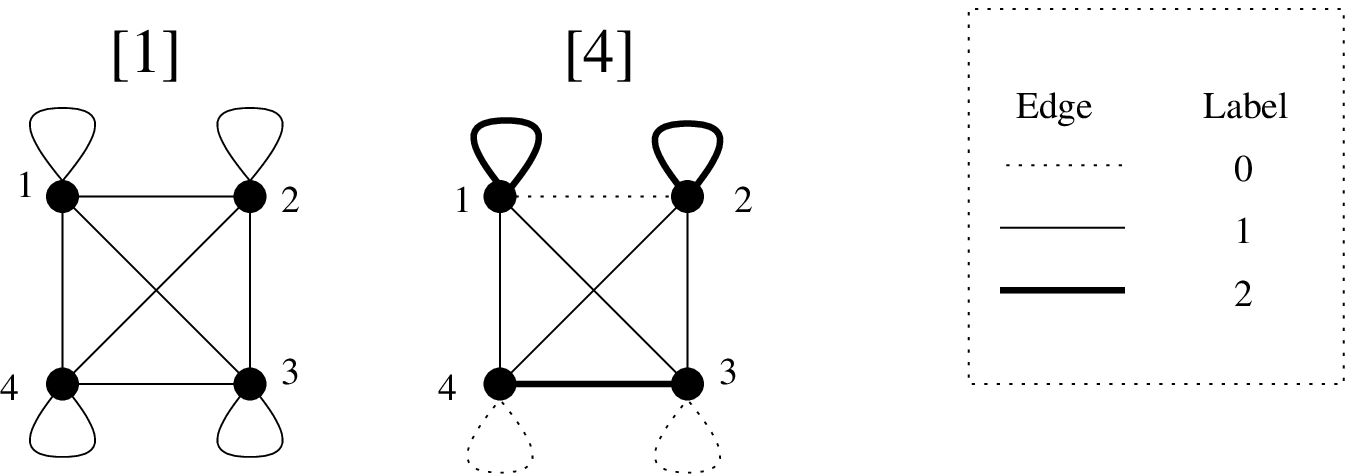}
\caption{Generators of the minimal Hilbert basis of $4 \times 4$
pandiagonal symmetric magic squares.}  \label{pansymm4hilb} \end{center}
\end{figure}

\begin{figure}[h]
 \begin{center} 
 \includegraphics[scale=0.4]{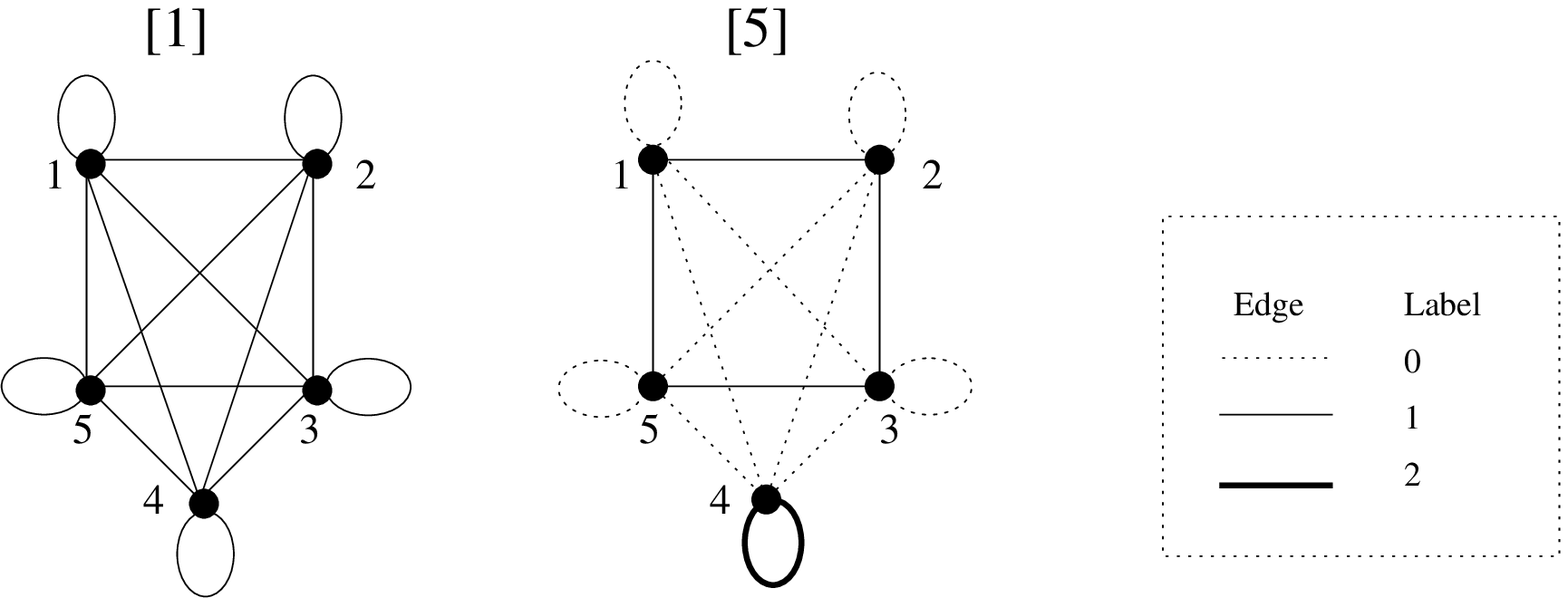}
\caption{Generators of the minimal Hilbert basis of $5 \times 5$
pandiagonal symmetric magic squares.}  \label{pansymm5hilb} \end{center}
\end{figure}

{\em Pandiagonal symmetric magic squares} are symmetric magic squares
such that all the pandiagonals also add to the magic sum (see Figure
\ref{pans}).  The generators of the Hilbert basis of $n \times n$ pandiagonal
symmetric magic squares for $n=3$, $4$, and  $5$ are given in Figures
\ref{pansymm3hilb}, \ref{pansymm4hilb}, and \ref{pansymm5hilb},
respectively.  The Hilbert basis of $6 \times 6$ pandiagonal symmetric
magic squares contain 4927 elements and can be computed using the
program 4ti2.

%%%% Example %%%%%%%%%%%%%%%%%%%%%%%

%%%% Example %%%%%%%%%%%%%%%%%%%%%%%

%%%% Example %%%%%%%%%%%%%%%%%%%%%%%

Let $P_n(r)$ denote the number of $n \times n$ pandiagonal symmetric magic squares
with magic sum $r$. We derive:

\[
P_3(r)  = 
 \left \{
\begin{array}{ll}
1 &  \mbox{ if 3 divides $r$,} \\ \\
0 & \mbox{ otherwise.}
\end{array}
\right .
\]

\[
P_4(r)  =  \left \{
\begin{array}{ll}
\frac{1}{8}r^2 + \frac{1}{2}r + 1 & 
\mbox{if 4 divides $r$,} \\ \\  
0 & \mbox{ otherwise.}
\end{array}
\right . 
\]

\[
P_5(r)   = \left \{
\begin{array}{ll}
\frac{1}{384}t^4 + \frac{5}{96}t^3 + \frac{35}{96}t^2 + \frac{25}{24}t + 1 &  \mbox{if 2 divides $r$,} \\ \\  
\frac{1}{384}t^4 - \frac{5}{192}t^2 + \frac{3}{128}
& \mbox{ otherwise.}
\end{array}
\right .
\]

%%%%%%%%%%%%%%%%%%%%%%%%%%%%%%%%%%%%%%%%%%%%%%%%%%%%%%%%%%%%%%%%%%%
\subsection{Magic labelings of Complete Graphs.}

The minimal Hilbert basis of the cone of magic labelings of the
complete graph $K_n$ corresponds to the set of elements of the minimal
Hilbert basis of $C_{\Pi_n}$ for which all the loops are labeled with a
0. ${\cal P}_{K_n}$ is an $\frac{n(n-3)}{2}$ dimensional polytope with
the description:
\[
{\cal P}_{K_n} = \{ X \in {\cal P}_{\Gamma_n} | x_{ii} = 0, i = 1, \dots,  n \}.
\]
\[
H_{K_1}(r) = 0
\]
\[
H_{K_2}(r) =  1
\]
\[
H_{K_3}(r) =  \left \{
\begin{array}{ll}
1 &  \mbox{if 2 divides $r$,} \\ \\ 
0 & \mbox{ otherwise.}
\end{array}
\right .
\]

\[
H_{K_4}(r) = \frac{1}{2}r^2 + \frac{3}{2}r + 1.
\]

\[
H_{K_5}(r) = \left \{
\begin{array}{ll}
{\frac {5}{256}}{r}^{5}+{\frac {25}{128}}{r}^{4}+{\frac {155}{192}
}{r}^{3}+{\frac {55}{32}}{r}^{2}+{\frac {47}{24}}r+1 
&  \mbox{if 2 divides $r$,} \\ \\ 
0 & \mbox{ otherwise.}
\end{array}
\right .
\]

\[
H_{K_6}(r) = \left \{
\begin{array}{l}
{\frac {19}{120960}}{r}^{9}+{\frac {19}{5376}}{r}^{8}+{\frac {143}
{4032}}{r}^{7}+{\frac {5}{24}}{r}^{6}+{\frac {4567}{5760}}{r}^{5
}+{\frac {785}{384}}{r}^{4}+{\frac {10919}{3024}}{r}^{3} \\ \\ +{\frac {
955}{224}}{r}^{2}+{\frac {857}{280}}r+1 \\ 
\hfill \mbox{if 2 divides $r$,} \\ \\ \\ 
{\frac {19}{120960}}{r}^{9}+{\frac {19}{5376}}{r}^{8}+{\frac {143}
{4032}}{r}^{7}+{\frac {5}{24}}{r}^{6}+{\frac {4567}{5760}}{r}^{5
}+{\frac {785}{384}}{r}^{4}+{\frac {10919}{3024}}{r}^{3} \\ \\ +{\frac {
955}{224}}{r}^{2}+{\frac {857}{280}}r+{\frac {251}{256}}
\\
\hfill \mbox{otherwise.}
\end{array}
\right .
\]

See \cite{stewartcomplete} for more aspects of magic labelings of 
complete graphs $K_n$. 

\subsection{Magic labelings of the Petersen graph.}
\begin{figure}[h]
 \begin{center} 
\includegraphics[scale=0.3]{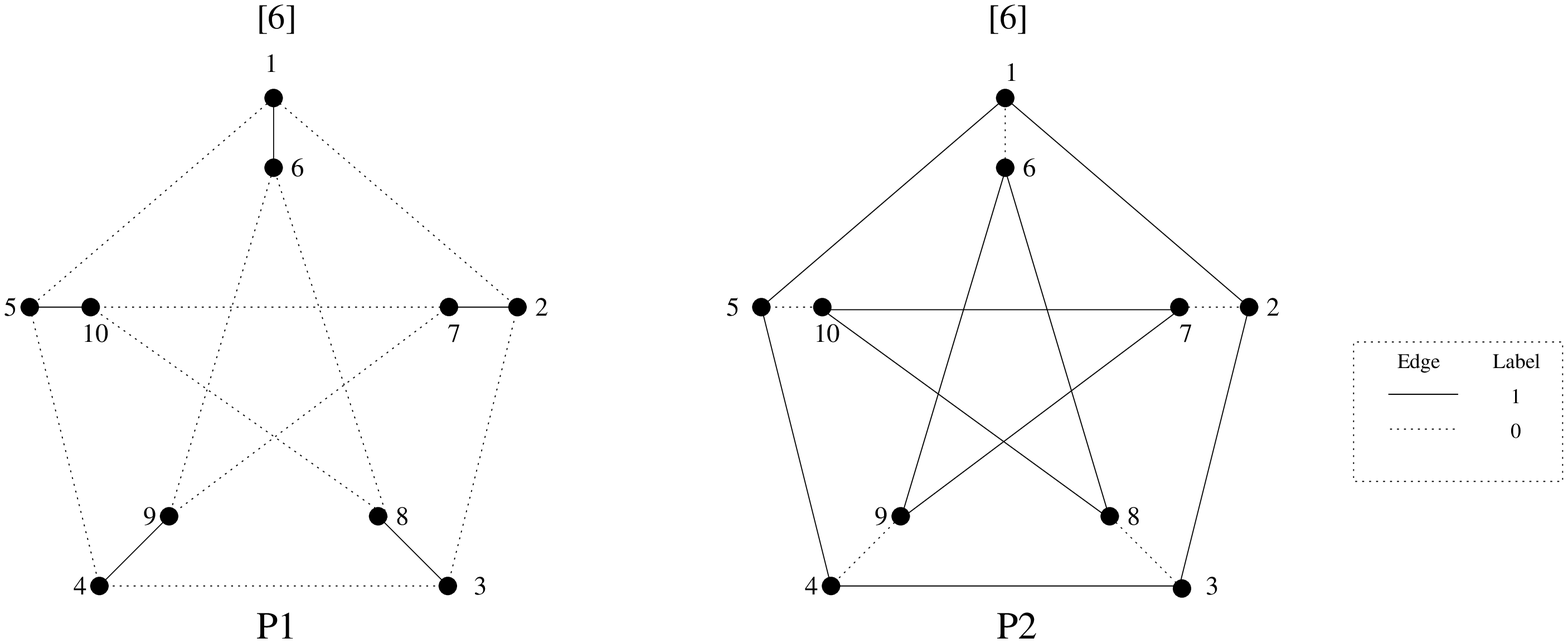}
\caption{Generators of the minimal Hilbert basis of magic labelings of the Petersen Graph.}
\label{petersenhilbert}
 \end{center} \end{figure}

The generators of the Hilbert basis are given in Figure
\ref{petersenhilbert} (numbers in square brackets are the number of
elements in the orbit of the generators).

%%%% Example %%%%%%%%%%%%%%%%%%%%%%%

Let $H_{Petersen}(r)$ denote the number of magic labelings of the
Petersen graph with magic sum $r$. The generating function $F(t)$ for
the Petersen graph (F(t) is also derived in \cite{stanley4}) is:
\[
F(t) = {\frac {{t}^{4}+{t}^{3}+6\,{t}^{2}+t+1}{ \left( 1-t \right)
 ^{6} \left( 1+t \right) }} =
 1+6t+27{t}^{2}+87{t}^{3}+228{t}^{4}+513{t}^{5}+1034{t}^{6
 }+1914{t}^{7}+ \dots
\]

Therefore, we get

\[
H_{Petersen}(r) = \left \{
\begin{array}{ll}

\frac{1}{24}{r}^{5}+{\frac {5}{16}}{r}^{4}+{\frac {25}{24}}{r}^{3}+{
\frac {15}{8}}{r}^{2}+{\frac {23}{12}}r+1
&  \mbox{if 2 divides $r$,} \\ \\ 

\frac{1}{24}{r}^{5}+{\frac {5}{16}}{r}^{4}+{\frac {25}{24}}{r}^{3}+{
\frac {15}{8}}{r}^{2}+{\frac {23}{12}}r+{\frac {13}{16}}

& \mbox{ otherwise.}
\end{array}
\right .
\]

\subsection{Magic labelings of the Platonic graphs.}

\subsubsection{Magic labelings of the Tetrahedral Graph.}

The minimal Hilbert basis of the cone of magic labelings of the
tetrahedral graph is given in Figure \ref{tetrahilb}. Since all the
elements of the minimal Hilbert basis have magic sum 1, it follows
that the vertices of the polytope of magic labelings of the
tetrahedral graph are integral points. Therefore, we get that
$H_{tetrahedral}(r)$ is a polynomial, where $H_{tetrahedral}(r)$
denotes the number of magic labelings of the Tetrahedral graph.
We derive an explicit formula.
\[
H_{tetrahedral}(r) = \frac{1}{2}r^2 + \frac{3}{2}r + 1.
\]

Theorem \ref{bippoly} states that if a graph $G$ is bipartite then
$H_G(r)$ is a polynomial. Thus, the tetrahedral graph is an example
that proves that $H_G(r)$ being a polynomial does not imply that $G$
is bipartite. 

\subsubsection{Magic labelings of the Cubical graph.}

\begin{figure}[h]
 \begin{center} 
 \includegraphics[scale=0.4]{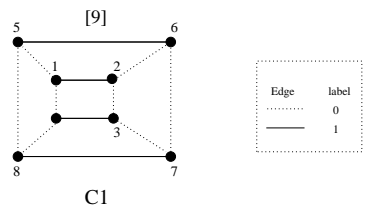}
\caption{Generator of the minimal Hilbert basis of magic labelings of the
Cubical Graph.} \label{cubehilb}
\end{center} \end{figure}

The minimal Hilbert basis of the cone of the magic labelings of the cubical
graph is the set consisting of C1 in Figure \ref{cubehilb} and the
eight magic labelings isomorphic to C1.

%%%% Example %%%%%%%%%%%%%%%%%%%%%%%

Let $H_{cube}(r)$ denote the number of magic labelings of the Cubical
graph of magic sum $r$. Observe that the cubical graph is bipartite,
therefore Theorem \ref{bippoly} applies, and $H_{cube}(r)$ is a
polynomial.  $H_{cube}(r)$ is also derived in \cite{stanley3}.

\[
H_{cube}(r) = \frac{1}{15}r^5 + \frac{1}{2}r^4 + \frac{5}{3}r^3 + 3r^2 
+ \frac{83}{30}r + 1   
\]

\subsubsection{The magic labelings of the Octahedral graph.}

There are 12 elements in the minimal Hilbert basis of the cone of
magic labelings of the Octahedral graph: the 8 magic labelings in the
$S_6$ orbit of the magic labeling O1, and the four magic labelings in
the orbit of O2 (see Figure \ref{octahilb}). Let $H_{octahedral}$
denote the number of magic labelings of the Octahedral graph of magic
sum $r$. The generating function of $H_{octahedral}(r)$ is given in
\cite{stanley3}. 
\[
H_{octahedral}(r) = \left \{
\begin{array}{ll}
{\frac {1}{120}}{r}^{6}+ \frac{1}{10}{r}^{5}+{\frac {25}{48}}{r}^{4}
+\frac{3}{2}{r}^{3}+{\frac {38}{15}}{r}^{2}+{\frac {12}{5}}r+1
&  \mbox{if 2 divides $r$,} \\ \\ 

{\frac {1}{120}}{r}^{6}+ \frac{1}{10}{r}^{5}+{\frac {25}{48}}{r}^{4}+
\frac{3}{2}{r}^{3}+{\frac {38}{15}}{r}^{2}+{\frac {12}{5}}r+{\frac {15}{16}
}

& \mbox{ otherwise.}
\end{array}
\right .
\]

\subsubsection{The magic labelings of the Dodecahedral graph.}

\begin{figure}[h]
% \begin{center} 
\includegraphics[scale=0.25]{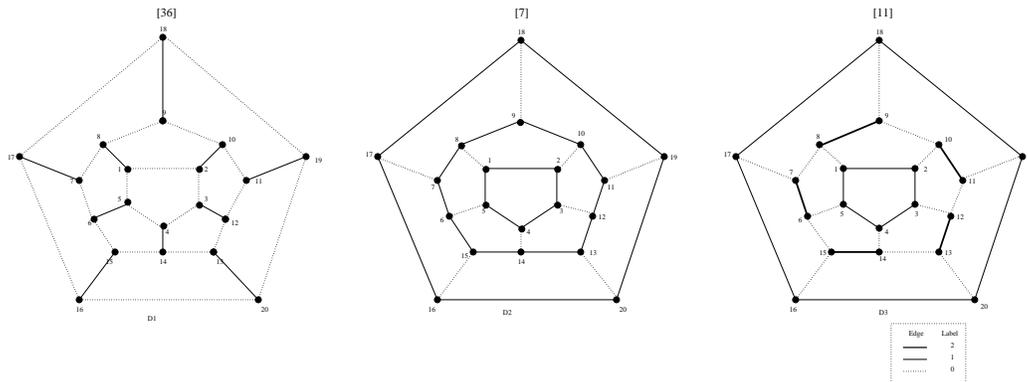}
\caption{Generators of the minimal Hilbert basis of magic labelings of the
dodecahedral Graph.} \label{dodecahilb}   
%\end{center} 
\end{figure}

The generators of the minimal Hilbert basis of the Dodecahedral graph are
given in Figure \ref{dodecahilb}.

%%%% Example %%%%%%%%%%%%%%%%%%%%%%%

Let $H_{dodecahedral}(r)$ denote the number of magic labelings of the 
Dodecahedral graph with magic sum $r$. We derive: 

\[
H_{dodecahedral}(r) = \left \{
\begin{array}{l}

{\frac {47}{40320}}{r}^{10}+{\frac {47}{2688}}{r}^{9}+{\frac {225}
{1792}}{r}^{8}+{\frac {9}{16}}{r}^{7}+{\frac {3361}{1920}}{r}^{6
}+{\frac {255}{64}}{r}^{5}+{\frac {27625}{4032}}{r}^{4} \\ \\ +{\frac {
1513}{168}}{r}^{3}+{\frac {4691}{560}}{r}^{2} 
+\frac{9}{2}r+1

\\ \hfill  \mbox{if 2 divides $r$,} \\ \\

{\frac {47}{40320}}{r}^{10}+{\frac {47}{2688}}{r}^{9}+{\frac {225}
{1792}}{r}^{8}+{\frac {9}{16}}{r}^{7}+{\frac {3361}{1920}}{r}^{6
}+{\frac {255}{64}}{r}^{5}+{\frac {27625}{4032}}{r}^{4} \\ \\ +{\frac {
1513}{168}}{r}^{3}+{\frac {4691}{560}}{r}^{2} +{\frac {567}{128}}
r+{\frac {229}{256}}

 \\  \hfill \mbox{ otherwise.}
\end{array}
\right .
\]

\subsubsection{The magic labelings of the Icosahedral graph.}

\begin{figure}[h]
 \begin{center} 
 \includegraphics[scale=0.4]{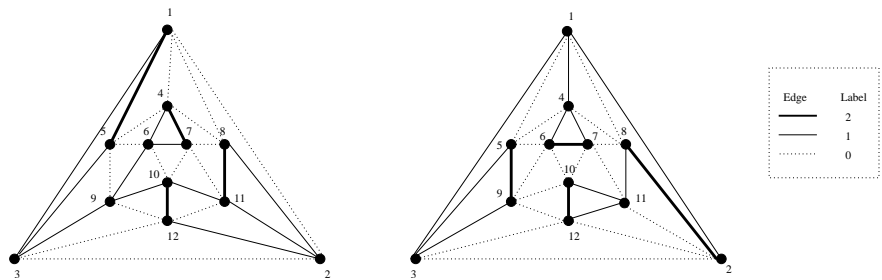}
\caption{Icosahedral graph : minimal Hilbert basis elements of magic sum 3.} 
\label{icohilb}
 \end{center}  \end{figure}

There are 4195 elements in the minimal Hilbert basis of the cone of
magic labelings of the Icosahedral graph which can be computed using
4ti2. It is interesting that unlike the other platonic graphs, all the
minimal Hilbert basis elements are not two-matchings (see Figure
\ref{icohilb}). The formula for the number of magic labelings of the
Icosahedral graph remains unresolved.

%%%% Example %%%%%%%%%%%%%%%%%%%%%%%

\appendix 
\newpage
\pagestyle{myheadings} 
\markright{  \rm \normalsize Appendix A \hspace{0.5cm}}
%Appendix------
\chapter{}
\label{appendixa}

\section{Proof of the minimal Hilbert basis Theorem.}

The {\em cone generated} by a set $X$ of vectors is the smallest cone
containing $X$ and is denoted by cone $X$; so
\[
\mbox{cone } X = \{ {\lambda}_1x_1+....+ {\lambda}_kx_k | k \geq 0;
x_1,\dots, x_k \in X;  {\lambda}_1, \dots, {\lambda}_k \geq 0 \}.
\]

{\em Proof of Theorem \ref{hilbtheorem}.}

Let $C$ be a rational polyhedral cone, generated by $b_1,b_2,...,b_k$.
Without loss of generality $b_1,b_2,...,b_k$ are integral vectors. Let
$a_1,a_2,...,a_t$ be all the integral vectors in the polytope $\cal
P$:
\[
{\cal P} = \{ {\lambda}_1b_1+....+ {\lambda}_kb_k | 
0 \leq  {\lambda}_i \leq 1 \,(i=1,..,k) \}
\]

Then $a_1,a_2,...,a_t$ generate $C$ as $b_1,b_2,...,b_k$ occur among
$a_1,a_2,...,a_t$ and as $\cal P$ is contained in $C$. We will now
show that $a_1,a_2,...,a_t$ also form a Hilbert basis. Let $b$ be an
integral vector in $C$.  Then there are ${\mu}_1,{\mu}_2,...,{\mu}_k
\geq 0$ such that
\begin{eqnarray} \label{bdecomp}
b = {\mu}_1b_1+{\mu}_2b_2+ \cdots +{\mu}_kb_k.
\end{eqnarray}
Then
\[
b = {\lfloor{\mu}_1\rfloor}b_1+{\lfloor{\mu}_2\rfloor}b_2+ \cdots
+{\lfloor{\mu}_k\rfloor}b _k+({\mu}_1 -
{\lfloor{\mu}_1\rfloor})b_1+({\mu}_2 - {\lfloor{\mu}_2\rfloor})b_2+
\cdots + ({\mu}_k - {\lfloor{\mu}_k\rfloor})b_k.
\]

Now the vector
\begin{eqnarray} \label{b}
b -{\lfloor{\mu}_1\rfloor}b_1- \dots -{\lfloor{\mu}_k\rfloor}b _k =
({\mu}_1 - {\lfloor{\mu}_1\rfloor})b_1+ \dots + ({\mu}_k -
{\lfloor{\mu}_k\rfloor})b_k
\end{eqnarray}

occurs among $a_1,a_2,...,a_t$ as the left side of the Equation
\ref{b} is clearly integral and the right side belong to $\cal P$.
Since also $b_1,b_2,...,b_k$ occur among $a_1,a_2,...,a_t$, it follows
that \ref{bdecomp} decomposes $b$ as a nonnegative integral
combination of $a_1,a_2,...,a_t$. So $a_1,a_2,...,a_t$ form a Hilbert
basis.

Next suppose $C$ is pointed. Consider $H$ the set of all irreducible
integral vectors. Then it is clear that any Hilbert basis must contain
$H$. So $H$ is finite because it is contained in $\cal P$. To see that
$H$ itself is a Hilbert basis generating $C$, let $b$ be a vector such
that $bx > 0$ if $x \in C \backslash \{ 0 \}$ ($b$ exists because $C$
is pointed). Suppose not every integral vector in $C$ is a nonnegative
integral combination of vectors in $H$.  Let $c$ be such a vector,
with $bc$ as small as possible (this exists, as $c$ must be in the set
$\cal P$). As $c$ is not in $H$, $c = c_1+c_2$ for certain nonzero
integral vectors $c_1$ and $c_2$ in $C$. Then $bc_1 < bc$ and  $bc_2 < bc$.
Therefore $c_1$ and $c_2$ are nonnegative integral combinations of vectors
in $H$, and therefore $c$ is also. $\square$ 

\section{Proof of the  Hilbert-Serre Theorem.}

Let $C$ be a class of $A$-modules and let $H$ be a function on $C$
with values in $\integers$. The function $H$ is called {\em additive}
if for each short exact sequence
\[
0 \rightarrow M^{\prime}  \stackrel{f}{\rightarrow} M  \stackrel{g} {\rightarrow} M^{\prime \prime} \rightarrow 0
\]

in which all the terms belong to $C$, we have 
\[
H( M^{\prime}) - H(M) + H(M^{\prime \prime}) = 0.
\] 

\begin{prop}[proposition 2.11, \cite{atiyah}]
Let $0 \rightarrow M_0 \rightarrow M_1 \rightarrow \dots \rightarrow
M_n \rightarrow 0$ be an exact sequence of $A$-modules in which all
the modules $M_i$ and the kernels of all the homomorphisms belong to
$C$. Then for any additive function $H$ on $C$ we have
\[
\sum_{i=0}^n (-1)^i H(M_i) = 0.
\]

\end{prop}

{\em Proof.} The proof follows because every exact sequence can be split
into short exact sequences (see \cite{atiyah}, Chapter 2). $\square$

For any $A$-module homomorphism $\phi$ of $M$ into $N$, we have an 
an exact
sequence, 
\[
0  {\rightarrow} \mbox{ker}(\phi) {\rightarrow} M  \stackrel{\phi}{\rightarrow} N  {\rightarrow} \mbox{coker}(\phi)  {\rightarrow} 0,
\]

where $\mbox{ker}(\phi) {\rightarrow} M$ is the inclusion map and 
$ N  {\rightarrow} \mbox{coker}(\phi)=  N/{\mbox{im}(\phi)}$
is the natural homomorphism onto the quotient module \cite{dummit}.

\begin{thm}[Theorem 11.1 \cite{atiyah},(Hilbert,Serre)] Let $A =
\bigoplus_{n=0}^{\infty} A_n$ be a graded Noetherian ring. Let $A$ be
generated as a $A_0$-algebra by say $x_1,x_2,...,x_s$, which are
homogeneous of degrees $k_1,k_2,..,k_s$ (all $> 0$). Let $H$ be an
additive function on the class of all finitely-generated
$A_0$-modules. Let $M$ be a finitely generated $A$-module.  Then the
Hilbert-Poincar\'e series of $M$, $H_{M}(t) =
\sum_{n=0}^{\infty}H(M_n)t^n$ is a rational function in $t$ of the
form ${p(t)}/{\Pi_{i=1}^s(1-t^{k_i})}$, where $p(t) \in {\integers}[t]
$.
\end{thm}

{\em Proof.} Let $M = \bigoplus M_n$, where $M_n$ are the graded
components of $M$, then  $M_n$ is finitely generated as a $A_0$-module.
The proof of the theorem is by induction on $s$, the number of
generators of $A$ over $A_0$. Start with $s=0$; this means that  $A_n = 0$
for all $n > 0$, so that $A = A_0$, and $M$ is a finitely-generated $A_0$ 
module, hence $M_n=0$ for all large $n$. Thus  $H_{M}(t)$ is a polynomial
in this case.

Now suppose $s > 0$ and the theorem true for $s-1$. Multiplication by $x_s$ is
an $A$-module homomorphism  of $M_n$ into $M_{n+k_s}$, hence it gives an exact
sequence, say
\begin{eqnarray} \label{shortexact}
0  {\rightarrow} K_n {\rightarrow} M_n  \stackrel{x_s}{\rightarrow} M_{n+k_s}  {\rightarrow} L_{n+k_s}  {\rightarrow} 0.
\end{eqnarray}

$K = \bigoplus_n K_n$, $L = \bigoplus_n L_n$ are both finitely
generated $A$-modules and both are annihilated by $x_s$, hence they are
$A_0[x_1, \dots, x_{s-1}]$-modules. Applying $H$ to \ref{shortexact}
we have
\[
H(K_n) - H(M_n) + H(M_{n+k_s}) - H(L_{n+k_s}) = 0;
\]

multiplying by $t^{n+k_s}$ and summing with respect to $n$ we get
\[
(1-t^{k_s})H(M,t) = H(L,t) - t^{k_s}H(K,t) + g(t),
\]

where $g(t)$ is a polynomial. Applying the inductive hypothesis the
result now follows. $\square$

\newpage
\pagestyle{myheadings} 
\markright{  \rm \normalsize Appendix B \hspace{0.5cm}}
%Appendix------
\chapter{}
\label{appendixb}

In this chapter, we provide some basic algorithms to compute Hilbert bases,
Hilbert-Poincar\'e series, and toric ideals. A knowledge of Gr\"obner
bases is assumed. An excellent introduction to Gr\"obner bases is
given in \cite{coxlittleoshea}. Many available computer algebra packages
(for example Maple and CocoA) can compute Gr\"obner bases.   
 
\section{Algorithms to compute Hilbert bases.}

We describe Algorithm 1.4.5 in \cite{sturmfels} to compute the Hilbert basis
of a cone $C_A= \{ {\bf x} : A{\bf x}= 0, {\bf x} \geq 0 \} $.

Let A be an $m \times n$ matrix. We introduce $2n+m$ variables
$t_1,t_2,..t_m$, $x_1,..,x_n$, $y_1,y_2,..,y_n$ and fix any
elimination monomial order  such that
$$\{ t_1,t_2,..t_m \} > \{x_1,..,x_n \} > \{y_1,y_2,..,y_n\}.$$

Let $I_A$ denote the kernel of the map
\[
\complex[x_1,\dots ,x_n,y_1, \dots ,y_n]
\rightarrow \complex[t_1, \dots, t_m, t_1^{-1}, \dots, t_m^{-1},y_1, \dots ,y_n],  \]

\[
x_1 \rightarrow y_1 \prod_{j=1}^m t_j^{a_{1j}}, \dots, 
x_n \rightarrow y_n \prod_{j=1}^m t_j^{a_{nj}},
y_1 \rightarrow y_1, \dots, y_n \rightarrow y_n.
\]

We can compute a Hilbert basis of  $C_A$ as follows.
\begin{algo} [Algorithm 1.4.5, \cite{sturmfels}]
\end{algo}

1. Compute the reduced Gr\"obner basis $\cal G$ with respect to $<$
for the ideal $I_A$.

2. The Hilbert basis of $C_A$ consists of all vectors $\beta$ such that
$x^{\beta} - y^{\beta}$ appears in $\cal G$.

For example, let
\[
A = \left[
\begin{array}{ccc}
1 & -1 \\
-2 & 2 
\end{array}
\right]
\]

To handle computations with negative exponents we introduce a new
variable $t$ and consider the lexicographic ordering
\[
t > t_1 > t_2 > x_1 > x_2 > y_1 > y_2.
\]

We compute the  Gr\"obner basis of     
$I_A = (x_1-y_1t_1^3t^2, \,  x_2-y_2t_2^3t, \,  t_1t_2t-1)$ with respect to 
the above ordering and get:
\[ I_A = 
(\underline{x_1x_2-y_1y_2}, \, t_1y_1-t_2^2x_1, \, t_1x_2-t_2^2y_2, \,
t_2^3ty_2-x_2, \, t_2^3tx_1-y_1, \, t_1t_2t-1)
\]

Therefore, the Hilbert basis is $\{ (1,1)\}$. 

See \cite{sturmfels} for more details about this algorithm.  See
\cite{raymond} for more effective algorithms to compute the Hilbert
basis.

\section{ Algorithms to compute toric ideals.}
Computing  toric ideals is the biggest challenge we face in applying
the methods we developed in this thesis. Many algorithms to compute
toric ideals exist and we present a few of them here. 

Let ${\cal A} = \{a_1,a_2,...,a_n\}$ be a subset of ${\integers}^d$.
Consider the map 
\begin{eqnarray}
{\pi}:k[x] \mapsto k[t^{\pm 1}] \\
x_i \mapsto t^{a_i} 
\end{eqnarray}

Recall that the kernel of $\pi$ is the toric ideal of $\cal A$
and we denote it by $I_{\cal A}$. The most basic method to compute
$I_{\cal A}$ would be the elimination method. Though this method is
computationally expensive and not recommended, it serves as a starting
point. Note that every vector $u \in {\integers}^n$ can be written
uniquely as $u = u^+ - u^{-}$ where $u^+$ and $u^-$ are non-negative
and have disjoint support.

\begin{algo} \label{elimalgo}
[Algorithm 4.5, \cite{sturmfelstoric}].
\end{algo} 
\begin{enumerate}
\item Introduce $n+d+1$ variables $t_0,t_1,..,t_d,x_1,x_2,...,x_n$.

\item Consider any elimination order with $\{t_i; i=0, \dots, d\} >
\{x_j; j=1, \dots, n\}$. Compute the reduced Gr\"obner basis $G$ for
the ideal
\[ 
(t_0t_1t_2...t_d - 1, x_1t^{a_1-} - t^{a_1+},....,x_nt^{a_n-} - t^{a_n+}).
\]

\item $G \cap k[x]$ is the reduced Gr\"obner basis for $I_{\cal A}$ with 
respect to the chosen elimination order. 
\end{enumerate}

If the lattice points $a_i$ have only non-negative coordinates, the
variable $t_0$ is unnecessary and we can use the ideal $( x_i- t^{a_i}
: i = 1, \dots, n)$ in the second step of the Algorithm \ref{elimalgo}.

To reduce the number of variables involved in the Gr\"obner basis
 computations, it is better to use an algorithm that operates entirely
 in $k[x_1, \dots, x_n]$.  We now present such an algorithm for
 homogeneous ideals. Observe that all the toric ideals we face in our
 computations in this thesis are homogeneous.

Recall that the {\em saturation of an ideal} $J$ denoted by
$(J:f^{\infty})$ is defined to be
\[
(J:f^{\infty}) = \{g \in k[x]: f^rg \in J \mbox{ for some } r \in
\naturals \}.
\]

Let $ker ({\cal A}) \in Z^n$ denote the integer kernel of the $d
\times n$ matrix with column vectors $a_i$. With any subset $\cal C$ of the
lattice $ker ({\cal A})$ we associate a subideal of $I_{\cal A}$:

\[
J_{\cal C} := (X^{u^+} - X^{u^-} : u \in {\cal C}).
\]

We now describe another algorithm to compute the toric ideal $I_{\cal
A}$.
\begin{algo} \label{stualgo}
[Algorithm 12.3 \cite{sturmfelstoric}].
\end{algo}
\begin{enumerate}
\item Find any lattice basis $L$ for  $ker ({\cal A})$.

\item Let $J_L := (X^{u^+} - X^{u^-}: u \in L)$.

\item Compute a Gr\"obner basis of $(J_L : (x_1x_2 \cdots x_n)^{\infty})$
which is also a Gr\"obner basis of the toric ideal $I_{\cal A}$.
\end{enumerate}

From the computational point of view, computing $(J_L : (x_1x_2 \cdots
x_n)^{\infty})$ is the most demanding step. The algorithms implemented
in CoCoA try to make this step efficient \cite{bigatti1}. For example,
one way to compute $(J_L : (x_1x_2 \cdots x_n)^{\infty})$, would be to
eliminate $t$ from the ideal $H := J_L + (tx_1x_2 \cdots x_n -1)$ but
this destroys the homogeneity of the ideal. It is well-known that
computing with homogeneous ideals have many advantages.  Therefore, it
is better to introduce a variable $u$ whose degree is the sum of the
degrees of the variables $x_i, i = 1, \dots , n$. We then compute the
Gr\"obner basis of the ideal $H := J_L + (x_1x_2 \cdots x_n - u)$
. Then a Gr\"obner basis for $(J_L : (x_1x_2 \cdots x_n)^{\infty})$ is
obtained by simply substituting $u = x_1x_2 \cdots x_n$ in the
Gr\"obner basis of $H$.

Another trick to improve the efficiency of the computation of
saturation ideals is to use the fact
\[
(J_L : (x_1x_2 \cdots x_n)^{\infty}) = 
((\dots ((J_L : x_1^{\infty}) : x_2^{\infty}) \dots ): x_n^{\infty}).
\]

Therefore we can compute the saturations sequentially one variable at
a time. See \cite{bigatti2} for other tricks. We refer the reader to
\cite{sturmfelstoric} for details and proofs of the concepts needed to
develop these algorithms and other algorithms.  We now illustrate
Algorithm \ref{stualgo} by applying it to an example.

Let ${\cal A} = \{(1,1), (2,2), (3,3)\}$. Consider the matrix
whose columns are the vectors of ${\cal A}$
\[
\left [
\begin{array}{ccc}
1  & 2  & 3\\
1 & 2 & 3 
\end{array}
\right ].
\]

Then $ker{\cal A} = \{[-2,1,0], [-3,0,1] \}$. A lattice basis of
$ker{\cal A}$ can be computed using the software Maple, and we get a
basis is
\[
\{ [-1,-1,1], [-2,1,0] \}.
\]

Therefore $J_L = (x_3 - x_1x_2, x_2 - x_1^2)$ and
 \[(J_L :(x_1x_2x_3)^{\infty}) = (x_3 - x_1x_2, x_2 - x_1^2, x_2^2 - x_1x_3)\]
which is also $I_A$. 

Note that many available computer algebra packages including CoCoA can
compute saturation of ideals.

\section{Algorithms to compute Hilbert Poincar\'e series.}
In this section, we will describe a pivot-based algorithm to compute
the Hilbert Poincar\'e series \cite{bigatti2}. Variations of this
algorithm is implemented in CoCoA.  Let $k$ be a field and
$R:=k[x_1,x_2,...,x_r]$ be a graded Noetherian ring. let
$x_1,x_2,...,x_r$ be homogeneous of degrees $k_1,k_2,..,k_r$ (all $>
0$). Let $M$ be a finitely generated $R$-module. Let $H$ be an
additive function on the class of $R$-modules with values in
$\integers$.  Then by the Hilbert-Serre theorem, we have
$$H_{M}(t)=\frac{p(t)}{\Pi_{i=1}^r(1-t^{deg x_i})}.$$
where $p(t) \in {\integers}[t] $. 

Let $I$ be an ideal of $R$, we will denote 
\[
H_{R/I}(t) = \frac{<I>}{\Pi_{i=1}^r(1-t^{deg x_i})}.
\]
Observe that we only need to calculate the numerator $<I>$ since the
denominator is already known.

Let $y$ be a monomial of degree $(d_1,...,d_r)$ called the {\em
pivot}.  The degree of the pivot is $d = \sum_{i=1}^r d_i$. 
Recall the definition of {\em ideal quotients} $(J:f)$ \cite{sturmfelstoric}
\[
(J:f) = \{g \in k[x]: fg \in J \}.
\]

 Consider the
following short exact sequence on graded $R$-modules.

\[
0 \rightarrow R/(I:y) \stackrel{y}{\rightarrow} R/I \rightarrow R/(I,y) 
\rightarrow 0
\]

which yields  (since $H$ is additive)   
\[ 
H_{R/I}(t) = H_{R/(I,y)}(t) + t^d(H_{R/(I:y)})(t). 
\]

This implies  
\begin{eqnarray} \label{recur}
<I> = <I,y> + t^d <I:y>. 
\end{eqnarray}

When $I$ is a homogeneous ideal,
\[
 H_{R/I}(t) = H_{R/\mbox{in}(I)}(t),
\]

where in$(I)$ denotes the ideal of initial terms of $I$
\cite{coxlittleoshea}.

The pivot $y$ is usually chosen to be a monomial that divides a
generator of $I$ so that the total degrees of 
$(I,y)$ and $(I:y)$ are lower than the total degree of $I$. 
The computation proceeds inductively.

We illustrate this algorithm with an example. Let $R=k[x_1,x_2, \dots,
x_n]$ be the polynomial ring. Let $R = {\bigoplus}_{d \in \naturals} R_d$
where each $R_d$ is minimally generated as a $k$-vector by all the
${n+d-1} \choose d$ monomials of degree $d$. Therefore,
\[
H_{R/(0)}(t) = H_R(t) = \sum_{d=0}^{\infty} \mbox{dim}R_dt^d =
\sum_{d=0}^{\infty}{{n+d-1} \choose d}t^d = 1/(1-t)^n.
\]

Therefore we get $<0>=1$. We will use this information to compute 
$H_{R/(I)}(t)$, where $I= (x_1, x_2, \dots, x_n)$.

Let $J= (x_2, \dots, x_n)$. Then, $(J:x_1) = J$. Therefore by Equation
\ref{recur}, we get 
\[
<(J,x_1)> = (1-t^{\mbox{deg}x_1})<J>.  
\]

That is, 
\[ <x_1, x_2, \dots, x_n > = (1- t^{\mbox{deg}x_1}) <x_2, \dots, x_n >.
\]

Now, choosing the pivot $x_2, x_3, \dots, x_n$  subsequently we get
\[
<x_1, x_2, \dots, x_n > = \prod_{i=1, \dots, n} (1-t^{\mbox{deg} x_i})<0>.
\]

Now since $<0> = 1$, we get $<x_1, x_2, \dots, x_n > = \prod_{i=1,
\dots, n} (1-t^{\mbox{deg} x_i})$.

Therefore $H_{R/(x_1, x_2, \dots, x_n)}(t) = 1$.

See \cite{bigatti2} for the effects of choosing different pivots
in the algorithm and also for other algorithms. LattE uses a different
algorithm from the CoCoA algorithms \cite{rudy}.

\newpage
\pagestyle{myheadings} 
\markright{  \rm \normalsize Appendix C \hspace{0.5cm}}
%Appendix------
\chapter{}
\label{appendixc}
{\small
\begin{quote}
And no grown-up will ever understand 
that this is a matter of so much importance! -- Antoine de Saint-Exup\'ery.
% in {\em The Little Prince}.
\end{quote}
}

\section{Constructing natural magic squares.}

%%%% Figure of Franklin squares %%%%%%%%%%%%%%%%%%%%%%%
\begin{figure}[h]
 \begin{center}
     \includegraphics[scale=0.6]{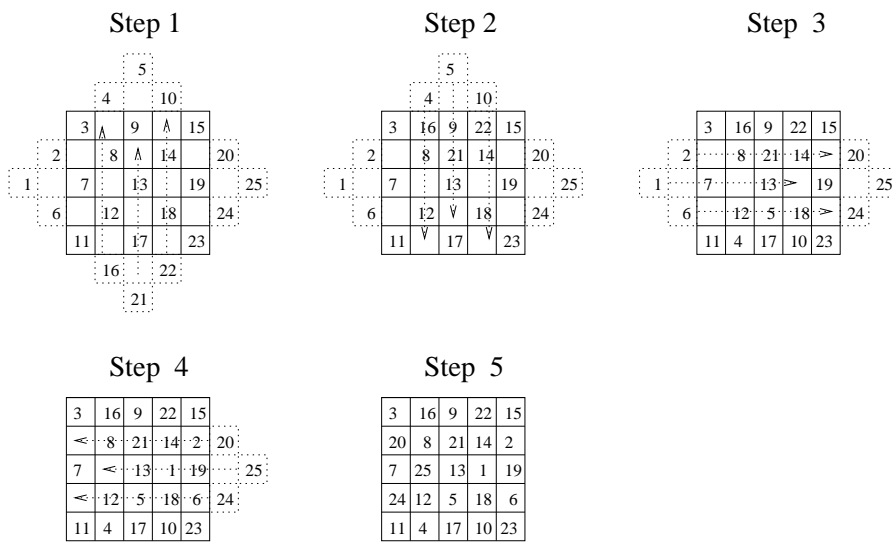}
\caption{Constructing natural magic squares with an odd number of cells \cite{andrews}.} \label{oddmagic}
 \end{center}
 \end{figure}

When the entries of an $n \times n$ magic square are $1,2,3,...,n^2$,
the magic square is called a {\em pure magic square} or a {\em natural
magic square}. Methods for constructing natural magic squares are
known for every order. One of my favorite methods of constructing a
natural magic square with an odd number of cells is as follows: the
numbers 1 to $n^2$ are written consecutively in diagonal columns as
shown in Figure \ref{oddmagic}. The numbers which are outside the
center square are then transferred to the empty cells on the opposite
sides of the latter without changing their order to get a magic
square. This method is said to have been originated by Bachet de
Me\'ziriac (see \cite{andrews}).

\begin{figure}[h]
 \begin{center}
     \includegraphics[scale=0.6]{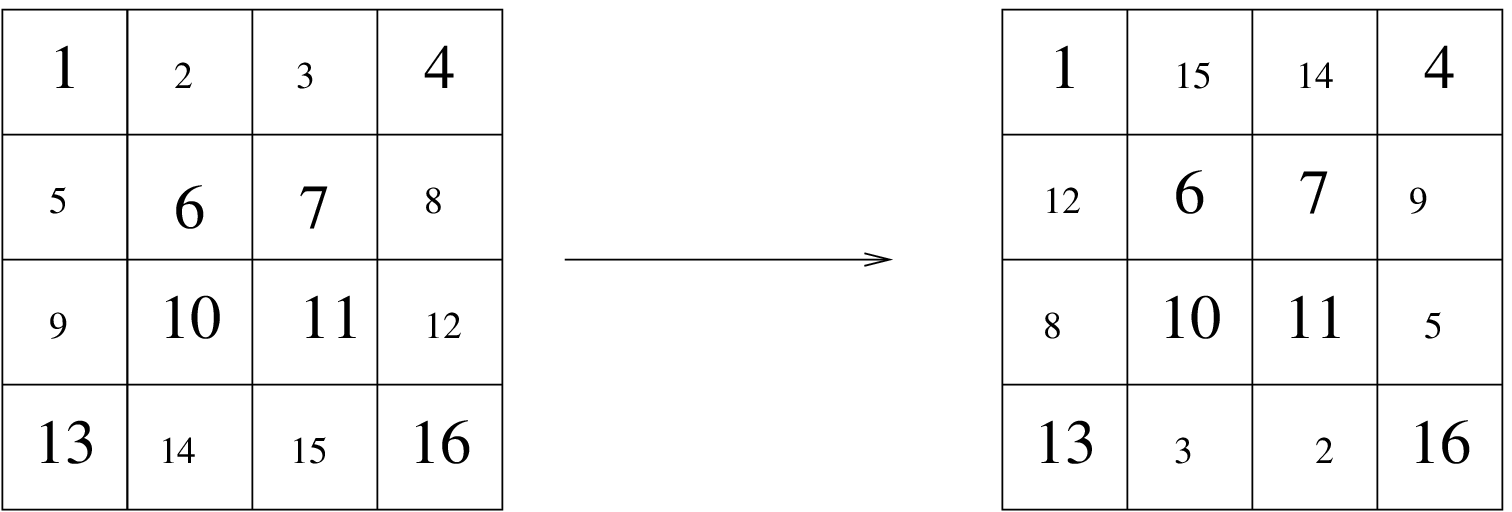}
\caption{Constructing magic squares with an even number of cells \cite{andrews}.} \label{evenmagic}
 \end{center}
 \end{figure}

A pair of numbers in the set $\{1,2,...,n^2\}$ which add to $n^2+1$ is
called {\em complementary}. A method for constructing magic squares
with an even number of cells using complementary pairs of numbers is
as follows: write the number $1 \dots n^2$ consecutively across
rows. All entries except the numbers in the two main diagonals are
then replaced by their complements to get a magic square (see Figure 
\ref{evenmagic}). See \cite{andrews} for details of these methods and 
other methods of constructing natural magic squares. For more recent
developments in the construction of natural magic squares, see 
\cite{beckzaslavsky}. 

%%%% Figure of even magic construction squares %%%%%%%%%%%%%%%%%%%%%%%

\section{Other magic figures.}

\begin{figure}[h]
 \begin{center}
     \includegraphics[scale=0.5]{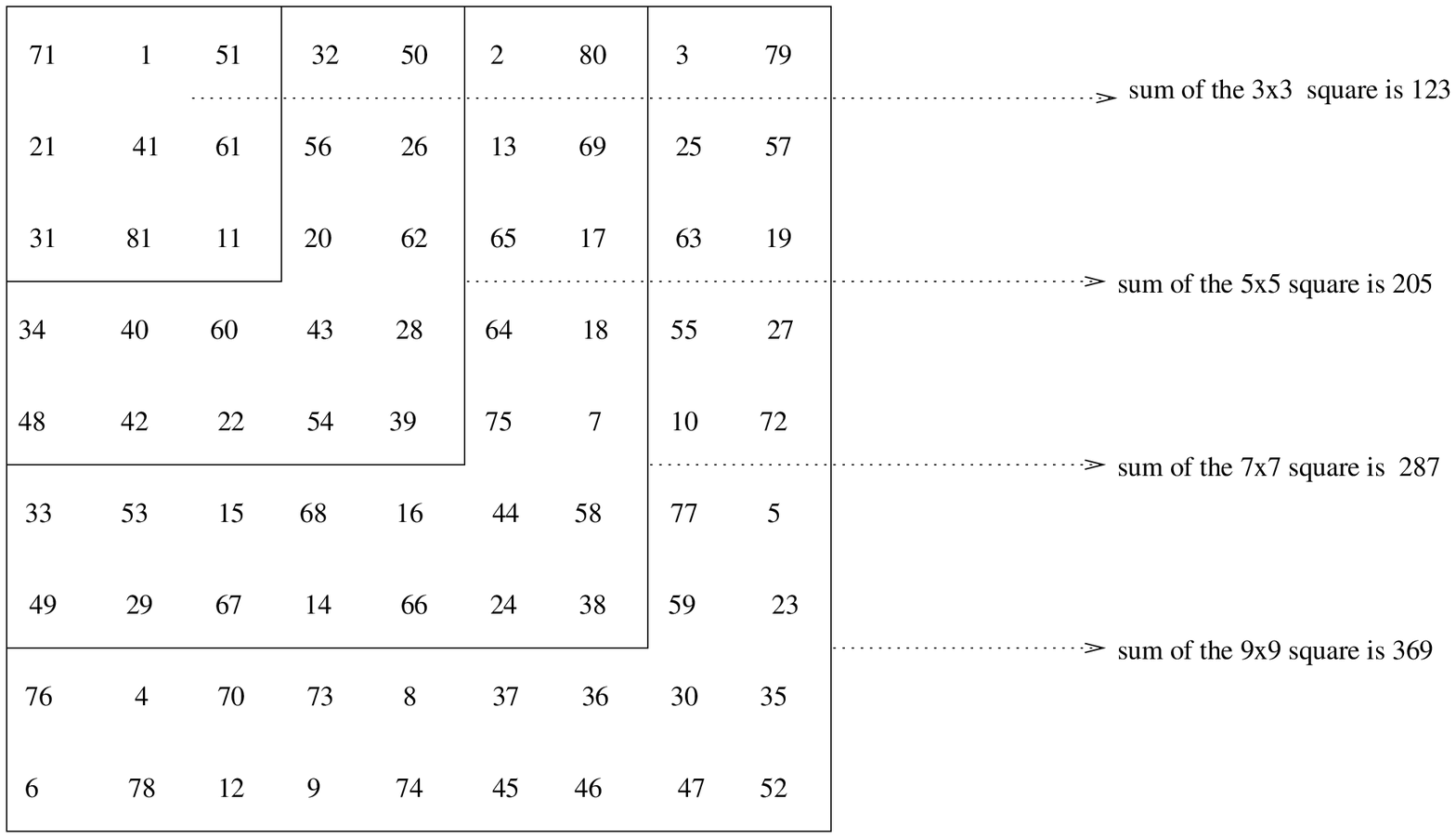}
\caption{Variation of the concentric magic square \cite{andrews}.}
\label{variationconcentric}
 \end{center}
 \end{figure}

A {\em composite magic square} is a magic square composed of a series of small
magic squares and an example is given in Table \ref{composite_square}.
A {\em concentric magic square} is a magic square that remains a magic square 
with borders removed.  An example of a modification of a concentric
magic square devised by Frierson is shown in figure \ref{variationconcentric}.

{\scriptsize
%%%%%%%%%%% TABLE  %%%%%%%%%%%%%%
\begin{table}[hptb]
\begin{center}
\begin{tabular}{|l|l|l||l|l|l||l|l|l|} \hline
71 & 64 & 69 & 8  & 1  & 6  & 53 & 46 & 51 \\ \hline
66 & 68 & 70 & 3  & 5  & 7  & 48 & 50 & 52 \\ \hline
67 & 72 & 65 & 4  & 9  & 2  & 49 & 54 & 47 \\ \hline \hline
26 & 19 & 24 & 44 & 37 & 42 & 62 & 55 & 60   \\  \hline
21 & 23 & 25 & 39 & 41 & 43 & 57 & 59 & 61 \\  \hline
22 & 27 & 20 & 40 & 45 & 38 & 58 & 63 & 56 \\  \hline  \hline
35 & 28 & 33 & 80 & 73 & 78 & 17 & 10 & 15 \\  \hline
30 & 32 & 34 & 75 & 77 & 79 & 12 & 14 & 16 \\  \hline
31 & 36 & 29 & 76 & 81 & 74 & 13 & 18 & 11\\  \hline
\end{tabular}
\end{center}
\caption{A composite magic  square \cite{andrews}.} \label{composite_square}
\end{table}
%%%%%%%%%%%%%%%%%%%%%%%%%%%%%%%%%%%%%%%%%%%%%%%%%%
}

\begin{figure}[h]
 \begin{center}
     \includegraphics[scale=0.5]{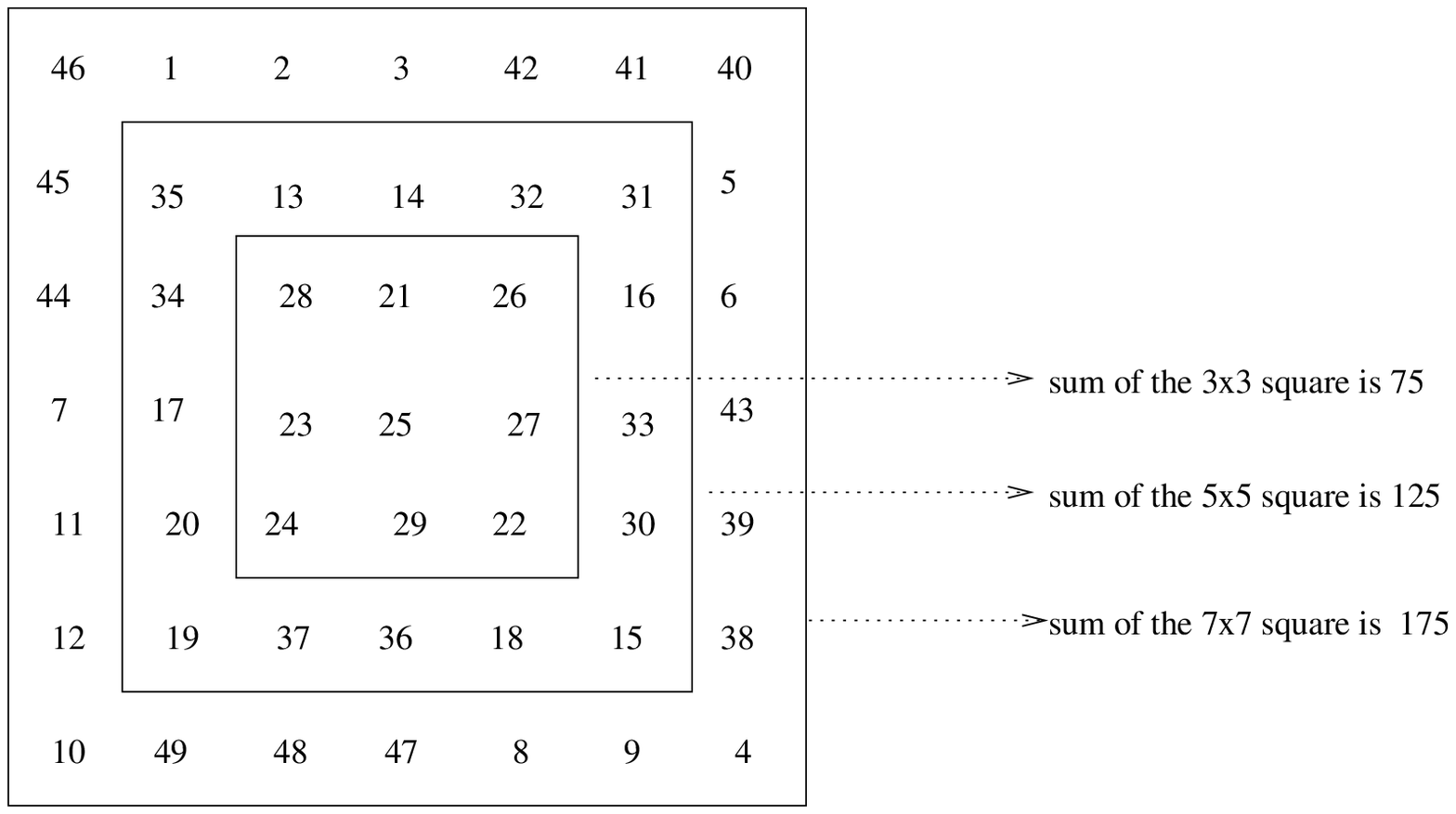}
\caption{A concentric magic square \cite{andrews}.} \label{concentric}
 \end{center}
 \end{figure}

\begin{figure}[h]
 \begin{center}
     \includegraphics[scale=0.5]{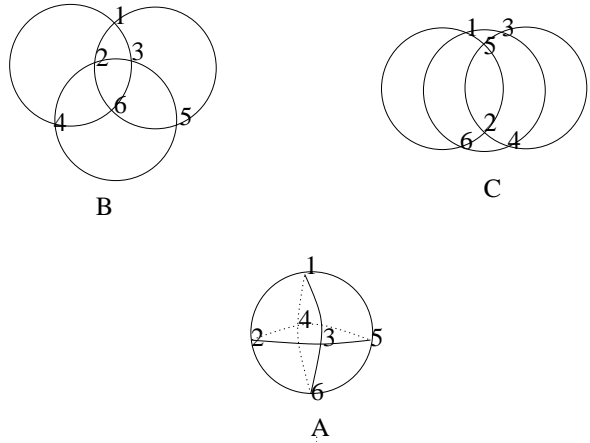}
\caption{Magic circles and magic sphere of dice \cite{andrews}.} \label{dice}
 \end{center}
 \end{figure}

A set of $n$ {\em magic circles} is a numbering of the intersections
of the $n$ circles such that the sum over all intersections is the
same constant for all circles. Consider the example of a die. It is
commonly known that the opposite faces of a die contain complementary
numbers that always add up to 7. Consequently any band of four numbers
encircling a die gives a summation of 14 (see Figure \ref{dice} A).
These bands form magic circles (see Figures \ref{dice} B or \ref{dice}
C). A {\em magic sphere} is a sphere that contains magic circles.  The
sphere in Figure \ref{dice} A is a magic sphere.

\begin{figure}[h]
 \begin{center}
     \includegraphics[scale=0.4]{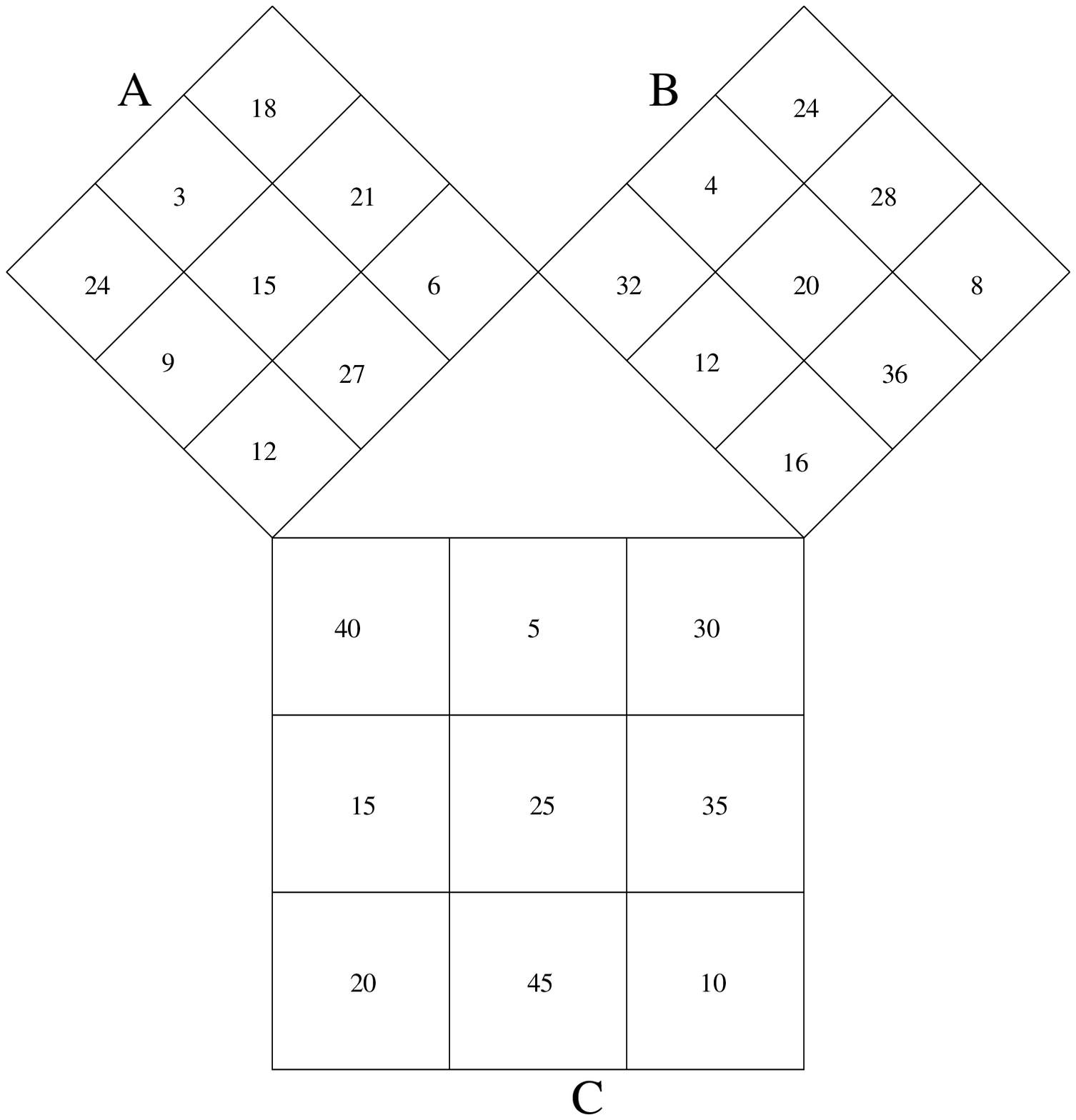}
\caption{A Magic Triangle.} \label{triangle}
 \end{center}
 \end{figure}

A {\em magic triangle} is composed of three magic squares A, B, and C
such that the square of any cell in C is equal to the sum of the
squares of the corresponding cells in A and B. In other words the
corresponding entries of the magic squares always form a Pythagorean
triple ($c^2 = a^2 + b^2$). C is called the hypotenuse, and A and B
are called the legs of the magic triangle. An example is given in
Figure \ref{triangle}.

\begin{figure}[h]
 \begin{center}
     \includegraphics[scale=0.4]{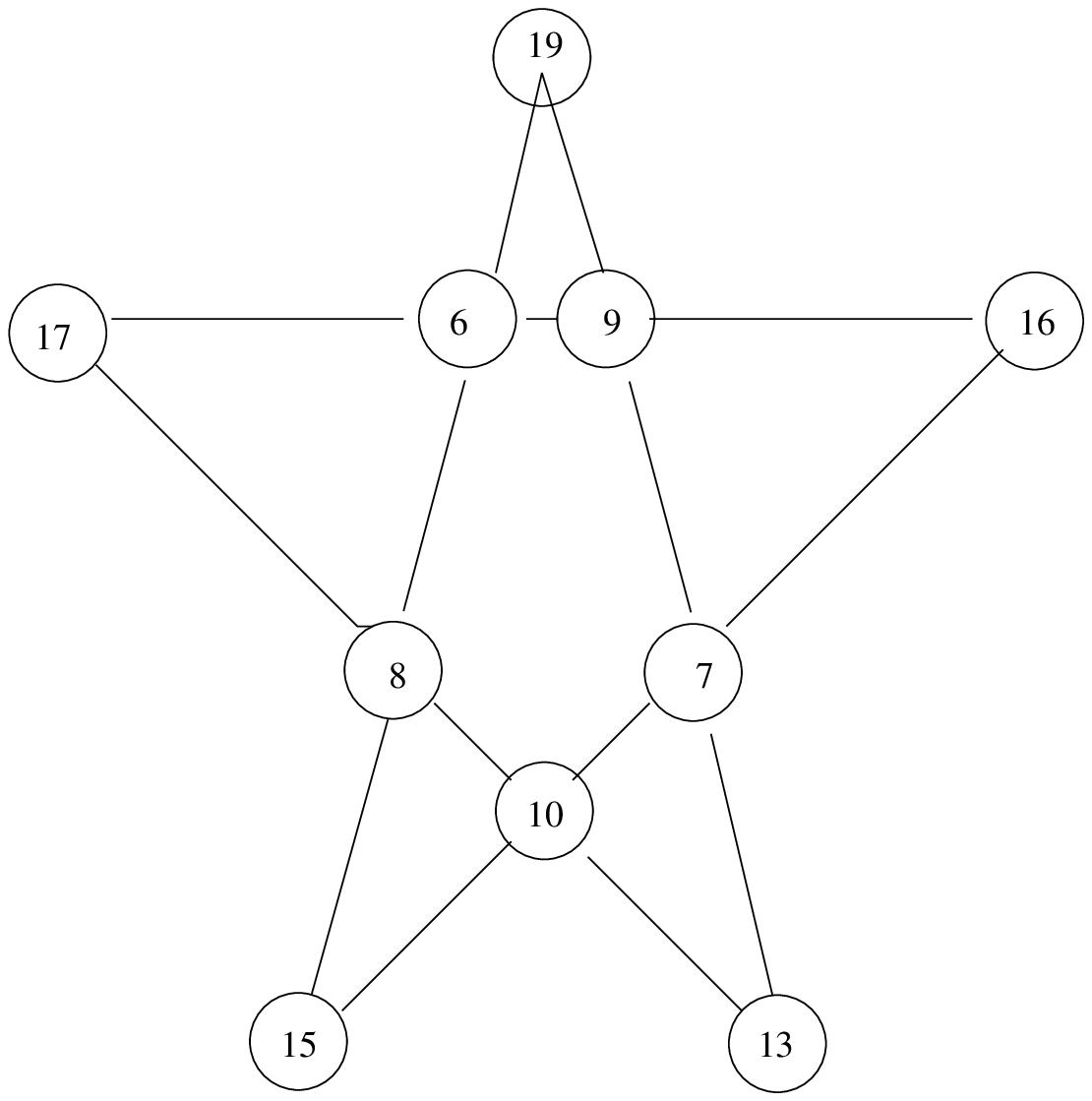}
\caption{A Magic Star.} \label{star}
 \end{center}
 \end{figure}

A {\em magic star} is a numbering of the intersections of a set of lines
that form a star such that the sum over every intersection is the same for
each line (see Figure \ref{star} for an example). 

A {\em magic carpet} is a magic square in which a limited range of
digits is used several times. 

{\scriptsize
%%%%%%%%%%% TABLE  %%%%%%%%%%%%%%
\begin{table}[hptb]

\begin{center}
\begin{tabular}{|l|l|l|l|} \hline
   0 & 0 & 1 & 1 \\ \hline
   1 & 1 & 0 & 0 \\   \hline
   0 & 0 & 1 & 1 \\  \hline
   1 & 1 & 0 & 0  \\ \hline
\end{tabular}
\end{center}
\caption{A magic carpet.} \label{carpet}
\end{table}
%%%%%%%%%%%%%%%%%%%%%%%%%%%%%%%%%%%%%%%%%%%%%%%%%%
}

{\scriptsize
%%%%%%%%%%% TABLE  %%%%%%%%%%%%%%
\begin{table}[hptb]
\begin{center}
\begin{tabular}{|l|l|l|l|l|l|l|l| } \hline
   0 & 0 & 1 & 1 & 0 & 0 & 1 & 1 \\ \hline
   1 & 1 & 0 & 0 & 1 & 1 & 0 & 0 \\   \hline
   0 & 0 & 1 & 1 & 0 & 0 & 1 & 1 \\  \hline
   1 & 1 & 0 & 0 & 1 & 1 & 0 & 0 \\ \hline
   0 & 0 & 1 & 1 & 0 & 0 & 1 & 1 \\ \hline
   1 & 1 & 0 & 0 & 1 & 1 & 0 & 0 \\   \hline
   0 & 0 & 1 & 1 & 0 & 0 & 1 & 1 \\  \hline
   1 & 1 & 0 & 0 & 1 & 1 & 0 & 0 \\ \hline
\end{tabular}
\end{center}
\caption{A magic carpet made from the magic carpet in Table
\ref{carpet}.} \label{bigcarpet}
\end{table}
%%%%%%%%%%%%%%%%%%%%%%%%%%%%%%%%%%%%%%%%%%%%%%%%%%
}

\begin{figure}[h]
 \begin{center}
     \includegraphics[scale=0.5]{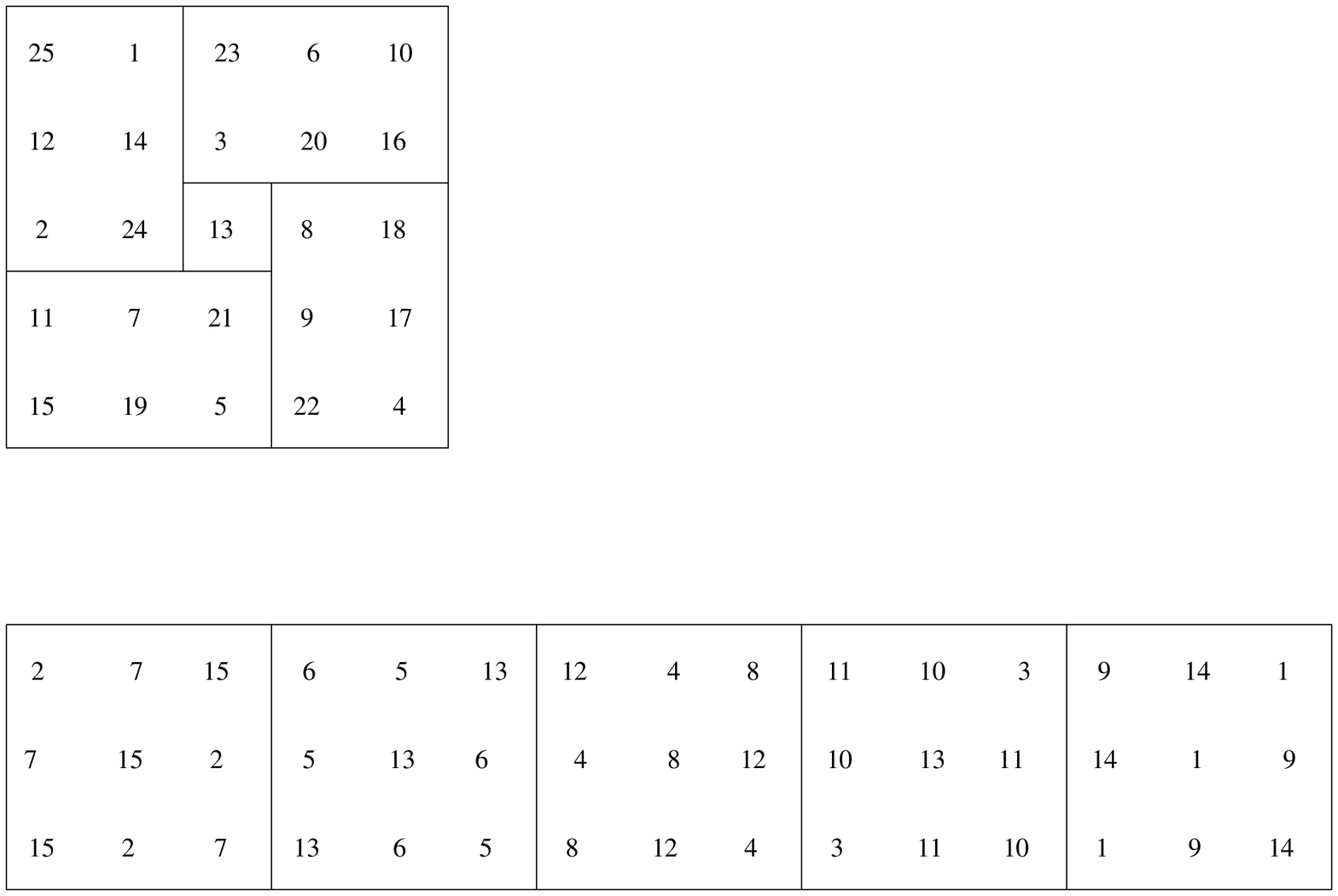}
\caption{Ornate magic squares \cite{andrews}.} \label{squaresandrectangle}
 \end{center}
 \end{figure}

A {\em magic rectangle} is an $m \times n$ matrix such that all its
rows add to a prescribed common sum, and all its columns add to the
another prescribed sum. Therefore, a magic rectangle has two magic
sums: A {\em magic column sum} and a {\em magic row sum}.  Figure
\ref{squaresandrectangle} shows embeddings of magic squares in magic
rectangles and vice versa and such patterns are called {\em ornate
magic squares} (\cite{andrews}).

 See \cite{andrews} for more examples of magic figures. Needless to
say, we sure can construct and enumerate most of these magic figures
with our methods.

%%%%%%%%%%%%%%%%%    BIBLIOGRAPHY    %%%%%%%%%%%%%%%%%%%%%%%%%%%

\newpage
\pagestyle{myheadings} 
\markright{  \rm \normalsize BIBLIOGRAPHY \hspace{0.5cm}}
\bibliographystyle{plain}
\bibliography{phd}

\begin{thebibliography}{99}
\thispagestyle{myheadings}
\addcontentsline{toc}{chapter}{\bf Bibliography}

\bibitem{alvin}{Alvis, D. and Kinyon, M.} {\em Birkhoff's theorem for
Panstochastic matrices}, Amer. Math. Monthly,  (2001),  Vol 108,  no.1,
28-37.

\bibitem{adh} Ahmed, M., De Loera, J., and Hemmecke, R., {\em Polyhedral cones
of magic cubes and squares}, New Directions in Computational Geometry,
The Goodman-Pollack Festschrift volume, Aronov et al., eds.,
Springer-Verlag, (2003), 25--41.

\bibitem{ma} Ahmed, M., {\em How many squares are there,
Mr. Franklin?: Constructing and Enumerating Franklin Squares},
   Amer. Math. Monthly, Vol. 111, 2004, 394--410.

\bibitem{maya} \underline{\makebox[.5in]{}}, {\em Magic graphs and the
 faces of the Birkhoff polytope}, arXiv:math.CO/0405181

\bibitem{alantarsi} Alon, N. and  Tarsi, M., {\em  A note on graph colorings and graph polynomials},  J. Combin. Theory Ser. B 70 (1997), no. 1, 197--201.

\bibitem{anandgupta} Anand, H., Dumir, V.C., and Gupta, H., {\em A 
combinatorial distribution problem}, Duke Math. J. 33, (1966),
757-769.


\bibitem{andrews} Andrews, W. S., {\em Magic Squares and Cubes},
 2nd. ed., Dover, New York, 1960.

\bibitem{atiyah} Atiyah, M.F.,  and  Macdonald, I.G.,  {\em Introduction to
Commutative Algebra}, Addison-Wesley, Reading, MA, 1969.

\bibitem{ball} Ball, W.W.R. and Coxeter, H.S.M, {\em Mathematical 
Recreations and Essays}, 13th edition, Dover Publications, Inc, New
York, 1987.

\bibitem{beckthesis}  Beck, M., {\em The arithmetic of rational polytopes},
Dissertation, Temple University (2000).

\bibitem{beck}  Beck, M. and Pixton, D., {\em The Ehrhart polynomial of the
Birkhoff polytope}, to appear in  Discrete and computational geometry, 
Springer-Verlag, New York.

\bibitem{becketal} Beck, M., Cohen, M., Cuomo, J., and Gribelyuk, P., 
{\em The number of magic squares, cubes and hypercubes}, Amer. Math. Monthly,
110, no.8, (2003), 707-717. 

\bibitem{beckzaslavsky} Beck, M. and  Zaslavsky, T., {\em Inside-Out 
Polytopes}, arXiv:math.CO/0309330.

\bibitem{bigatti1} Bigatti, A.M., La Scala, R., and Robbiano, L.,
{\em Computing toric ideals}, J. Symbolic Computation, 27, (1999),
351-365.

\bibitem{bigatti2}{Bigatti, A.M}, {\em Computation of Hilbert-Poincar\'e Series}, 
J. Pure Appl. Algebra, 119/3, (1997), 237--253. 


\bibitem{billerasarang} Billera, L.J. and  Sarangarajan, A., {\em 
 The combinatorics of permutation polytopes}, Formal power series and
 algebraic combinatorics (New Brunswick, NJ, 1994), 1--23, DIMACS
 Ser. Discrete Math. Theoret. Comput. Sci., 24, Amer. Math. Soc.,
 Providence, RI, 1996.

\bibitem{bona} {Bona, M.}, {\em Sur l'enumeration des cubes magiques},
C. R. Acad. Sci. Paris Ser. I Math., 316, (1993), no.7, 633-636.

\bibitem{brualdigibson} Brualdi, A. R. and Gibson, P., {\em Convex polyhedra
of doubly stochastic matrices: I, II, III}, Journal of combinatorial Theory,
A22, (1977), 467-477.

\bibitem{brualdi} Brualdi, A.R., {\em Introductory combinatorics},
3 rd ed., Prentice hall, New Jersey, 1999. 

\bibitem{brunskoch} {Bruns, W. and Koch, R.}, {\em {\tt NORMALIZ},  Computing
normalizations of affine semigroups}, Available via anonymous ftp from
ftp//ftp.mathematik.uni-onabrueck.de/pub/osm/kommalg/software/


\bibitem{cocoa} Capani, A., Niesi, G., and Robbiano, L., {\em {\tt CoCoA}, A
System for Doing Computations in Commutative Algebra}, available via
anonymous ftp from {\tt cocoa.dima.unige.it} (2000).

\bibitem{carlitz} Carlitz, L., {\em Enumeration of symmetric arrays},
Duke Math. J., 33,  (1966), 771-782.  

\bibitem{chanrobbins} Chan, S. C. and Robbins, D. P., {\em On the volume of
the polytope of doubly stochastic matrices}, Experiment. Math. 8
(1999), no.3, 291-300.

\bibitem{contejean}{Contejean,  E. and Devie, H.}, {\em Resolution de
systemes lineaires d'equations diophantienes}, C. R. Acad. Sci. Paris
 S\'er. I Math., 313, (1991), no. 2, 115--120. 

\bibitem{coxlittleoshea}{Cox, D., Little, J., and O'Shea, D.}, {\em
Ideals,  varieties, and Algorithms}, Springer Verlag,
Undergraduate Text, 2nd Edition, 1997.

\bibitem{clo2} \underline{\makebox[.5in]{}}, {\em Using
Algebraic Geometry}, Springer-Verlag, New York, 1998.

\bibitem{rudy} De Loera, J., Hemmecke, R., Tauzer, J., and Yoshida,
R., {\em Effective lattice point counting in rational convex
polytopes}, to appear in the Journal of Symbolic Computation.

\bibitem{deloerasturmfels}{De Loera, J.A. and Sturmfels, B.}, {\em 
Algebraic unimodular counting}, Algebraic and geometric methods in
discrete optimization. Math. Program. 96 (2003), no. 2, Ser. B,
183--203.

\bibitem{dummit} Dummit, D. S. and Foote, R. M., {\em Abstract Algebra},
Prentice Hall, New Jersey, 1991. 

\bibitem{Ehrhart2} {Ehrhart, E.}, {\em Sur un probl\'eme de
g\'eom\'etrie diophantienne lin\'eaire II}, J. Reine Angew. Math., 227,
(1967),  25-49.

\bibitem{Ehrhart1} \underline{\makebox[.5in]{}}, {\em Figures magiques et methode des polyedres}, J. Reine Angew. Math., 299/300, (1978), 51-63.

\bibitem{Ehrhart0} \underline{\makebox[.5in]{}}, {\em Sur les carr\'es magiques}, 
C. R. Acad. Sci., Paris, 227 A,  (1973),  575-577.

\bibitem{gardner} Gardner, M., {\em Martin Gardner's New mathematical
Diversions from Scientific American}, Simon and Schuster, New York,
(1966), 162-172.

\bibitem{gilespulley} Giles, F.R. and Pulleyblank, W.R.,  Total dual
integrality and integer polyhedra, {\em Linear Algebra Appl.}, 25,
(1979), 191-196.

\bibitem{gupta} Gupta, H., {\em Enumeration of symmetric matrices},
Duke Math. J., 35, (1968), 653-659.

\bibitem{hagstorm} {Hagstorm R.}, {\em Superlatively Regular $8 \times 8$
Magic Squares}, Personal communication (2004). 

\bibitem{halleck} {Halleck, E.Q.}, {\em Magic squares subclasses as
linear Diophantine systems},  Ph.D. dissertation,  Univ. of California
San Diego,  (2000),  187 pages.

\bibitem{raymond} Hemmecke, R.,  {\em On the computation of Hilbert bases of
cones}, in Proceedings of First International Congress of Mathematical
Software, A. M. Cohen, X.S. Gao, and N. Takayama, eds., Beijing,
(2002); software implementation 4ti2 is available from
http://www.4ti2.de.

\bibitem{henkweismantel}{Henk, M. and Weismantel, R.}, {\em On Hilbert bases of
polyhedral cones}, Results in Mathematics, 32, (1997), 298-303.

\bibitem{hilton} Hilton, H., {\em An introduction to the theory of
groups of finite order}, Oxford, 1908.
 
\bibitem{konig} K\"onig, D., {\em Theory of finite and infinite
graphs}, Birkh\"auser Boston, 1990.

\bibitem{lovasz} Lo\'vasz, L. and Plummer, M. D., {\em Matching
Theory}, North-Holland, Amsterdam, 1986.

\bibitem{macmahon} MacMahon, P.A., {\em Combinatorial Analysis},  
Chelsea, 1960.

\bibitem{pak} Pak, I.,  {\em On the number of faces of certain transportation polytopes}, European J. Combinatorics, vol. 21 (2000), 689-694.

\bibitem{pasles} Pasles, P. C., {\em The lost squares of Dr. Franklin:
Ben Franklin's missing squares and the secret of the magic circle},
Amer. Math. Monthly, 108, (2001), 489-511.

\bibitem{pasles2}  \underline{\makebox[.5in]{}}, {\em Franklin's other 8-square},  J. Recreational Math.,  31, (2003), 161-166.

\bibitem{patel} L. D. Patel, {\em The secret of Franklin's $8 \times 8$
magic square}, J.Recreational Math., 23, (1991), 175-182.

\bibitem{zen} Pickover, C.A., {\em The Zen of magic squares, circles, and
stars}, Princeton University Press, New Jersey, 2002.

\bibitem{pottier} {Pottier, L.}, {\em Bornes et algorithme de calcul
des g\'en\'erateurs des solutions de syst\'emes diophantiens
lin\'eaires}, C. R. Acad. Sci. Paris,  311,  (1990), no. 12,
813-816.

\bibitem{pottier2} {Pottier, L.}, {\em Minimal solutions of linear
Diophantine systems: bounds and algorithms}, in Rewriting techniques
and applications (Como, 1991), 162--173, Lecture Notes in
Comput. Sci., 488, Springer, Berlin, 1991.

\bibitem{schrijver} {Schrijver, A.}, {\em Theory of Linear and
Integer Programming}, Wiley-Interscience,  1986.

\bibitem{schubert} {Schubert, H.}, {\em Mathematical Essays and Recreations}, The open Court Publishing Co., 1899.

\bibitem{stanley}{Stanley, R.P.}, {\em Enumerative Combinatorics}, Volume
I, Cambridge, 1997.

\bibitem{stanley2} \underline{\makebox[.5in]{}}, {\em Combinatorics and commutative algebra}, Progress in Mathematics, 41,  Birkha\"user Boston, MA,
1983.

\bibitem{stanley3}) \underline{\makebox[.5in]{}}, {\em Linear Homogeneous Diophantine Equations and Magic Labelings Of Graphs}, Duke Mathematical Journal,
Vol. 40, September 1973, 607-632.

\bibitem{stanley4} \underline{\makebox[.5in]{}}, {\em Magic Labelings of Graphs, Symmetric Magic Squares, Systems of Parameters and Cohen-Macaulay
Rings}, Duke Mathematical Journal, Vol. 43, No.3, September 1976, 511-531.

\bibitem{stewart} Stewart, B. M., {\em Magic graphs}, Canad. J. Math., vol. 18, (1966), 1031-1059.

\bibitem{stewartcomplete} \underline{\makebox[.5in]{}}, {\em Supermagic complete graphs},  Canad. J. Math., vol. 19, (1967), 427-438.

\bibitem{sturmfels} {Sturmfels, B.,} {\em Algorithms in invariant
theory}, Springer-Verlag, Vienna, 1993.

\bibitem{sturmfelstoric} \underline{\makebox[.5in]{}}, {\em Gr\"obner Bases and Convex Polytopes}, University Lecture Series, no. 8, American Mathematical Society, Providence, 1996.

\bibitem{thiery} Thi\'ery, N. M., {\em Algebraic invariants of graphs;
a study based on computer exploration}, SIGSAM Bulletin (2000), 9-20.

\bibitem{wallis} Wallis, D., {\em Magic Graphs}, Birkh\"auser Boston, 2001.

\bibitem{lintwilson} Van Lint, J.H. and Wilson, R.M., {\em A course in
Combinatorics}, 2 nd edition, Cambridge University Press, Cambridge,
2001.


\bibitem{vergne}{ Vergne, M. and Baldoni-Silva, W.}, {\em Residues
formulae for volumes and Ehrhart polynomials of convex polytopes}, 
manuscript 81 pages available at math.ArXiv, CO/0103097.


\end{thebibliography}

\end{document}